\documentclass[a4paper,fleqn]{my-cas-sc}

\usepackage[numbers]{natbib}
 
\usepackage{enumerate}
\usepackage{graphicx} 
\usepackage{multirow,cases,subfigure}
\usepackage{longtable}
\usepackage{gensymb}
\usepackage{amssymb}
\usepackage{lipsum}
\usepackage{stfloats}
\usepackage{float}
\usepackage{subfigure}
\usepackage{caption}
\restylefloat{table}
\usepackage{amssymb,amsfonts,amsmath,mathtools,amsthm,enumerate}

\newtheorem{theorem}{Theorem}
\newtheorem{proposition}[theorem]{Proposition}

\def\tsc#1{\csdef{#1}{\textsc{\lowercase{#1}}\xspace}}
\tsc{WGM}
\tsc{QE}
\tsc{EP}
\tsc{PMS}
\tsc{BEC}
\tsc{DE}

\ExplSyntaxOn
\keys_set:nn { stm / mktitle } { nologo }
\ExplSyntaxOff

\begin{document}

\let\WriteBookmarks\relax
\def\floatpagepagefraction{1}
\def\textpagefraction{.001}
\shorttitle{}
\shortauthors{M. Charles et~al.}

\title[mode=title]{Mathematical model to assess the impact of contact rate 
and environment factor on transmission dynamics of rabies in humans and dogs}   


\author[1,2]{Mfano Charles}[orcid=0009-0008-8294-940X]
\ead{mfanoc@nm-aist.ac.tz}
\cormark[1]
\cortext[cor1]{Corresponding author: mfanoc@nm-aist.ac.tz}
\credit{Conceptualization, 
Methodology, 
Writing -- original draft, 
Writing -- review \& editing}


\author[1]{Verdiana G. Masanja}[orcid=0000-0003-4844-9133]
\ead{verdiana.masanja@nm-aist.ac.tz}
\credit{Conceptualization, 
Methodology, 
Writing -- original draft, 
Writing -- review \& editing, 
and made insightful suggestions}


\author[3]{Delfim F. M. Torres}[orcid=0000-0001-8641-2505]
\ead{delfim@ua.pt}
\credit{Conceptualization, 
Methodology, 
Writing -- original draft, 
Writing -- review \& editing, 
and made insightful suggestions}


\author[4]{Sayoki G. Mfinanga}[orcid=]
\ead{gsmfinanga@yahoo.com}
\credit{Conceptualization, 
Methodology, 
Writing -- original draft, 
Writing -- review \& editing, 
and made insightful suggestions}

	
\author[1]{G.~A. Lyakurwa}[orcid=]
\ead{geminpeter.lyakurwa@nm-aist.ac.tz}
\credit{Conceptualization, 
Methodology, 
Writing -- original draft, 
Writing -- review \& editing, 
and made insightful suggestions}


\address[1]{School of Computational and Communication Science and Engineering, 
The Nelson Mandela African Institution of Science and Technology (NM-AIST), 
P.O. BOX 447 Arusha, Tanzania}

\address[2]{Department of ICT and Mathematics, 
College of Business Education (CBE), P.O. BOX 1968 Dar es Salaam, Tanzania}

\address[3]{Center for Research and Development in Mathematics and Applications (CIDMA),
Department of Mathematics, University of Aveiro, 3810-193 Aveiro, Portugal}

\address[4]{NIMR Chief Research Scientist Fellow}

	
\begin{abstract}
\noindent This paper presents a mathematical model to understand 
how rabies spreads among humans, free-range, and domestic dogs. 
By analyzing the model, we discovered that there are equilibrium points 
representing both disease-free and endemic states. We calculated the basic 
reproduction number, $\mathcal{R}_{0}$, using the next generation matrix method. 
When $\mathcal{R}_{0}<1$, the disease-free equilibrium is globally stable, 
whereas when $\mathcal{R}_{0}>1$, the endemic equilibrium is globally stable. 
To identify the most influential parameters in disease transmission, we used 
the normalized forward sensitivity index. Our simulations revealed that the 
contact rates between the infectious agent and humans, free-range dogs, 
and domestic dogs have the most significant impact on rabies transmission. 
The study also examines how periodic changes in transmission rates affect 
the disease dynamics, emphasizing the importance of transmission frequency 
and amplitude on the patterns observed in rabies spread.
Therefore, the study proposes that to mitigate the factors most strongly 
linked to disease sensitivity, effective disease control measures should 
primarily prioritize on reducing the population of both free-range 
and domestic dogs in open environments.
\end{abstract}



\maketitle


\section{Introduction}

Rabies is a viral disease that affects mammals, including humans, caused by the rabies virus 
($\textit{Rabies}\; \text{lyssavirus}$) that travels from the site of infection to the brain, 
causing inflammation and damage to the nervous system \cite{de2022importance,kumar2023canine}. 
Although dogs are the primary source of the virus, causing more than 99\% of human rabies 
infections worldwide, other animals such as bats, raccoons, skunks, and foxes can also carry 
and transmit the virus through bites or scratches 
\cite{de2022importance,kumar2023canine,slathia2023rabies}. 
The transmission dynamics of rabies are influenced 
by environmental factors such as changes in habitat, 
land use patterns, and wildlife populations 
which create varied and frequent interactions 
between infected and susceptible individuals 
\cite{mcmahon2018ecosystem}. 

Symptoms of rabies include fever, headaches, fatigue, restlessness, anxiety, hallucinations, 
hydrophobia (fear of water), difficulty swallowing, and paralysis that manifest between 20 days 
and 3 months after exposure, but may vary from 1 week to 1 year after exposure, depending 
on the location of entry of the virus and the viral load. In rare cases, the incubation period 
can last up to 7 years, and if left untreated and without appropriate medical care such 
as vaccination, the disease progresses into a coma state and ultimately leads to death
\cite{nigg2009overview,johnson2010immune,hailemichael2022effect}.

Rabies is responsible for 60,000 human deaths every year globally \cite{bilal2021rabies}, 
and the effective management of rabies in low- and middle-income countries (LMICs) 
in Asia and Africa is often hindered by the lack of timely and accurate information 
on rabies cases in both humans and animals \cite{mbilo2021dog,tian2018transmission}. 
The actual number of fatalities caused by rabies virus (RABV) infections in LMICs 
is believed to be underestimated, and the dynamics of rabid dog populations are 
poorly understood \cite{tian2018transmission}. For example, a study conducted by 
\cite{sambo2013burden} revealed that the mortality rate attributed to human rabies 
in the United Republic of Tanzania was significantly higher than what has been officially 
reported. By analyzing the active surveillance data on bite incidence, the researchers 
estimated an annual mortality rate of 1499 deaths, with a confidence interval spanning 
from 891 to 2238 deaths. This indicates an annual incidence of 4.9 deaths per 100,000 
population, ranging from 2.9 to 7.2 deaths per 100,000.

The understanding and control of contagious diseases, such as rabies, have been 
significantly enhanced through the application of mathematical models. This analytical 
tool has proven to be instrumental in the prediction and analysis of various phenomena, 
enabling medical experts to structure their approach towards disease management. 
Researchers in different countries have developed numerous mathematical models 
to study the transmission dynamics of rabies in dog populations
\cite{abdulmajid2021analysis,amoako2021rabies,abrahamian2022rhabdovirus,%
ruan2017modeling,ayoade2019saturated,chapwanya2022environment}, 
as well as interactions between dog and human populations 
\cite{kadowaki2018risk,ayoade2023modeling,tulu2017impact}, 
and even among dog, human, and other animal populations
\cite{abdulmajid2021analysis,pantha2021modeling,ega2015sensitivity}. 
These studies have identified several factors that influence the dynamics 
of rabies within their respective countries, and have addressed various strategies 
to control the disease. The results of these studies can help guide public 
health officials and medical professionals in developing effective 
measures to prevent and manage rabies outbreaks.

However, most of these studies have not given enough attention to the impact 
of environmental factors such as urbanization, deforestation, and encroachment 
into wildlife habitats. These factors can result in increased interactions between 
humans, domestic animals, and wildlife, which together make up a population 
of approximately 900 million worldwide. This, in turn, increases the risk 
of rabies transmission \cite{hailemichael2022effect,esposito2023impact}. 
The movement of infected animals, whether domestic or wild, also contributes 
to the spread of the virus. Infected animals can travel long distances, 
introducing the virus to new areas or re-establishing it in regions 
where it was previously under control. 
Studies by \cite{hampson2019potential,rulli2021land} 
have documented how these factors can result in increased contact between 
humans and domestic animals, further facilitating the spread of the virus.

This paper is structured as follows. Section~\ref{sec:2} 
outlines a deterministic mathematical model describing 
the dynamics of rabies. Section~\ref{sec:3} focuses 
on analyzing the proposed model. Real-world data is then used 
in Section~\ref{sec:model:fit:parm:est} for
model fitting and parameter estimation,
while Section~\ref{sec:4} covers 
the implementation of numerical simulations. 
Sections~\ref{sec:5} and \ref{sec:6} are dedicated 
to the discussion and conclusion of the study, respectively.


\section{Model Formulation and Model Analysis}
\label{sec:2}

In this study, we formulate a mathematical model using ordinary differential 
equations by incorporating the interactions between susceptible individuals, 
infectious dogs, and the environment 
\cite{amoako2021rabies,chapwanya2022environment,peng2022transmission,allen2008basic}. 
The model considers the influence of contact rate which represents the frequency 
of encounters between susceptible individuals and infectious dogs.
In particular,  using mathematical simulations, we examine how rabies spreads 
and persists under various contact rate and environmental impact scenarios, 
which helps us better comprehend and predict how rabies spreads. 

	
\subsection{Model formulation}

We divide the  model  into three settings: the human population,  
the dog population, and rabies viruses in the environment.  
The dog population  is divided into two subgroups: free-range dogs,  
defined as  dogs in a public area and not under direct human control, 
such as stray dogs, street dogs, feral dogs, village dogs  and wild dogs, 
as  well as  domestic dogs, defined as dogs who are under human control. 
	
The free-range  dog  population, defined by $N_{F}$, was divided into three 
compartments. First compartment were susceptible free-range dogs denoted 
as $S_{F}$, that were group of free-range dogs free of infection but could 
get infected after adequate contact with either infected free-range dogs 
$\left(I_{F}\right)$, infected domestic dogs $\left(I_{D}\right)$ 
or the  environment containing rabies  virus $\left(M\right)$. The 
second compartment was exposed free-range dogs $\left(E_{F}\right)$,    
which included all the free-range dogs who had contracted the disease 
but could not transmit the disease  and  had no symptoms of rabies infections. 
The third compartment were infected free-range dogs $\left(I_{F}\right)$   
made up of dogs who had contracted the disease with fully developed 
rabies symptoms and were infectious.
	
The domestic dog  population, defined by $N_{D}$, was divided into four  compartments. 
First compartment  is susceptible domestic dogs, denoted by $S_{D}$,  
who were not infected but could get infected after  adequate contact with either 
infected free-range dogs, infected  domestic dogs or environment that 
contains rabies virus. The second compartment were exposed domestic dogs 
$\left(E_{D}\right)$, who had contracted the disease but could not transmit 
the disease and did not have disease symptoms. The third compartment is infectious 
domestic dogs $\left(I_{D}\right)$, who had contracted the disease and were infectious. 
The fourth compartment is recovered domestic dogs, denoted by $R_{D}$,  
who get post-exposure prophylaxis after contract with rabies infectious 
agent and develop temporary immunity.

The  human population, denoted by $N_{H}$, was divided into four compartments.
First compartment were humans who had not acquired the rabies infection 
but could get if they adequately contacted with infectious free-range dogs, 
domestic dogs or contaminated environment. The second compartment  were   
exposed human $\left(E_{H}\right)$, who had contracted the disease but 
could not transmit it and do not have symptoms of infection. The third 
compartment is the infected human $\left(I_{H}\right)$, who had contracted 
the rabies virus and were showing all the symptoms of rabies  and  were, 
therefore, infectious. The  fourth compartment are recovered human beings, 
denoted by  $R_{H}$, who get post-exposure prophylaxis after 
they contact and developed temporary immunity. 

The environment represents the virus causing rabies, that are in the 
physical object or any other materials in the environment. These virus 
containing material are considered as virus transmitting media 
or the infectious agent denoted as $M(t)$.

	
\subsubsection{Description of Model Interaction} 

Susceptible humans are constantly recruited at a rate of $\theta_{1}$. 
When they come into adequate contact with $I_{F}$, $I_{D}$, or the virus 
in the environment, individuals in the $S_{H}$ category become infected 
at rates $\tau_{1}$, $\tau_{2}$, and $\tau_{3}$, respectively:
\begin{equation*}
\begin{aligned}
\chi_{1}=\left(\tau_{1}I_{F}+\tau_{2}I_{D}
+\tau_{3} \lambda \left(M\right)\right)S_{H}.
\end{aligned}
\end{equation*}
After contracting rabies infections, susceptible humans become exposed individuals, 
denoted by $E_{H}$, typically for 1 to 3 months. Those in the $E_{H}$ category, 
who receive post-exposure prophylaxis, recover at the rate $\beta_{2}$. 
Since post-exposure prophylaxis does not confer permanent immunity, individuals 
in the $R_{H}$ category can lose immunity and become susceptible again 
at the rate $\beta_{3}$. The remaining proportion of the exposed class 
progresses to the infectious state $I_{H}$ at a rate $\beta_{1}$. 
Infected humans can die due to the disease at a rate $\sigma_{1}$. 
All human compartments experience natural death at a rate of $\mu_{1}$.
	
Susceptible free-range dogs are constantly recruited at a rate of $\theta_{2}$. 
Dogs in the $S_{F}$ category become infected when they come into adequate contact 
with $I_{F}$, $I_{D}$, or the rabies virus in the environment, 
at the rates $\kappa_{1}$, $\kappa_{2}$, and $\kappa_{3}$, respectively:
\begin{equation*}
\begin{aligned}
\chi_{2}=\left(\kappa_{1}I_{w}+\kappa_{2} I_{D}+\kappa_{3} \lambda \left(M\right)\right)S_{F}.
\end{aligned}
\end{equation*}
After contracting the rabies infections, susceptible free-range dogs 
progress to $E_{F}$ for 1 to 3 months. The exposed free-range dogs then 
progress to infectious state $I_{F}$ at a rate $\gamma$. The infected 
free-range dogs can die due to the disease at a rate $\sigma_{2}$. 
All free-range dog compartments experience natural death at a rate of $\mu_{2}$.
		
Susceptible domestic dogs are recruited constantly at a rate $\theta_{3}$; 
$S_{D}$  become infected when they either adequately contact $I_{F}$, $I_{D}$,  
or the virus in the environment at the rates $\psi_{1}$, $\psi_{2}$ and $\psi_{3}$, 
respectively, since domestic dogs are under human control:
\begin{equation*}
\begin{aligned}
\chi_{3}=\left(\dfrac{\psi_{1}I_{F}}{1+\rho_{1}}+\dfrac{\psi_{2}I_{D}}{1
+\rho_{2}}+\dfrac{\psi_{3}}{1+\rho_{3}}\lambda \left(M\right)\right) S_{D}.
\end{aligned}
\end{equation*}
After contracting a rabies infection, domestic dogs enter a state  $E_{D}$, 
where they stay for several months at a rate of $\beta_{1}$. If dogs in the 
$E_{D}$ category receive post-exposure prophylaxis, they move to the $R_{D}$ 
category at a rate of $\gamma_{2}$. However, post-exposure prophylaxis does 
not provide permanent immunity, and dogs in the $R_{D}$ category can lose 
their immunity and become susceptible again at a rate of $\gamma_{3}$. 
Meanwhile, the remaining proportion of the exposed class progresses 
to the infectious state $I_{D}$ at a rate of $\gamma_{1}$. Infected domestic 
dogs can die from the disease at a rate of $\sigma_{3}$. All compartments 
of domestic dogs experience natural death at a rate of $\mu_3$.

The virus causing rabies that are in the environment are recruited 
through shedding from the infectious free-range dogs, domestic dogs 
and humans at the rates $\nu_2$, $\nu_3$ and $\nu_1$, respectively: 
\begin{equation*}
\begin{aligned}
\theta_{4}=\left(\nu_1I_H+\nu_2I_F+\nu_3I_D\right)M.
\end{aligned}
\end{equation*}	
The viruses are removed from the environment at a rate $\mu_{4}$.
 
\begin{figure}[!ht]
\centering
\fbox{\includegraphics[width=1.0\linewidth, height=0.5\textheight]{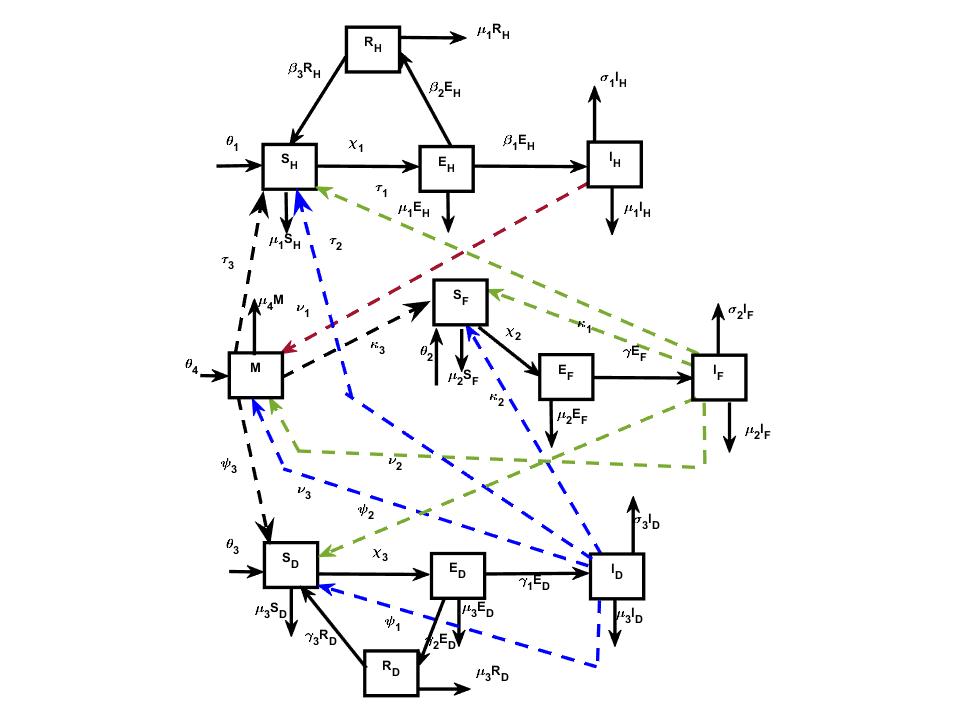}}
\caption{Flow diagram for rabies transmission among human, free range and domestic dogs with parameters.}
\label{fig:2}
\end{figure}

Figure~\ref{fig:2} provides a concise representation of the interactions 
among  the environment, humans, free-range dogs, and domestic dogs, 
which depict the dynamics of the rabies disease.

Based on the description of the model parameters and their connections to the state variables, 
we formulate a model in the form of a system of ordinary differential equations, 
as presented in the model equation (\ref{eqn1}):
\begin{center}
\begin{equation}
\left.
\begin{array}{llll}
\dot{S_{H}}&=\theta_{1}+\beta_{3}R_{H}-\mu_{1} S_{H}-\chi_{1}\\
\dot{E_{H}}&=\chi_{1}-\left(\mu_{1}+\beta_{1}+\beta_{2}\right)E_{H}\\
\dot{I_{H}}&=\beta_{1}E_{H}-\left(\sigma_{1}+\mu_{1}\right) I_{H}\\
\dot{R_{H}}&=\beta_{2} E_{H}-\left(\beta_{3}+\mu_{1} \right) R_{H}\\\\
\dot{S_{F}}&=\theta_{2}-\chi_{2}-\mu_{2}S_{F}\\
\dot{E_{F}}&=\chi_{2}-\left(\mu_{2}+\gamma\right)E_{F}\\
\dot{I_{F}}&=\gamma E_{F}-\left(\mu_{2}+\sigma_{2}\right)I_{F}\\\\
\dot{S_{D}}&=\theta_{3}-\mu_{3}S_{D}-\chi_{3}+\gamma_{3}R_{D}\\
\dot{E_{D}}&=\chi_{3}-\left(\mu_{3}+\gamma_{1}+\gamma_{2}\right) E_{D}\\
\dot{I_{D}}&=\gamma_{1}E_{D}-\left(\mu_{3}+\sigma_{3}\right) I_{D}\\
\dot{R_{D}}&=\gamma_{2}E_{D}-\left(\mu_{3}+\gamma_{3}\right)R_{D}\\\\
\dot{M}&=\left(\nu_1I_H+\nu_2I_F+\nu_3I_D\right)-\mu_4M	
\end{array}
\right\}
\label{eqn1}
\end{equation}
\end{center}
subject to non-negative conditions
\begin{equation*}
\begin{gathered}
S_{H}(0) > 0, \quad E_{H}(0) \geq 0,\quad I_{H}(0) \geq 0, \quad R_{H}(0) \geq 0, 
\quad S_{F}(0) > 0, \quad E_{F}(0) \geq 0, \quad I_{F}(0) \geq 0, \\
S_{D}(0) \geq 0, \quad E_{D}(0) \geq 0, \quad I_{D}(0) \geq 0, \quad R_{D}(0) \geq 0, 
\end{gathered}
\end{equation*}
and where
\begin{align*}
\chi_{1} &= \left(\tau_{1}I_{F} + \tau_{2}I_{D} + \tau_{3} \lambda(M)\right)S_{H},\;\\ 
\chi_{2} &= \left(\kappa_{1}I_{F} + \kappa_{2}I_{D} + \kappa_{3} \lambda(M)\right)S_{F},\\
\chi_{3} &= \left(\dfrac{\psi_{1}I_{F}}{1+\rho_{1}} + \dfrac{\psi_{2}I_{D}}{1+\rho_{2}} 
+ \dfrac{\psi_{3}}{1+\rho_{3}}\lambda(M)\right) S_{D},\;\\
\lambda\left(M\right) &= \dfrac{M}{M+C}.
\end{align*}
	
		
\subsection{Model analysis}

Since the model $\left(\ref{eqn1}\right)$ monitors the human and dog populations, 
we assume that the model's state variable and parameters are non-negative 
for $\forall  t\geq 0$. 
	
\subsubsection{Model's invariant region}

Theorem~\ref{The1} gives the invariant region of the rabies model.
	
\begin{theorem}
\textup{The solution of the rabies model system 
$\left( \ref{eqn1}\right)$ is uniformly bounded if 
$\Omega\in \mathbb{R}_{+}^{12}$ and $\Omega=\Omega_{H} \cup \Omega_{D}
\cup \Omega_{F} \cup \Omega_{M}\in\mathbb{R}_{+}^{4}
\times\mathbb{R}_{+}^{3}\times\mathbb{R}_{+}^{4}\times\mathbb{R}_{+}^{1}$, where}
\begin{align*}
\textup{$\Omega_{H}$} 
&= \left\{\left(S_{H}, E_{H}, I_{H}, R_{H}\right)\in\mathbb{R}_{+}^{4}: 
0\leq N_{H}\leq\dfrac{\theta_{1}}{\mu_{1}}\right\}, \\
\textup{$\Omega_{F}$} 
&= \left\{\left(S_{F}, E_{F}, I_{F} \right)\in\mathbb{R}_{+}^{3}: 
0\leq N_{F}\leq\dfrac{\theta_{2}}{\mu_{2}}\right\}, \\
\textup{$\Omega_{D}$} 
&= \left\{\left(S_{D}, E_{D}, I_{D}, R_{D}\right)\in\mathbb{R}_{+}^{4}: 
0\leq N_{D}\leq\dfrac{\theta_{3}}{\mu_{3}}\right\},
\end{align*}
\textup{and $\Omega$ is the positive invariant region.}
\label{The1}
\end{theorem}

\begin{proof}
Consider the population of the human from  equation $\left(\ref{eqn1} \right)$ as
\begin{eqnarray}
\dfrac{dN_{H}}{dt}
&=\dfrac{dS_{H}}{dt}+\dfrac{dE_{H}}{dt}+\dfrac{dI_{H}}{dt}+\dfrac{dR_{H}}{dt}.
\label{eqn7}
\end{eqnarray}
Then, the sum of the total population in  equation 
$\left(\ref{eqn7}\right)$  implies that
\begin{eqnarray}
\dfrac{dN_{H}}{dt}&= \theta_{1}
-\left(S_{H}+E_{H}+I_{H}+R_{H}\right)\mu_{1}-\sigma_{1}I_{H}.
\label{eqn8}
\end{eqnarray}  
Thus, the equation $\left(\ref{eqn8}\right)$ becomes
\begin{eqnarray}
\dfrac{dN_{H}}{dt}
&= \theta_{1}-\left(S_{H}+E_{H}+I_{H}+R_{H}\right)\mu_{1}.
\label{eqn10}
\end{eqnarray} 
But $N_{H}=S_{H}+E_{H}+I_{H}+R_{H}$ 
and equation $\left(\ref{eqn10}\right)$ is then  written as
\begin{eqnarray*}
\dfrac{dN_{H}}{dt}= \theta_{1}-N_{H}\mu_{1}.
\end{eqnarray*}
From the integrating  factor we have
\begin{eqnarray}
N_{H}\left(t\right) 
&=e^{\int\limits_{0}^{t} \mu_{1}dt}&= e^{\mu_{1} t}.
\label{13}
\end{eqnarray}
Then, for $t\rightarrow 0$ equation $\left(\ref{13} \right)$ gives
\begin{eqnarray}
N_{H}(0)\le\dfrac{\theta_{1}}{\mu_{1}}
+Ce^{0} \implies N_{H}(0)-\dfrac{\theta_{1}}{\mu_{1}}\le C.
\label{16}
\end{eqnarray}
By simplifying equation $\left(\ref{16}\right)$ 
and performing simple manipulations, we have
\begin{eqnarray*}
\Omega_{H} = \left\{\left(S_{H}, E_{H}, I_{H}, R_{H}\right)
\in\mathbb{R}_{+}^{4} : 0 \leq N_{H} \leq \dfrac{\theta_{1}}{\mu_{1}} \right\}.
\end{eqnarray*}
Using the same procedures for free-range and domestic dogs, we have   
\begin{eqnarray*}
\Omega_{w} = \left\{\left(S_{F}, E_{F}, I_{F}\right)\in\mathbb{R}_{+}^{3} 
: 0 \leq N_{F} \leq \dfrac{\theta_{2}}{\mu_{2}} \right\},\;\;
\Omega_{D} = \left\{\left(S_{D}, E_{D}, I_{D}, R_{D}\right)
\in\mathbb{R}_{+}^{4} : 0 \leq N_{D} \leq \dfrac{\theta_{3}}{\mu_{3}} \right\}. 
\end{eqnarray*}

Regarding the environment that contains the rabies virus, 
we consider equation $\left(\ref{16}\right)$ 
in the model system $\left(\ref{eqn1}\right)$ as
\begin{eqnarray}
\dot{M}&=\left(\nu_1I_H+\nu_2I_F+\nu_3I_D\right)-\mu_4M.
\label{eqn25}
\end{eqnarray}
Since $N_H\leq\dfrac{\theta_1}{\mu_1}$, 
$N_F\leq\dfrac{\theta_2}{\mu_3}$, 
and $N_D\leq\dfrac{\theta_3}{\mu_3}$, it follows that 
$I_H\leq\dfrac{\theta_1}{\mu_1}$,  
$I_F\leq\dfrac{\theta_2}{\mu_2}$, 
and $I_D\leq\dfrac{\theta_3}{\mu_3}$.
Thus, equation \eqref{eqn25} becomes
\begin{eqnarray*}
\dot{M}&\leq\left(\dfrac{\nu_1\theta_1}{\mu_1}
+\dfrac{\nu_2\theta_2}{\mu_2}+\frac{\nu_3\theta_3}{\mu_3}\right)-\mu_4M.
\end{eqnarray*}
Now, let $Y$ be the unique solution to the initial value problem, such that
\begin{equation}
\left.
\begin{array}{llll}
\dot{Y}&\leq\left(\dfrac{\nu_1\theta_1}{\mu_1}+\dfrac{\nu_2\theta_2}{\mu_2}
+\dfrac{\nu_3\theta_3}{\mu_3}\right)-\mu_4M,\; \textit{for}\;\;\; t>0,\\
Y(0)&=M(0).
\end{array}
\right\}
\label{eqn28}
\end{equation}
By using the integration factor, equation $\left(\ref{eqn28}\right)$ becomes
\begin{eqnarray}
M\left(t\right)\le \left(\dfrac{\nu_1\theta_1}{\mu_1}
+\dfrac{\nu_2\theta_2}{\mu_2}+\dfrac{\nu_3\theta_3}{\mu_3}\right)
\dfrac{1}{\mu_{4}}+\left(M(0)-\left(\dfrac{\nu_1\theta_1}{\mu_1}
+\dfrac{\nu_2\theta_2}{\mu_2}+\dfrac{\nu_3\theta_3}{\mu_3}\right)\dfrac{1}{\mu_{4}}\right).
\label{7}
\end{eqnarray}
As $t\rightarrow \infty$, the expression 
$\left(M(0)-\left(\frac{\nu_1\theta_1}{\mu_1}+\frac{\nu_2\theta_2}{\mu_2}
+\frac{\nu_3\theta_3}{\mu_3}\right)\frac{1}{\mu_{4}}e^{\mu_{4} t}\right)$  
in equation $\left(\ref{7}\right)$ goes to zero, and we have
\begin{equation}
M\le \left(\dfrac{\nu_1\theta_1}{\mu_1}
+\dfrac{\nu_2\theta_2}{\mu_2}
+\dfrac{\nu_3\theta_3}{\mu_3}\right)\dfrac{1}{\mu_{4}}.
\label{eqn31}
\end{equation}
Equation $\left(\ref{eqn31}\right)$ gives
\begin{eqnarray*} 
M\left(t\right)\le\Omega_{M}  
= \max\left\{\dfrac{\theta_1 \nu_1}{\mu_1\mu_4}
+\dfrac{\theta_2 \nu_2}{\mu_2\mu_4}
+\dfrac{\theta_3 \nu_3}{\mu_3\mu_4},M\left(0\right)\right\}. 
\end{eqnarray*}
Thus, the model system $\left(\ref{eqn1}\right)$ 
is biologically and mathematically meaningful 
with its solution relying on the region $\Omega$.
\end{proof}


\subsubsection{Positivity of the solution}	

For the model system $\left(\ref{eqn1}\right)$ to be epidemiologically meaningful 
and well-posed, we need to prove that the state variables are non-negative for all $t \geq 0$.

\begin{theorem}
\textup{Let} $\left\{\textup{S}_{H}(0), \textup{E}_{H}(0), 
\textup{I}_{H}(0), \textup{R}_{H}(0), \textup{S}_{F}(0), 
\textup{E}_{F}(0), \textup{I}_{F}(0), \textup{S}_{D}(0), 
\textup{E}_{D}(0), \textup{I}_{D}(0), \textup{R}_{D}(0),M(0)\right\} 
\in \mathbb{R}_{+}^{12}$. \textup{Then the set of solutions}
\[
\left\{\textup{S}_{H}(t), \textup{E}_{H}(t), 
\textup{I}_{H}(t), \textup{R}_{H}(t), \textup{S}_{F}(t), 
\textup{E}_{F}(t), \textup{I}_{F}(t), 
\textup{S}_{D}(t), \textup{E}_{D}(t), 
\textup{I}_{D}(t), \textup{R}_{D}(t),M(t)\right\}
\]
\textup{of the model system $\left(\ref{eqn1}\right)$ is positive} $\forall t > 0$.
\label{Th2}
\end{theorem}

\begin{proof}
Consider the human subpopulation of the model system $\left(\ref{eqn1}\right)$. 
We have
\begin{eqnarray}
\dfrac{dS_{H}}{dt}
&=\theta_{1}+\lambda_{3}R_{H}-\mu_{1} S_{H}-\left(\tau_{1}I_{F}
+\tau_{2}I_{D}+\tau_{3} \lambda \left(M\right)\right)S_{H}.
\label{eqn17}
\end{eqnarray}
Then equation $\left(\ref{eqn17}\right)$ is presented as
\begin{eqnarray}
\dfrac{dS_{H}}{dt}\ge-\left(\tau_{1}I_{F}+\tau_{2}I_{D}
+\tau_{3} \lambda \left(M\right)+\mu_{1}\right) S_{H}
\label{eqn19}
\end{eqnarray}
which, separating and integrating equation $\left(\ref{eqn19}\right)$ on both sides, gives
\begin{eqnarray}
\int\dfrac{dS_{H}}{ S_{H}}\ge\int\limits_{0}^{t}-\left(\tau_{1}I_{F}
+\tau_{2}I_{D}+\tau_{3} \lambda \left(M\right)+\mu_{1}\right)  ds.
\label{eqn20}
\end{eqnarray}
Therefore, from equation $\left(\ref{eqn20}\right)$, we have
\begin{eqnarray}
S_{H}\ge S_{H}\left(0\right)e^{\int\limits_{0}^{t}-\left(\tau_{1}I_{F}
+\tau_{2}I_{D}+\tau_{3} \lambda \left(M\right)+\mu_{1}\right)  ds} >0.
\label{21}
\end{eqnarray}
Thus, $S_{H}$ is positive for all $t > 0$. Using the same procedure 
from equations $\left(\ref{eqn17}\right)$ to $\left(\ref{21}\right)$, we have
\begin{align*}
E_{H} &\ge E_{H}(0)e^{-(\mu_{1} + \beta_{1} + \beta_{2})t} > 0 ,\;
I_{H} \ge I_{H}(0)e^{-(\sigma_{1} + \mu_{1})t} > 0,\;
R_{H} \ge R_{H}(0)e^{-(\mu_{1} + \beta_{3})t} > 0, \\
S_{F} &\ge S_{F}(0)e^{\int_{0}^{t} -((\kappa_{1}I_{F} 
+ \kappa_{2}I_{D} + \kappa_{3}\lambda(M)) + \mu_{2})ds} > 0 ,\;
E_{F} \ge E_{F}(0)e^{-(\mu_{2} + \gamma)t} > 0, \;
I_{F} \ge I_{F}(0)e^{-(\mu_{2} + \sigma_{2})t} > 0, \\
S_{D} &\ge S_{D}(0)e^{\int_{0}^{t} -((\frac{\psi_{1}I_{F}}{1+\rho_{1}} 
+ \frac{\psi_{2}I_{D}}{1+\rho_{2}} 
+ \frac{\psi_{3}}{1+\rho_{3}}\lambda(M) + \mu_{3}))ds} > 0,\;
E_{D} \ge E_{D}(0)e^{-(\mu_{3} + \gamma_{1} + \gamma_{2})t} > 0, \\
R_{D} &\ge R_{D}(0)e^{-(\mu_{3} + \gamma_{3})t} > 0 ,\;
I_{D} \ge I_{D}(0)e^{-(\mu_{3} + \sigma_{3})t} > 0,\;
M \ge M(0)e^{-\mu_{4}t} > 0. 
\end{align*}
Thus, the set of solutions $\left\{S_{H}(t),E_{H}(t),I_{H}(t),R_{H}(t),S_{F}(t),
E_{F}(t),I_{F}(t),S_{D}(t),E_{D}(t),I_{D}(t),R_{D}(t),M(t)\right\}$  
of the model system $\left(\ref{eqn1}\right)$ is positive $\forall t>0$.
\end{proof}


\subsubsection{Disease free equilibrium point (DFE)}

The DFE, denoted as ${\cal E}_{0}$, is defined as the point at which 
there is no disease in a given population. To obtain ${\cal E}_{0}$ in the model 
system $\left(\ref{eqn1}\right)$, we consider all infectious 
compartments to be equal to zero:  
$$
E_{H}=I_{H}=E_{F}=I_{F}=E_{D}=I_{D}=M=0  
\implies 
\dfrac{dE_{H}}{dt}=\dfrac{dI_{H}}{dt}
=\dfrac{dE_{F}}{dt}=\dfrac{dI_{F}}{dt}
=\dfrac{dE_{D}}{dt}=\dfrac{dI_{D}}{dt}
=\dfrac{d M}{dt}=0
$$  
such that  
\begin{equation}
\left.
\begin{array}{llll}
0&=\theta_{1}-\mu_{1}S_{H},\\
0&=\theta_{2}-\mu_{2}S_{F},\\
0&=\theta_{3}-\mu_{3}S_{D}.
\end{array}
\right\}
\label{eqn22}
\end{equation}
Hence, the DFE point ${\cal E}_{0}$ is derived through the mathematical 
rearrangement of the equation $\left(\ref{eqn22}\right)$ as 
$$
{\cal E}_{0} =\left(S_{H}^{0},E_{H}^{0},I_{H}^{0},R_{H}^{0},
S_{F}^{0},E_{F}^{0},I_{F}^{0},S_{D}^{0},E_{D}^{0},I_{D}^{0},R_{H}^{0},M^{0}\right) 
=\left(\frac{\theta_{1}}{\mu_{1}},0,0,0,
\frac{\theta_{2}}{\mu_{2}},0,0,\frac{\theta_{3}}{\mu_{3}},0,0,0,0\right).
$$


\subsubsection{The basic reproduction number ${\cal R}_0$}

The ${\cal R}_0$ predicts whether rabies will spread across the community or die out. 
If ${\cal R}_0 < 1$, it means that every infectious individual will cause less than 
one secondary infection, hence the disease will die out in the community. 
If ${\cal R}_0 > 1$, it means that every infectious individual will cause more 
than one secondary infection, leading to the persistence of rabies 
in the entire population. In order to determine ${\cal R}_0 > 1$, the next 
generation operator, as applied by 
\cite{lasalle1976stability,yang2014basic,saha2021dynamics}, 
and the Jacobian Matrix are used:
\begin{eqnarray*}
\frac{dx_{i}}{dt} 
&=\mathcal{F}_{i}\left(x\right)-\left(\mathcal{V}_{i}^{+}
\left(x\right)-\mathcal{V}_{i}^{-}\left(x\right)\right),
\end{eqnarray*} 
where $\mathcal{F}_{i}$ is the new infections in the compartment $i$  
while  $\mathcal{V}_{i}^{+}$ and $\mathcal{V}_{i}^{-}$ are the transfer 
terms in and out of the compartment $i$, respectively. From equation 
$\left(\ref{eqn1} \right)$ we define
$\mathcal{F}_{i}$ and $\mathcal{V}_{i}$ by
\begin{eqnarray}
\mathcal{F}_{i}=
\left(
\begin{array}{c}
\left(\tau_{1}I_{F}+\tau_{2}I_{D}+\tau_{3} \lambda \left(M\right)\right)S_{H}   \\
0\\
\left(\kappa_{1}I_{F}+\kappa_{2} I_{D}+\kappa_{3} \lambda \left(M\right)\right)S_{F}\\
0\\
\left(\dfrac{\psi_{1}I_{F}}{1+\rho_{1}}+\dfrac{\psi_{2}I_{D}}{1+\rho_{2}}
+\dfrac{\psi_{3}}{1+\rho_{3}}\lambda \left(M\right)\right) S_{D}\\
0\\
0
\end{array}
\right),\quad
\mathcal{V}_{i}=
\left(
\begin{array}{c}
\left(\mu_{1}+\beta_{1}+\beta_{2}\right)E_{H} \\
\left(\sigma_{1}+\mu_{1}\right)I_{H}-\beta_{1}E_{H}\\
\left(\mu_{2}+\gamma\right)E_{F}\\
\left(\mu_{2}+\sigma_{2}\right)I_{F}- \gamma E_{F}\\
\left(\mu_{3}+\gamma_{1}+\gamma_{2}\right) E_{D}\\
\left(\mu_{3}+\delta_{3}\right) I_{D}-\gamma_{1}E_{D}\\
\mu_4M-\left(\nu_1I_H+\nu_2I_F+\nu_3I_D\right)
\end{array}
\right).
\label{M}
\end{eqnarray}
The Jacobian Matrices $F$ and $V$ at the disease free equilibrium point 
${\cal E}_{0}$ is given by equation $\left(\ref{Next}\right)$:
\begin{eqnarray}
\begin{aligned}
F&=\dfrac{\partial\mathcal{F}_{i}\left({\cal E}_{0}\right)}{\partial x_{j}}, 
\quad V&=\dfrac{\partial\mathcal{V}_{i}\left({\cal E}_{0}\right)}{\partial x_{j}}.
\label{Next}
\end{aligned}
\end{eqnarray}
From  equation $\left(\ref{M}\right)$,  $F$ at DFE point ${\cal E}_0$ 
is given in equation $\left(\ref{job1}\right)$:
\begin{eqnarray}
\dfrac{\partial\mathcal{F}_{i}\left({\cal E}_{0}\right)}{\partial x_{j}}=F=
\left(
\begin{array}{ccccccc}
0&0&0&\dfrac{\tau_{1}\theta_{1}}{\mu_{1}}&0&\dfrac{\tau_{2}\theta_{1}}{\mu_{1}}&0\cr
0&0&0&0&0&0&0\cr
0&0&0&\dfrac{\kappa_{1}\theta_{2}}{\mu_{2}}&0&\dfrac{\kappa_{2}\theta_{2}}{\mu_{2}}&0\cr
0&0&0&0&0&0&0\cr
0&0&0&\dfrac{\psi_{1}\theta_{3}}{\left(1+\rho_{1}\right)\mu_{3}}
&0&\dfrac{\psi_{2}\theta_{3}}{\left(1+\rho_{2}\right)\mu_{3}}&0\cr
0&0&0&0&0&0&0\cr
0&0&0&0&0&0&0
\end{array}
\right).
\label{job1}
\end{eqnarray} 

In the linearized system $\left(\ref{job1}\right)$, the entry $F_{ij}$ 
represents the rate at which individuals in the infected state $j$ give 
rise to or develop new infections in individuals in the infected state $i$, 
with reference to the infected states indexed by $i$ and $j$ 
for $i,j\in\left\{1,2,3,4,5,6,7\right\}$. As a result, when a person 
in an infected condition $j$ does not instantly produce any 
new instances in an infected state $i$, we have $F_{ij}=0$. 
Similarly, $V$ at DFE point ${\cal E}_0$ 
is given in equation $\left(\ref{jobi2}\right)$:
\begin{eqnarray}
\dfrac{\partial\mathcal{V}_{i}\left({\cal E}_{0}\right)}{\partial x_{j}}
=V= \setlength{\arraycolsep}{1.5pt}
\left(
\begin{array}{ccccccc}
\mu_{1}+\beta_{1}+\beta_{2}&0&0&0&0&0&0\cr
-\beta_{1}&\sigma_{1}+\mu_{1}&0&0&0&0&0\cr
0&0&\mu_{2}+\gamma&0&0&0&0\cr
0&0&-\gamma &\mu_{2}+\sigma_{2}&0&0&0\cr
0&0&0&0&\mu_3+\gamma_{1}+\gamma_{2}&0&0\cr
0&0&0&0&-\gamma_{1}&\mu_{3}+\sigma_{3}&0\cr
0&-\nu_{1}&0&-\nu_{2}&0&-\nu_{3}&\mu_{4}
\end{array}
\right).
\label{jobi2}
\end{eqnarray} 
 
The inverse of matrix $V$ is easily obtained 
using \textsf{Maple} software, and its result is given as
\begin{eqnarray}
V^{-1}=
\setlength{\arraycolsep}{1.5pt}
\left(
\begin{array}{ccccccc}
\frac{1}{\mu_{1}+\beta_{1}+\beta_{2}}&0&0&0&0&0&0\cr
\noalign{\medskip}{\frac {\beta_{{1}}}{ \left( \sigma_{{1}}
+\mu_{{1}} \right)  \left( \mu_{{1}}+\beta_{{1}}
+\beta_{{2}} \right) }}& \frac{1}{\sigma_{{1}}+\mu_{{1}} }&0
&0&0&0&0\cr 
\noalign{\medskip}0&0& \frac{1}{\mu_{{2}}+\gamma}&0&0&0&0\cr 
\noalign{\medskip}0&0&{\frac {\gamma}{ \left( \mu_{{2}}+
\gamma \right)  \left( \mu_{{2}}+\sigma_{{2}} \right) }}& \dfrac{1}{\mu_{
{2}}+\sigma_{{2}}}&0&0&0\cr \noalign{\medskip}0&0&0&0&
\frac{1}{\mu_{{3}}+\gamma_{{1}}+\gamma_{{2}}}&0&0
\cr \noalign{\medskip}0&0&0&0&{\dfrac {\gamma_{{1}}}{ \left( \mu_{{3}}+
\sigma_{{3}} \right)  \left( \mu_{{3}}+\gamma_{{1}}+\gamma_{{2}}
\right) }}& \dfrac{1}{\mu_{{3}}+\sigma_{{3}}}&0
\cr \noalign{\medskip}{\frac {\nu_{{1}}\beta_{{1}}}{ \left( \sigma_{{1}}
+\mu_{{1}} \right)  \left( \mu_{{1}}+\beta_{{1}}+\beta_{{2}} \right) 
\mu_{{4}}}}&{\frac {\nu_{{1}}}{ \left( \sigma_{{1}}+\mu_{{1}} \right) 
\mu_{{4}}}}&{\frac {\nu_{{2}}\gamma}{ \left( \mu_{{2}}+\gamma \right) 
\left( \mu_{{2}}+\sigma_{{2}} \right) \mu_{{4}}}}&{\frac {\nu_{{2}}}{
\left( \mu_{{2}}+\sigma_{{2}} \right) \mu_{{4}}}}&{\frac {\nu_{{3}}
\gamma_{{1}}}{ \left( \mu_{{3}}+\sigma_{{3}} \right)  \left( \mu_{{3}}
+\gamma_{{1}}+\gamma_{{2}} \right) \mu_{{4}}}}&{\frac {\nu_{{3}}}{
\left( \mu_{{3}}+\sigma_{{3}} \right) \mu_{{4}}}}&\frac{1}{\mu_{{4}}}
\end{array}
\right).
\label{Jacob3}
\end{eqnarray}
In the context of computing the basic reproduction number ${\cal R}_{0}$ 
in epidemiology, the $(V^{-1})_{ij}$ obtained in equation 
$\left(\ref{Jacob3}\right)$ represents the generation matrix. 
The generation matrix describes the expected number of newly 
infected individuals that a single infectious individual 
in each of the different susceptible classes generates. Meanwhile, 
the diagonal elements represent the rate of leaving the corresponding 
susceptible class due to other causes such as recoveries caused by 
administration of post-exposure prophylaxis (PEP). In particular, 
$\dfrac{1}{\mu_{1}+\beta_{1}+\beta_{2}}$, $\dfrac{1}{\mu_{2}+\gamma}$, 
and $\dfrac{1}{\mu_{3}+\gamma_{1}+\gamma_{2}}$, respectively, represent 
the average incubation period for rabies in humans, free-range dogs, 
and domestic dogs. Meanwhile, $\dfrac{1}{\mu_{1}+\sigma_{1}}$, 
$\dfrac{1}{\mu_{2}+\sigma_{2}}$, and $\dfrac{1}{\mu_{3}+\sigma_{3}}$, 
respectively, represent the average time spent by an infective human, 
free-range dog, and domestic dog in the infectious state, and 
$\dfrac{1}{\mu_{4}}$ is the average time the rabies virus 
spends in the environment. The next-generation matrix is then calculated by
\begin{eqnarray*}
F V^{-1}= \setlength{\arraycolsep}{1.5pt}
\left(
\begin{array}{ccccccc}
0&0&{\frac {\tau_{{1}}\theta_{{1}}
\gamma}{\mu_{{1}} \left( \mu_{{2}}+\gamma \right)  \left( \mu_{{2}}+
\sigma_{{2}} \right) }}&{\frac {\tau_{{1}}\theta_{{1}}}{\mu_{{1}}
\left( \mu_{{2}}+\sigma_{{2}} \right) }}&{\dfrac {\tau_{{2}}
\theta_{{1}}\gamma}{\mu_{{1}} \left( \mu_{{3}}+\gamma_{{1}}+\gamma_{{2}}
\right)  \left( \mu_{{3}}+\sigma_{{3}} \right) }}&{\frac {\tau_{{2}}
\theta_{{1}}}{\mu_{{1}} \left( \mu_{{3}}+\sigma_{{3}} \right) }}&0
\cr \noalign{\medskip}0&0&0&0&0&0&0\cr \noalign{\medskip}0&0&{\frac {
\kappa_{{1}}\theta_{{2}}\gamma}{\mu_{{2}} \left( \mu_{{2}}+\gamma
\right)  \left( \mu_{{2}}+\sigma_{{2}} \right) }}&{\frac {\kappa_{{1}}
\theta_{{2}}}{\mu_{{2}} \left( \mu_{{2}}+\sigma_{{2}} \right) }}&{
\frac {\kappa_{{1}}\theta_{{2}}\gamma}{\mu_{{2}} \left( \mu_{{3}}
+\gamma_{{1}}+\gamma_{{2}} \right)  \left( \mu_{{3}}+\sigma_{{3}}
\right) }}&{\frac {\kappa_{{1}}\theta_{{2}}}{\mu_{{2}} \left( \mu_{{3}}
+\sigma_{{3}} \right) }}&0\cr \noalign{\medskip}0&0&0&0&0&0&0
\cr \noalign{\medskip}0&0&{\dfrac {\psi_{{1}}\theta_{{3}}\gamma}{
\left( 1+\rho_{{1}} \right) \mu_{{3}} \left( \mu_{{2}}+\gamma
\right)  \left( \mu_{{2}}+\sigma_{{2}} \right) }}&{\frac {\psi_{{1}}
\theta_{{3}}}{ \left( 1+\rho_{{1}} \right) \mu_{{3}} \left( \mu_{{2}}+
\sigma_{{2}} \right) }}&{\dfrac {\psi_{{2}}\theta_{{3}}\gamma}{
\left( 1+\rho_{{2}} \right) \mu_{{3}} \left( \mu_{{3}}+\gamma_{{1}}+
\gamma_{{2}} \right)  \left( \mu_{{3}}+\sigma_{{3}} \right) }}&{\frac 
{\psi_{{2}}\theta_{{3}}}{ \left( 1+\rho_{{2}} \right) \mu_{{3}}
\left( \mu_{{3}}+\sigma_{{3}} \right) }}&0\cr \noalign{\medskip}0&0&0&0
&0&0&0\cr \noalign{\medskip}0&0&0&0&0&0&0
\end{array}
\right).
\end{eqnarray*}
The expression of matrix $F V^{-1}$ can be presented as
\begin{equation}
\begin{aligned}
\setlength{\arraycolsep}{10pt}
F V^{-1} = \begin{pmatrix}
0 & 0 & R_{13} & R_{14} & R_{15} & R_{16} & 0 \\
0 & 0 & 0 & 0 & 0 & 0 & 0 \\
0 & 0 & R_{33} & R_{34} & R_{35} & R_{36} & 0 \\
0 & 0 & 0 & 0 & 0 & 0 & 0 \\
0 & 0 & R_{53} & R_{54} & R_{55} & R_{56} & 0 \\
0 & 0 & 0 & 0 & 0 & 0 & 0 \\
0 & 0 & 0 & 0 & 0 & 0 & 0 \\
\end{pmatrix},
\label{eqn24}
\end{aligned}
\end{equation}
where
\begin{equation}
\left.
\begin{array}{llll}
R_{13}={\dfrac {\tau_{{1}}\theta_{{1}}
\gamma}{\mu_{{1}} \left( \mu_{{2}}+\gamma \right)  \left( \mu_{{2}}
+\sigma_{{2}} \right) }},\;\; R_{14}={\frac {\tau_{{1}}\theta_{{1}}}{\mu_{{1}}
\left( \mu_{{2}}+\sigma_{{2}} \right) }},\;\;\;
R_{15}={\dfrac {\tau_{{2}}\theta_{{1}}\gamma}{\mu_{{1}} 
\left( \mu_{{3}}+\gamma_{{1}}+\gamma_{{2}}\right)  
\left( \mu_{{3}}+\sigma_{{3}} \right) }},\;\;\; R_{16}= {\dfrac {\tau_{{2}}
\theta_{{1}}}{\mu_{{1}} \left( \mu_{{3}}+\sigma_{{3}} \right) }},\\       
R_{33}=\dfrac{\kappa_{1}\theta_{2}\gamma}{\mu_{2}\left(\mu_{2}
+\gamma\right)\left(\mu_{2}+\sigma_{2}\right)},\;\;\; R_{34}={\dfrac {\kappa_{{1}}
\theta_{{2}}}{\mu_{{2}} \left( \mu_{{2}}+\sigma_{{2}} \right) }},\;\;\; 
R_{35}=\dfrac{\kappa_{1}\theta_{2}\gamma}{\mu_{2}\left(\mu_{3}
+\gamma_{1}+\gamma_{2}\right)\left(\mu_{3}+\sigma_{3}\right)},\;\;
R_{36}={\dfrac {\kappa_{{1}}\theta_{{2}}}{\mu_{{2}} \left( \mu_{{3}}
+\sigma_{{3}} \right) }}, \\ 
R_{53}=\dfrac{\psi_{1}\theta_{3}\gamma}{\left(1+\rho_{1}\right)
\left(\mu_{2}+\gamma\right)\left(\mu_{2}+\sigma_{2}\right)\mu_{3}} , \;\;\; 
R_{54}={\dfrac {\psi_{{1}}
\theta_{{3}}}{ \left( 1+\rho_{{1}} \right) \mu_{{3}} \left( \mu_{{2}}+
\sigma_{{2}} \right) }},\;\;\;  
R_{55}=\dfrac{\psi_{2}\theta_{3}\gamma}{\left(1+\rho_{2}\right)\left(\mu_{3}
+\gamma_{1}+\gamma_{2}\right)\left(\mu_{3}+\sigma_{3}\right)\mu_{3}},\\
R_{56}={\dfrac 
{\psi_{{2}}\theta_{{3}}}{ \left( 1+\rho_{{2}} \right) \mu_{{3}}
\left( \mu_{{3}}+\sigma_{{3}} \right) }}.
\end{array}
\right\}
\label{eqn6}
\end{equation}
From equation $\left(\ref{eqn24}\right)$, 
we obtain the eigenvalues as
\begin{gather*}
\begin{aligned}
&\left(\begin {array}{c}
0 \\ 
0 \\
0 \\
0 \\
0 \\
\frac{1}{2} R_{{55}} + \frac{1}{2} R_{{33}} + \frac{1}{2} \sqrt{R_{{33}}^2 
- 2 R_{{33}} R_{{55}} + 4 R_{{35}} R_{{53}} + R_{{55}}^2} \\
\frac{1}{2} R_{{55}} + \frac{1}{2} R_{{33}} - \frac{1}{2} 
\sqrt{R_{{33}}^2 - 2 R_{{33}} R_{{55}} + 4 R_{{35}} R_{{53}} + R_{{55}}^2}
\end{array}\right).
\end{aligned}
\end{gather*}   
The $\left(i,\ k\right)$ element of the $FV^{-1}$ of the next generation matrix   
represents the expected number of secondary infections in the compartment $i$ 
caused  by individuals in compartment $k$, assuming that the individual's 
environments remain consistent throughout the infection. It is worth noting that 
the $FV^{-1}$ matrix is non-negative, meaning it has a non-negative eigenvalue. 
This non-negative eigenvalue corresponds to a non-negative eigenvector that represents 
the distribution of infected individuals who generate the highest number of secondary 
infections per generation, also known as ${\cal R}_0$. 
According to \cite{dharmaratne2020estimation},
the basic reproduction number ${\cal R}_0$ is the largest eigenvalue 
of the next generating matrix, being given by
\begin{equation*}
{\cal R}_0=\rho \left(FV^{-1}\right).
\end{equation*}  
Therefore, the spectral radius of the next generation matrix is 
\begin{equation}
\rho\left(FV^{-1}\right)=\frac{\left(R_{55}+R_{33}\right)
+\sqrt{R_{33}\left(R_{33}-2R_{55}\right)+4R_{35}R_{53}+R_{55}^{2}}}{2}.
\label{eqn5}
\end{equation}
 
 
\subsubsection{Local sensitivity analysis}
 
Local sensitivity analysis seeks to determine how each parameter $\mathcal{P}_{i}$
affects the $\mathcal{R}_{0}$ and is determined by normalizing the sensitivity indices 
as approached by \cite{zhang2011analysis,asamoah2017modelling}: 
\begin{equation}
\gamma_{\mathcal{P}_{i}}^{\mathcal{R}_{0}}
=\frac{\partial \mathcal{R}_{0}}{\partial \mathcal{P}_{i}}
\times\frac{\mathcal{P}_{i}}{\mathcal{R}_{0}},
\label{eq22}
\end{equation}
where $\mathcal{R}_{0}$ is the rabies basic reproduction number. 
Therefore, utilizing equation $\left(\ref{eq22}\right)$ and the parameter values 
from Table~\ref{table3}, 
\begin{center}
\begin{longtable}{llll}
\caption{Model parameters, their description and values.}\\ \hline 
\textbf{Parameters}&\textbf{Description}&\textbf{Value (Year$^{-1}$)}&\textbf{Source}\\
\hline\hline
\endfirsthead
\multicolumn{2}{c}%
{\tablename\ \thetable\ -- \textit{Continued from previous page}} \\
\hline
\textbf{Parameters}&\textbf{Descriptions}\\[0.5ex] 
\hline
\endhead
\hline \multicolumn{2}{r}{\textit{Continued on next page}} \\
\endfoot
\hline
\endlastfoot
$\theta_{1}$ & Recruitment rate $S_{H}$ &2000& (Estimated) \\
$\tau_{1}$ & The rate that $S_{H}$ gets infection from $I_F$ & 0.0004
& \cite{tian2018transmission}\\
$\tau_{2}$ & The rate that $S_{H}$ gets infection from $I_D$& 0.0004
& \cite{tian2018transmission}\\
$\tau_{3}$ & The rate that $S_{H}$ gets infection from $M$
&$\left[0.0003\;\;  0.0100\right]$ & (Estimated)\\
$\beta_{1}$ & Progression rate out of $E_{H}$ to $I_H$ &$\frac{1}{6}$ 
& \cite{tian2018transmission,zhang2011analysis}\\
$\beta_{2}$ & Recovery rate of  $E_{H}$ &$\left[0.54 \;\; 1\right]$ 
& \cite{abdulmajid2021analysis,zhang2011analysis}\\
$\beta_{3}$ & Rate of immunity loss of humans &1& (Estimated)\\
$\mu_1$ &  Natural death rate  of humans & 0.0142
& \cite{world2010working,world2013expert}\\
$\sigma_1$ & Disease induced death rate for $I_H$ & 1
& \cite{abdulmajid2021analysis,zhang2011analysis}\\
$\theta_{2}$ &  Recruitment rate of free-range dogs &1000& (Estimated)\\
$\kappa_{1}$ & The rate that $S_{F}$ gets infection from $I_F$&  0.00006 & (Estimated)\\
$\kappa_{2}$ & The rate that $S_{F}$ gets infection from $I_D$ & 0.00005 & (Estimated)\\
$\kappa_{3}$ & The rate that $S_{F}$ gets infection from $M$
&$\left[0.00001 \;\; 0.00003\right]$&(Estimated)\\
$\gamma$ & The rate that $S_{F}$ gets infection from $I_F$&$\frac{1}{6}$ 
& \cite{tian2018transmission,abdulmajid2021analysis,zhang2011analysis}\\
$\sigma_2$ &  Disease induced death rate of  $I_F$ & 0.09 
& \cite{zhang2011analysis,addo2012seir}\\
$\mu_{2}$ &  Natural mortality rate of  free-range dogs &0.067 & (Estimated)\\
$\theta_3$ &  Recruitment rate of domestic dog population & 1200& (Estimated)\\
$\psi_{1}$ & The rate that $S_{D}$ gets infection from $I_D$ &0.0004
& \cite{hampson2019potential,addo2012seir}\\
$\psi_{2}$ & The rate that $S_{D}$ gets infection from $I_F$ & 0.0004 
& \cite{hailemichael2022effect}\\
$\psi_{3}$ & The rate that $S_{D}$ gets infection from $M$ & 0.0003 & (Estimated)\\
$\mu_3$ & Natural death rate for domestic dog population&0.067& (Estimated)\\
$\sigma_3$ & Disease induced death rate for $I_D$&0.08& \cite{zhang2011analysis}\\
$\gamma_1$ & The rate at which $E_{D}$ becomes $I_D$ & $\frac{1}{6}$ 
& \cite{tian2018transmission,zhang2011analysis}\\
$\gamma_2$ & Recovery rate of  $E_{D}$&0.09 & \cite{zhang2011analysis}\\
$\gamma_3$ & Rate of loss of temporary immunity for $R_{D}$ &0.05&(Estimasted)\\
$\nu_1$ & Environmental virus shedding rate from $I_H$ &0.001& (Estimated)\\
$\nu_2$ & Environmental virus shedding rate from $I_F$ &0.006& (Estimated)\\
$\nu_3$ & Environmental virus shedding rate from  $I_D$ &0.001& (Estimated)\\
$\mu_4$  & Natural removal rate of rabies from the environment&0.08& (Estimated)\\
$\rho_{1}$ &  The deterrent coefficient of domestic dog from  $I_F$&10
& \cite{ruan2017spatiotemporal}\\
$\rho_{2}$ & The deterrent coefficient of domestic dog from $I_D$ & 8& (Estimated)\\
$\rho_{3}$ & The deterrent coefficient of domestic dog from $M$&15 & (Estimated)\\
$C$  &  Concentration of rabies in the environment&0.003  (PFU)/mL & (Estimated)
\label{table3}
\end{longtable}
\end{center}
we calculate the sensitivity indices for each parameter, 
as indicated in Table~\ref{table4}.
\begin{table}[H]
\caption{Sensitivity indices for $\mathcal{R}_{0}$.} 
\centering 
\begin{tabular}{c c c c} 
\hline\hline 
Parameter & Sensitivity Index & Parameter & Sensitivity Index \\ [0.5ex] \hline 
$\gamma_{1}$ & -0.105552 & $\psi_{1}$ & +0.051422 \\
$\gamma_{2}$ & -0.056998 & $\psi_{2}$  & +0.005436 \\
$\kappa_{1}$ & +0.897120 & $\kappa_{2}$ & +0.051422 \\
$\mu_{2}$ & -1.616021 & $\rho_{1}$ & -0.046747\\
$\mu_{3}$ & -0.105358 & $\rho_{2}$ & -0.004832 \\
$\sigma_{2}$ & -0.540654 & $\sigma_{3}$ & -0.05144\\ 
$\theta_{2}$ & +0.941420 & $\theta_{3}$ & +0.056858\\[1ex] \hline 
\end{tabular}
\label{table4} 
\end{table}

The findings from the research presented in Table~\ref{table4} 
are significant and provide crucial insights. This study highlights 
that an increase in the values of the rabies model parameters, 
including $\psi_{1}$, $\psi_{2}$, $\kappa_{1}$, and $\kappa_{2}$, 
results in a proportional rise in the magnitude of $\mathcal{R}_{0}$. 
This indicates that careful consideration of these parameters is essential 
in predicting the transmission dynamics of the disease.
Furthermore, the study also reveals that an increase in the values of 
parameters such as $\gamma_{1}$, $\gamma_{2}$, $\mu_{2}$, $\mu_{3}$, 
$\rho_{1}$, $\rho_{2}$, $\sigma_{2}$, and $\sigma_{3}$ leads to a decrease 
in the magnitude of $\mathcal{R}_{0}$. This implies that controlling 
the spread of the disease can be achieved by adjusting these parameters.
It is noteworthy that a 20\% increase in any of the parameters $\psi_{1}$, 
$\psi_{2}$, $\kappa_{1}$, or $\kappa_{2}$ corresponds exactly to 
a 20\% increase in the value of $\mathcal{R}_{0}$. 

 
\subsection{Existence of the steady state solution}
\label{sec:EEp}

The endemic equilibrium point is the steady state where rabies is present in humans, 
free-range dogs, and domestic dogs. To find this point, we set the equations of the model system 
$\left(\ref{eqn1}\right)$ to zero and solve the resulting system simultaneously. 
The state variables for each compartment are represented by
\begin{equation*}
{\mathbb E}\left(S_{H}^{*},\; E_{H}^{*},\; I_{H}^{*},\; 
R_{H}^{*},\; S_{F}^{*},\; E_{F}^{*},\; I_{F}^{*},\; S_{D}^{*},\; 
E_{D}^{*},\; I_{D}^{*},\;  R_{D}^{*},\; M^{*}\right).
\end{equation*}    

Solving the first and the second equations in the model system (\ref{eqn1}) we have
\begin{equation}
\begin{aligned}
E^{*}_{{H}}={\dfrac {R_{{H}}\beta_{{3}}-S_{{H}}\mu_{{1}}+\theta_{{1}}}{\mu_
{{1}}+\beta_{{1}}+\beta_{{2}}}}.
\end{aligned}
\label{eqn2}
\end{equation}
By substituting \eqref{eqn2} into the third 
and fourth equations of the model system (\ref{eqn1}), 
and solving for $I_{H}$ and $R_{H}$, we obtain the following results: 
\begin{equation}
\begin{aligned}
I^{*}_{{H}} &= {\dfrac {\beta_{{1}} \left( R_{{H}}\beta_{{3}}
-S_{{H}}\mu_{{1}}+\theta_{{1}} \right) }{ \left( \mu_{{1}}+\beta_{{1}}
+\beta_{{2}} \right) \left( \sigma_{{1}}+\mu_{{1}} \right) }}, \\
R^{*}_{{H}} &= {\dfrac {\beta_{{2}} \left( -S_{{H}}\mu_{{1}}
+\theta_{{1}} \right) }{ \beta_{{1}}\beta_{{3}}+\beta_{{1}}\mu_{{1}}
+\beta_{{2}}\mu_{{1}}+\mu_{{1}}\beta_{{3}}+{\mu_{{1}}}^{2}}}.
\end{aligned}
\label{eqn3}
\end{equation}
Upon substitution of $R^{*}_{{H}}$ of \eqref{eqn3} into \eqref{eqn2}, 
the following results are obtained:
\begin{equation}
\begin{aligned}
E^{*}_{{H}}={\dfrac { \left( S_{{H}}\mu_{{1}}+\beta_{{1}}
+\beta_{{2}}+\mu_{{1}}-\theta_{{1}} \right)  
\left( \mu_{{1}}+\beta_{{3}} \right) }{\beta_{{3}}\beta_{{2}}}},\;\;
I^{*}_{{H}}={\dfrac {\beta_{{1}} \left( \beta_{{3}}+\mu_{{3}} \right) 
\left( S_{{H}}\mu_{{1}}-\theta_{{1}} \right) }{ \left( \sigma_{{1}}
+\mu_{{1}} \right) ^{2} \left(  \left( \beta_{{1}}+\beta_{{2}}
+\beta_{{3}} \right) \mu_{{1}}+\beta_{{1}}\beta_{{3}} \right)}}.
\end{aligned}
\label{eqn4}
\end{equation}
From $\left(\ref{eqn2}\right)$, $\left(\ref{eqn3}\right)$ 
and $\left(\ref{eqn4}\right)$, it can be derived that
\begin{equation*}
\left.
\begin{aligned}
R^{*}_H &= \dfrac{\beta_2(\theta_1(\beta_3 + \mu_3)(\sigma_1 
+ \mu_1)^2 + (\beta_3 + \mu_3)(\sigma_1 + \mu_1)^2(\beta_1 + \beta_2 
+ \beta_3))}{\mu_1(\beta_3 + \mu_3)(\sigma_1 + \mu_1)^2 
+ \mu_1^2(\beta_1 + \beta_2 + \beta_3) + \mu_1\beta_1\beta_3},\\\\
I^{*}_{H} &= \dfrac{\beta_1(\beta_3 + \mu_3)(\sigma_1 + \mu_1)^2(\beta_1 
+ \beta_2 + \beta_3)\mu_1 + \beta_1\beta_3(\sigma_1 + \mu_1)^2}{(\sigma_1 
+ \mu_1)^2((\beta_1 + \beta_2 + \beta_3)\mu_1 + \beta_1\beta_3)} \\
&\quad - \dfrac{\beta_1(\beta_3 + \mu_3)(\sigma_1 + \mu_1)^2\beta_3 
- \theta_1(\beta_3 + \mu_3)(\sigma_1 + \mu_1)^2}{(\sigma_1 
+ \mu_1)^2((\beta_1 + \beta_2 + \beta_3)\mu_1 + \beta_1\beta_3)},\\\\
E^{*}_H &= \dfrac{(\mu_1 + \beta_3)(\theta_1(\beta_3 + \mu_3)(\sigma_1 
+ \mu_1)^2 + (\beta_3 + \mu_3)(\sigma_1 + \mu_1)^2(\beta_1 + \beta_2 
+ \beta_3) }{\mu_1^2 + (\beta_1 + \beta_2 + \beta_3)\mu_1 + \beta_1\beta_3}\\
&\quad + \dfrac{\beta_3(\sigma_1 + \mu_1)^2 - (\beta_3 + \mu_3)(\sigma_1 
+ \mu_1)^2\beta_3}{\mu_1^2 + (\beta_1 + \beta_2 + \beta_3)\mu_1 + \beta_1\beta_3},\\
S^{*}_{{H}} &={\dfrac {-{\mu_{{1}}}^{2}+ \left( -\beta_{{1}}-\beta_{{2}}-
\beta_{{3}}+\theta_{{1}} \right) \mu_{{1}}+ \left(  \left( E^{*}_{{H}}-1
\right) \beta_{{2}}-\beta_{{1}}+\theta_{{1}} \right) \beta_{{3}}}{\mu_{{1}} 
\left( \mu_{{1}}+\beta_{{3}} \right) }}. 
\end{aligned}
\right\}
\end{equation*}
Rewriting  $\left(\ref{eqn5}\right)$ as  
\begin{equation*}
\begin{aligned}
{R_{{0}}}^{2}- \left( R_{{0}}-1 \right)  \left( R_{{33}}-R_{{55}}
\right) +R_{{33}}+R_{{55}}+R_{{35}}R_{{53}},
\end{aligned}
\end{equation*}
where $R_{33}$, $R_{55}$, $R_{35}$ and $R_{53}$ 
are defined in $\left(\ref{eqn6}\right)$,
\begin{equation*}
\left.
\begin{aligned}
I^{*}_{D} & = \dfrac{(\mu_3 + \gamma_3) \left( 
\frac{{R_0}^2 - (R_0 - 1) \frac{\kappa_1 \theta_2 \gamma}{\mu_2 (\mu_2 
+ \gamma) (\mu_2 + \sigma_2)} + 1 + \frac{\kappa_1 \theta_2 \gamma}{\mu_2 
(\mu_2 + \gamma) (\mu_2 + \sigma_2)}}{Q_1 + Q_2} \right) \gamma_1}{\mu_3^3 
+ (\gamma_1 + \gamma_2 + \gamma_3 + \sigma_3 + 1) \mu_3^2 \ + Q_{3}},\\
E^{*}_{{D}} &=\dfrac{(\mu_{{3}}+\sigma_{{3}}) \left(
(R_{0}^2-(R_{0}-1)\frac{\kappa_{{1}}\theta_{{2}}\gamma}{\mu_{{2}}(\mu_{{2}}
+\gamma)(\mu_{{2}}+\sigma_{{2}})}+1+\frac{\kappa_{{1}}\theta_{{2}}\gamma}{\mu_{{2}}(\mu_{{2}}
+\gamma)(\mu_{{2}}+\sigma_{{2}})}) (\mu_{{3}}+\gamma_{{3}})\right) }{{\mu_{{3}}}^{3}
+ (\gamma_{{1}}+\gamma_{{2}}+\gamma_{{3}}+\sigma_{{3}}+1) {\mu_{{3}}}^{2}
+ ({\gamma_{{3}}}^{2}+ (1+\sigma_{{3}}) \gamma_{{3}}+ (1+\sigma_{{3}}) 
(\gamma_{{1}}+\gamma_{{2}})) \mu_{{3}}+\gamma_{{3}} (\gamma_{{3}} \sigma_{{3}}
+\gamma_{{1}})},\\
R^{*}_{{D}}&={\frac {\left(
(R_{0}^2-(R_{0}-1)\dfrac{\kappa_{{1}}\theta_{{2}}\gamma}{\mu_{{2}}(\mu_{{2}}
+\gamma)(\mu_{{2}}+\sigma_{{2}})}+1+\frac{\kappa_{{1}}\theta_{{2}}\gamma}{\mu_{{2}}(\mu_{{2}}
+\gamma)(\mu_{{2}}+\sigma_{{2}})}) (\mu_{{3}}+\sigma_{{3}}) \theta_{{3}} \right) 
\gamma_{{2}}}{{\mu_{{3}}}^{3}+ \left( \gamma_{{1}}+\gamma_{{2}}+\gamma_{{3}}
+\sigma_{{3}}+1 \right) {\mu_{{3}}}^{2}+Q_{4}
+\gamma_{{3}} \left( \gamma_{{3}}\sigma_{{3}}+\gamma_{{1}} \right)}},
\end{aligned}
\right\}
\end{equation*}
\begin{equation*}
\left.
\begin{aligned}
S^{*}_{{D}}&={\frac {-E^{*}_{{D}}{\mu_{{3}}}^{2}+ \left(  \left( -\gamma_{{1}}
-\gamma_{{2}}-\gamma_{{3}} \right) E^{*}_{{D}}+\theta_{{3}} \right) 
\mu_{{3}}-\gamma_{{3}} \left( \gamma_{{1}}E^{*}_{{D}}-\theta_{{3}} \right) }{
\mu_{{3}} \left( \mu_{{3}}+\gamma_{{3}} \right) }},\\
I^{*}_{{F}}&={\frac {\gamma\,E^{*}_{{F}}}{\mu_{{2}}+\sigma_{{2}}}},\;\;
S^{*}_{F}={\dfrac { \left( -\gamma-\mu_{{2}} \right) E^{*}_{{F}}+\theta_{{2}}}{
\mu_{{2}}}},\;\;
M^{*}={\dfrac {\nu_{{3}}I^{*}_{{D}}+\nu_{{2}}I^{*}_{{F}}+\nu_{{1}}I^{*}_{{H}}}{
\mu_{{4}}}},
\end{aligned}
\right\}
\end{equation*}
\begin{equation*}
\left.
\begin{aligned}
Q_{{1}} &= \frac{\kappa_{{1}}\theta_{{2}}\gamma_{{1}}\psi_{{1}}\gamma}{
\mu_{{2}} \left( \mu_{{3}}+\gamma_{{1}}+\gamma_{{2}} \right)  
\left( \mu_{{3}}+\sigma_{{3}} \right) \mu_{{3}} \left( 1+\rho_{{1}} \right)  
\left( \mu_{{2}}+\gamma \right)  \left( \mu_{{2}}+\sigma_{{2}} \right)} \\
&\quad - \frac{\psi_{{2}}\gamma_{{1}}R_{{0}}}{\left( 1+\rho_{{2}} \right)  
\left( \mu_{{3}}+\gamma_{{1}}+\gamma_{{2}} \right)  \left( \mu_{{3}}
+\sigma_{{3}} \right)}, \\
Q_{{2}} &= \frac{\psi_{{2}}\gamma_{{1}}}{\left( 1+\rho_{{2}} \right)  
\left( \mu_{{3}}+\gamma_{{1}}+\gamma_{{2}} \right)  \left( \mu_{{3}}
+\sigma_{{3}} \right)}, \\
Q_{{3}} &= (\gamma_3^2 + (1 + \sigma_3) \gamma_3 + (1 + \sigma_3)(\gamma_1 
+ \gamma_2)) \mu_3 + \gamma_3 (\gamma_3 \sigma_3 + \gamma_1),\\
Q_{{4}} &= \left( {\gamma_{{3}}}^{
2}+ \left( 1+\sigma_{{3}} \right) \gamma_{{3}}+ \left( 1+\sigma_{{3}}
 \right)  \left( \gamma_{{1}}+\gamma_{{2}} \right)  \right) \mu_{{3}}.
\end{aligned}
 \right\}
\end{equation*}
The endemic equilibrium point of the rabies disease persists if 
$E_{H},\; E_{F},\; E_{D} > 0$ and ${\cal R}_0 \geq 1$, 
as summarized in Theorem~\ref{The}.

\begin{theorem}
The system model $\left(\ref{eqn1}\right)$ has a unique endemic equilibrium 
$\mathbb{E}^*$ if $\mathcal{R}_0 \geq 1$ and $E_{H},\; E_{F},\; E_{D} > 0$.
\label{The}
\end{theorem}


\section{Stability Analysis}
\label{sec:3}

We begin by studying the local stability of the disease free equilibrium (DFE) point
(Section~\ref{subsec:LS:DFE}). Then, we investigate 
its global stability (Section~\ref{subsec:GS:DFE}).
The stability of the endemic equilibrium is given in Appendix~A.


\subsection{Local stability of the disease free equilibrium point} 
\label{subsec:LS:DFE}

We prove local stability of the disease free equilibrium point
with the help of the Routh--Hurwitz criterion.

\begin{theorem}
\textup{The DFE point} ${\cal E}_{0}$ \textup{is locally asymptotically 
stable if} ${\cal R}_{0} < 1$ \textup{and all eigenvalues of the Jacobian matrix} 
$\left(J\left({\cal E}_{0}\right)\right)$ \textup{evaluated at} 
${\cal E}_{0}$ \textup{have negative real parts.}
\label{The2}
\end{theorem}

\begin{proof}
We need to show that all the eigenvalues of the matrix $J\left({\cal E}_{0}\right)$ 
in equation $\left(\ref{Jacob4}\right)$ of the model system $\left(\ref{eqn1}\right)$ 
at the DFE point have negative real parts. Subsequently, since the endemic equilibrium 
exists if, and only if, ${\cal R}_{0} < 1$, we utilize the Jacobian matrix at 
the disease-free state $J\left({\cal E}_{0}\right)$, which is expressed as
\begin{equation}
\setlength{\arraycolsep}{1.5pt}
J\left({\cal E}_0\right)
=\left( 
\begin{array}{cccccccccccc} 
-\mu_{{1}}&0&0&0&0&0&-{\dfrac {
\tau_{{1}}\theta_{{1}}}{\mu_{{1}}}}&0&0&-{\frac {\tau_{{2}}\theta_{{1}}}{
\mu_{{1}}}}&0&0\cr \noalign{\medskip}0&-a_{1} 
&0&0&0&0&{\dfrac {\tau_{{1}}\theta_{{1}}}{\mu_{{1}}}}&0&0&{\dfrac {
\tau_{{2}}\theta_{{1}}}{\mu_{{1}}}}&0&0\cr \noalign{\medskip}0&\beta_{{
1}}&-a_{2}&0&0&0&0&0&0&0&0&0\cr \noalign{\medskip}0&
\beta_{{2}}&0& -a_{3}&0&0&0&0&0&0&0&0
\cr \noalign{\medskip}0&0&0&0&-\mu_{{2}}&0&-{\dfrac {\kappa_{{1}}\theta_
{{2}}}{\mu_{{2}}}}&0&0&-{\frac {\kappa_{{2}}\theta_{{2}}}{\mu_{{2}}}}&0
&0\cr \noalign{\medskip}0&0&0&0&0& -a_{4}&{\dfrac {\kappa_{{1}}
\theta_{{2}}}{\mu_{{2}}}}&0&0&{\dfrac {\kappa_{{2}}\theta_{{2}}}{\mu_{{2}}}}
&0&0\cr \noalign{\medskip}0&0&0&0&0&\gamma&-a_{5}&0&0&0&0&0\cr 
\noalign{\medskip}0&0&0&0&0&0&-{\dfrac {\psi_{{1}}\theta_
{{3}}}{\mu_{{3}} \left( 1+\rho_{{1}} \right) }}&-\mu_{{3}}&0&-{\dfrac {
\psi_{{2}}\theta_{{3}}}{\mu_{{3}} \left( 1+\rho_{{2}} \right) }}&
\gamma_{{3}}&0\\ \noalign{\medskip}0&0&0&0&0&0&{\dfrac {\psi_{{1}}
\theta_{{3}}}{\mu_{{3}} \left( 1+\rho_{{1}} \right) }}&0& -a_{6}
&{\dfrac {\psi_{{2}}\theta_{{3}}}{\mu_{{3}}
\left( 1+\rho_{{1}} \right) }}&0&0\cr \noalign{\medskip}0&0&0&0&0&0&0&0
&\gamma&-a_{7}&0&0\cr \noalign{\medskip}0&0&0&0&0&0&0&0
&\gamma_{{2}}&0& -a_{8}&0\cr \noalign{\medskip}0&0&\nu_{{1}}
&0&0&0&\nu_{{2}}&0&0&\nu_{{3}}&0&-\mu_{{4}}
\end {array} 
\right),
\label{Jacob4}
\end{equation} 
where
\begin{align*}
a_{1} = \mu_{1} + \beta_{1} + \beta_{2},\;
a_{2} = \mu_{1} + \sigma_{1},\;
a_{3} = \mu_{3} + \beta_{3},\;
a_{4} = \mu_{2} + \gamma,\\
a_{5} = \mu_{2} + \sigma_{2},\;
a_{6} = \mu_{3} + \gamma_{1} + \gamma_{2},\;
a_{7} = \mu_{3} + \sigma_{3},\;
a_{8} = \mu_{3} + \gamma_{3}.
\end{align*}
The first, fourth, fifth, eighth, and twelfth columns of the matrix 
$\left( J({\cal E}0) \right)$ in equation $\left(\ref{Jacob4}\right)$ 
contain the diagonal terms. It is obvious from the eigenvalues 
$\lambda_{1}=-\mu_{1}$, $\lambda_{2}=-a_3$, 
$\lambda_{3}=-\mu_{2}$, $\lambda_{4}=-\mu_{3}$, 
and $\lambda_{5}=-\mu_{4}$, respectively. Thus, the matrix 
$\left( J({\cal E}_0) \right)$ reduces to
\begin{equation}
\setlength{\arraycolsep}{1.5pt}
J\left({\cal E}_0\right)=
\left(\begin {array}{ccccccc} -a_{{1}}&0&0&{\frac {\tau_{{1}}
\theta_{{1}}}{\mu_{{1}}}}&0&{\frac {\tau_{{2}}\theta_{{1}}}{\mu_{{1}}}}&0\\ 
\noalign{\medskip}\beta_{{1}}&-a_{{2}}&0&0&0&0&0\\ 
\noalign{\medskip}0&0&-a_{{4}}&{\frac {\kappa_{{1}}\theta_{{2}}}{
\mu_{{2}}}}&0&{\frac {\kappa_{{2}}\theta_{{2}}}{\mu_{{2}}}}&0\\ 
\noalign{\medskip}0&0&\gamma&-a_{{5}}&0&0&0\\ \noalign{\medskip}0&0
&0&{\frac {\psi_{{1}}\theta_{{3}}}{\mu_{{3}} \left( 1+\rho_{{1}}
\right) }}&-a_{{6}}&{\frac {\psi_{{2}}\theta_{{3}}}{\mu_{{3}} \left( 
1+\rho_{{1}} \right) }}&0\\ \noalign{\medskip}0&0&0&0&\gamma&-a_{{7}}&0\\ 
\noalign{\medskip}0&0&0&0&\gamma_{{2}}&0&-a_{{8}} 
\end {array} 
\right).
\label{Jacob5}
\end{equation} 
Again, the second and seventh columns of the matrix $\left( J({\cal E}0) \right)$ 
in $\left(\ref{Jacob5}\right)$ contain the diagonal terms. It is obvious 
from the eigenvalues $\lambda_{6}=-a_2$ and $\lambda_{7}=-a_8$. Thus, 
the matrix $\left( J({\cal E}_0) \right)$ reduces to
\begin{equation}
\setlength{\arraycolsep}{1.5pt}
J\left({\cal E}_0\right)=
\left(\begin {array}{ccccccc} -a_{{1}}&0
&{\frac {\tau_{{1}}\theta_{{1}}}{\mu_{{1}}}}&0&{\frac {\tau_{{2}}\theta_{{1}}}{\mu_{{1}}}}\\ 
\noalign{\medskip}0&-a_{{4}}&{\frac {\kappa_{{1}}\theta_{{2}}}{\mu_{{2}}}}
&0&{\frac {\kappa_{{2}}\theta_{{2}}}{\mu_{{2}}}}\\ 
\noalign{\medskip}0&\gamma&-a_{{5}}&0&0\\ \noalign{\medskip}0
&0&{\frac {\psi_{{1}}\theta_{{3}}}{\mu_{{3}} \left( 1+\rho_{{1}} \right) }}
&-a_{{6}}&{\frac {\psi_{{2}}\theta_{{3}}}{\mu_{{3}} \left( 1+\rho_{{1}} 
\right) }}\\ \noalign{\medskip}0&0&0&\gamma&-a_{{7}}
\end {array} 
\right).
\label{Jacob6}
\end{equation} 
Again, the  first column of the matrix $\left(  J({\cal E}_0) \right)$ 
in equation $\left(\ref{Jacob6}\right)$ contains  the diagonal term, 
and it is obvious from  eigenvalues $\lambda_{8}=-a_1$. Thus,  
the matrix $\left(  J({\cal E}_0) \right)$ reduces to
\begin{equation}
\setlength{\arraycolsep}{1.5pt}
J\left({\cal E}_0\right)=
\left(\begin {array}{ccccccc} -a_{{4}}
&{\frac {\kappa_{{1}}\theta_{{2}}}{\mu_{{2}}}}&0&{\frac {\kappa_{{2}}\theta_{{2}}}{\mu_{{2}}}}\\ 
\noalign{\medskip}\gamma&-a_{{5}}&0&0\\ \noalign{\medskip}0
&{\frac {\psi_{{1}}\theta_{{3}}}{\mu_{{3}} \left( 1+\rho_{{1}} \right) }}
&-a_{{6}}&{\frac {\psi_{{2}}\theta_{{3}}}{\mu_{{3}} \left( 1+\rho_{{1}}\right) }}\\ 
\noalign{\medskip}0&0&\gamma&-a_{{7}}
\end {array} \right).
\label{Jacob7}
\end{equation}
Computing the eigenvalues of the given matrix $J({\cal E}_0)$ 
in equation $\eqref{Jacob7}$ involves solving the characteristic polynomial equation
\begin{equation*}
P\left(\lambda\right)=\det(J({\cal E}_0) - \lambda I)=0, 
\end{equation*}
where \(I\) is the identity matrix and \(\lambda\) represents the eigenvalues. Thus, 
\begin{equation}
{\lambda}^{4}+C_{{1}}{\lambda}^{3}+C_{{2}}{\lambda}^{2}+C_{{3}}\lambda +C =0.
\label{Poly1}
\end{equation} 
Equations $\left(\ref{Jacob4}\right)$, $\left(\ref{Jacob5}\right)$, 
and $\left(\ref{Jacob6}\right)$ evidence that $\lambda_1$, $\lambda_2$, $\lambda_3$, 
$\lambda_4$, $\lambda_5$, $\lambda_6$, $\lambda_7$, and $\lambda_8$ exhibit 
negative real parts. By applying the Routh--Hurwitz criterion, the other four 
eigenvalues of the matrix  equation $\left(\ref{Jacob7}\right)$  will also have 
negative real parts if all coefficients in equation $\left(\ref{Poly1}\right)$ 
are greater than zero. Then,
\begin{align*}
\left\{
\begin{aligned}
C_{1} &= a_{7} + a_{6} + a_{5} + a_{4}>0, \\
C_{2} &= -\gamma\,b_{1} - \gamma\,b_{4} + a_{5}a_{4} 
+ a_{6}a_{4} + a_{7}a_{4} + a_{6}a_{5} + a_{7}a_{5} + a_{7}a_{6}>0,\\
C_{3} &= \left( (a_{6} + a_{7})a_{5} - \gamma\,b_{4} + a_{7}a_{6} \right) 
a_{4} + \left( -\gamma\,b_{4} + a_{7}a_{6} \right) a_{5} - \gamma\,b_{1} (a_{6} + a_{7})>0,\\ 
C &= \gamma^{2}b_{1}b_{4} - \gamma^{2}b_{2}b_{3} - \gamma a_{4}a_{5}b_{4} 
- \gamma a_{6}a_{7}b_{1} + a_{4}a_{5}a_{6}a_{7}>0,
\end{aligned}
\right.
\end{align*}
where
\begin{align*}
\left\{
\begin{aligned}
b_{1} &= \dfrac{\kappa_{1}\theta_{2}}{\mu_{2}}, \;\;
b_{2} &= \dfrac{\psi_{1}\theta_{3}}{\mu_{3} (1+\rho_{1})}, \\
b_{3} &= \dfrac{\kappa_{2}\theta_{2}}{\mu_{2}}, \;\;
b_{4} &= \dfrac{\psi_{2}\theta_{3}}{\mu_{3} (1+\rho_{1})}.
\end{aligned}
\right.
\end{align*}
Since  \textup{all eigenvalues of the Jacobian matrix} $\left(J\left({\cal E}_{0}\right)\right)$ 
\textup{evaluated at} ${\cal E}_{0}$ \textup{have negative real parts, the model system  
$\left(\ref{eqn1}\right)$ at the ${\cal E}_0$  
\textup{is locally asymptotically stable if} ${\cal R}_{0} < 1$.}
\end{proof}
     

\subsection{Global stability of the DFE point}
\label{subsec:GS:DFE}

We prove the global stability of the DFE point ${\cal E}_{0}$  
of the rabies model \eqref{eqn1} using the theorem described by 
\cite{kamgang2008computation}. To apply the theorem, 
we write the model system $\left(\ref{eqn1}\right)$ as
\begin{equation*}
\begin{aligned}
\left.
\begin{array}{llll}
\dfrac{dY_{s}}{dt} &= G_{0}\left(Y_{s}-Y\left({\cal E}_{0}\right)\right)+G_{1}Y_{i}, \\\\
\dfrac{dY_{i}}{dt} &= G_{2}Y_{i},
\end{array}
\right\}
\end{aligned}
\end{equation*}
where $Y_{s}$ is the vector representing the compartments that do not 
transmit the rabies disease, and $Y_{i}$ symbolizes the rabies-transmitting 
vector compartments. In the case of $G_{2}$, if $G_{2}$ is a Metzler matrix 
(i.e., the off-diagonal entries of $G_{2}$ are non-negative), and $G_{0}$ 
has real negative eigenvalues, the rabies-free equilibrium is globally 
asymptotically stable. Based on the model system 
$\left(\ref{eqn1}\right)$, we have 
$Y_{s}=\left(S_{H} \;,R_{H} \;,S_{F} \;,S_{D} \;, R_{D} \right)^{T}$,  
$Y_{i}=\left(E_{H} \;,I_{H} \;,E_{F} \;,I_{F} \;, E_{D}\;,I_{D} \;,M \right)^{T}$ 
and 
\begin{equation*}
Y_{s}-Y\left({\cal E}_{0}\right)
= \left(
\begin{array}{c}
S_{H}-\dfrac{\theta_{1}}{\mu_{1}}\cr
R_{H}\cr
S_{F}-\dfrac{\theta_{2}}{\mu_{2}}\cr
S_{D}-\dfrac{\theta_{3}}{\mu_{3}}\cr
R_{D} \cr
\end{array}
\right),
\end{equation*}
\begin{equation*}
G_{0}=
\left(
\begin{array}{ccccc}
-\mu & \beta_{3}&0&0&0\cr
0&-\left(\beta_{3}+\mu_{1}\right)&0&0&0\cr
0&0&-\mu_{2}&0&0\cr
0&0&0&-\mu_{3}&\gamma_{3}\cr
0&0&0&0&-\left(\mu_{3}+\gamma_{3}\right)
\end{array}
\right).
\end{equation*}
The eigenvalues of the matrix $G_{0}$ are $\lambda_{1}=\mu_{3}$,  
$\lambda_{2}=\mu_{2}$,  $\lambda_{3}=\mu_{1}$, 
$\lambda_{4}=-\left(\mu_{3}+\gamma_{3}\right)$, 
$\lambda_{5}=-\left(\beta_{3}+\mu_{1}\right)$, while
\begin{equation*}
G_{1}=
\left(
\begin {array}{ccccccc} 
0&0&0&{\dfrac {\tau_{{1}}\theta_{{1}}}{
\mu_{{1}}}}&0&{\dfrac {\tau_{{2}}\theta_{{1}}}{\mu_{{1}}}}&0\\ 
\noalign{\medskip}\beta_{{2}}&0&0&0&0&0&0\cr \noalign{\medskip}0&0&0
&{\frac {\kappa_{{1}}\theta_{{2}}}{\mu_{{2}}}}&0&{\dfrac {\kappa_{{2}}
\theta_{{2}}}{\mu_{{2}}}}&0\cr \noalign{\medskip}0&0&0&{\dfrac {\psi_{{1}}
\theta_{{3}}}{\mu_{{3}} \left( 1+\rho_{{1}} \right) }}&0
&{\dfrac {\psi_{{2}}\theta_{{3}}}{\mu_{{3}} \left( 1+\rho_{{2}} \right) }}
&0\cr \noalign{\medskip}0&0&0&0&\gamma_{{2}}&0&0
\end{array}
\right),
\end{equation*}
\begin{equation*}
G_{2} =	
\left(
\begin {array}{ccccccc} 	
-\mu_{{1}}-\beta_{{1}}-\beta_{{2}}&0&0
&{\frac {\tau_{{1}}\theta_{{1}}}{\mu_{{1}}}}&0&{\frac {\tau_{{2}}
\theta_{{1}}}{\mu_{{1}}}}&0\\ \noalign{\medskip}\beta_{{1}}
&-\sigma_{{1}}-\mu_{{1}}&0&0&0&0&0\\ \noalign{\medskip}0&0
&-\mu_{{2}}-\gamma
&{\frac {\kappa_{{1}}\theta_{{2}}}{\mu_{{2}}}}&0
&{\frac {\kappa_{{1}}\theta_{{2}}}{\mu_{{2}}}}&0\\ 
\noalign{\medskip}0&0&\gamma&-\mu_{{2}}-
\sigma_{{2}}&0&0&0\\ \noalign{\medskip}0&0&0&{\frac {\psi_{{1}}
\theta_{{3}}}{\mu_{{3}} \left( 1+\rho_{{1}} \right) }}
&-\mu_{{3}}-\gamma_{{1}}-\gamma_{{2}}
&{\frac {\psi_{{2}}\theta_{{3}}}{\mu_{{3}} \left( 
1+\rho_{{1}} \right) }}&0\\ \noalign{\medskip}0&0&0
&0&\gamma&-\mu_{{3}}-\sigma_{{3}}&0\\ \noalign{\medskip}0
&\nu_{{1}}&0& \nu_{{2}}&0&\nu_{{3}}&-\mu_{{4}}
\end{array} \right). 
\end{equation*}
Since the eigenvalues of the $G_{0}$ are negative and the off diagonal 
of the Metzler matrix $G_{2}$ are non-negative, then the rabies 
DFE point is globally asymptotically stable. 


\section{Model fitting and parameter estimation} 
\label{sec:model:fit:parm:est}

After conducting model analysis of the dynamics and qualitative outcomes of the rabies model, 
it is essential to accurately determine the model's parameters for making quantitative 
predictions within a limited time frame using real-world data \citep{li2018introduction}. 
In this study, we employed the non-linear least squares method (NLSM) 
to estimate the parameters of model equation $\left(\ref{eqn1}\right)$. To achieve this, 
we generated synthetic data that represented the expected disease spread patterns 
at various time points, denoted as $t_i$ \citep{myung2003tutorial}. These patterns 
were computed by numerically solving equation  $\left(\ref{eqn1}\right)$ 
with a fifth-order Runge--Kutta method in the \textsf{MATLAB} environment, 
initializing the parameters  with values  from  literature  denoted as 
$\Theta_i$ and  initial condition for the number of   
$S_H\left(0\right) = 142000 $, $E_H\left(0\right) = 40$, 
$I_H\left(0\right) = 0$, $R_H\left(0\right) = 0$,  
$S_D\left(0\right) = 15000$, $E_D\left(0\right) = 25$,
$I_D\left(0\right) = 0$, $R_D\left(0\right) = 0$, 
$S_F\left(0\right) = 12500$, $E_F\left(0\right) = 20$, 
$I_F\left(0\right) = 0$, and $M\left(0\right) = 90$. 
In order to generate the  rabies  dataset $ RD\left(t_i\;\;\Theta_i\right)$ 
we added random Gaussian noise $\eta_i\left(t_i\;\;\Theta_i\right)$ 
measurements to the data, simulating real-world dynamics where 
measurement errors are common. Thus the  observed/actual  
dependent data were given as
\begin{equation*}
\begin{array}{llll}
Y_{i}=RD\left(t_i\;\;\Theta_i
\right)+\eta_i\left(t_i\;\;\Theta_i\right) \;\; 
\text{for each time}\;\; t_i\in[1,\;\;n].
\end{array}
\end{equation*}
 
\begin{figure}[!ht]
\begin{minipage}[b]{0.45\textwidth}
\includegraphics[height=5.0cm, width=7.5cm]{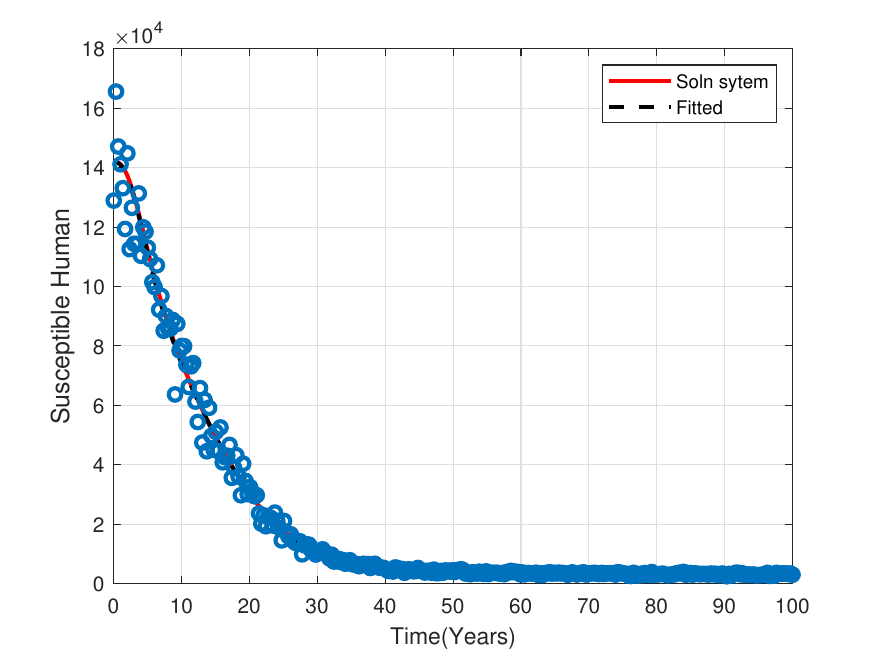}
\centering{(a)}
\end{minipage}
\begin{minipage}[b]{0.45\textwidth}
\includegraphics[height=5.0cm, width=7.5cm]{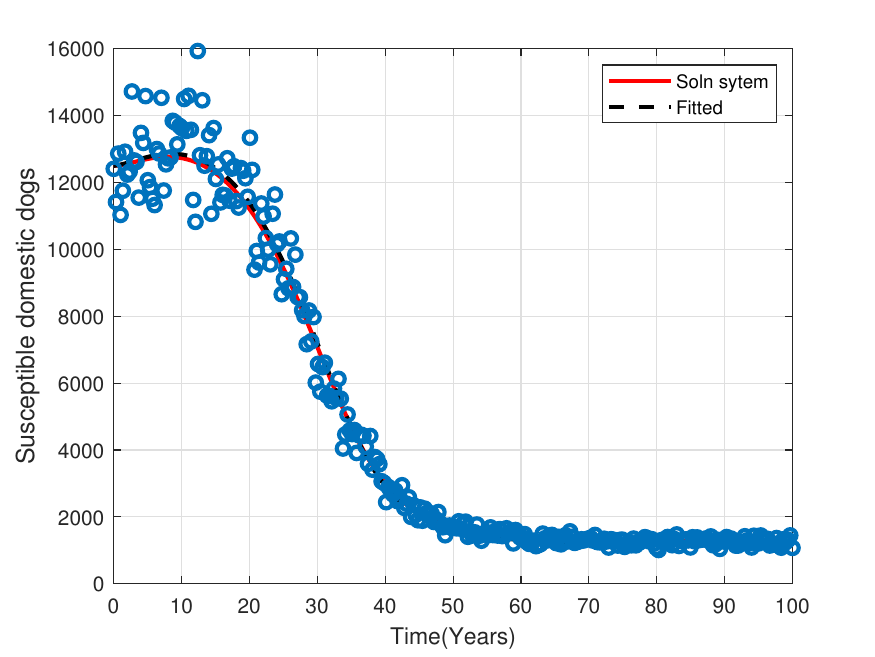}
\centering{(b)}
\end{minipage}
\begin{minipage}[b]{0.45\textwidth}
\includegraphics[height=5.0cm, width=7.5cm]{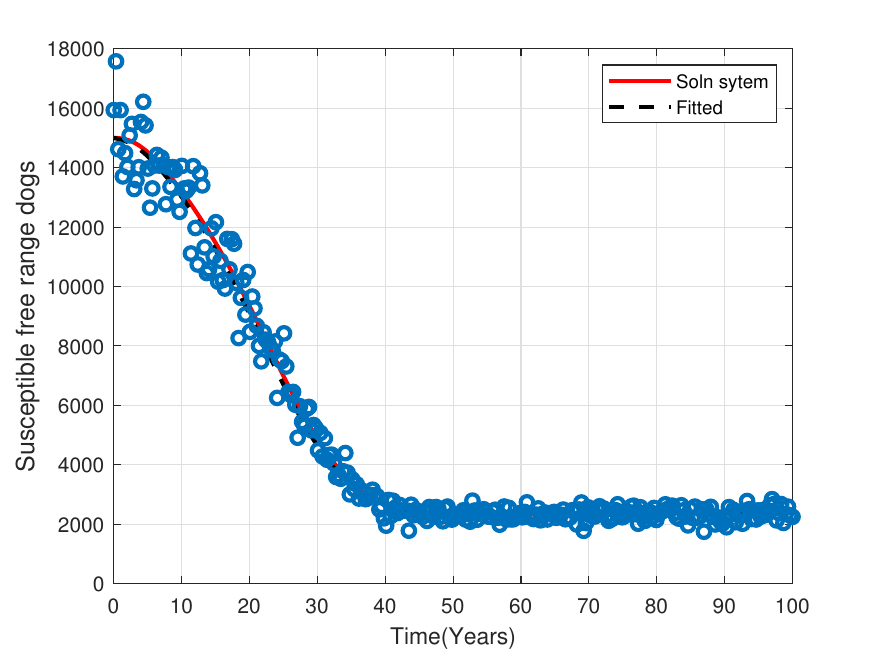}
\centering{(c)}
\end{minipage}
\begin{minipage}[b]{0.45\textwidth}
\includegraphics[height=5.0cm, width=7.5cm]{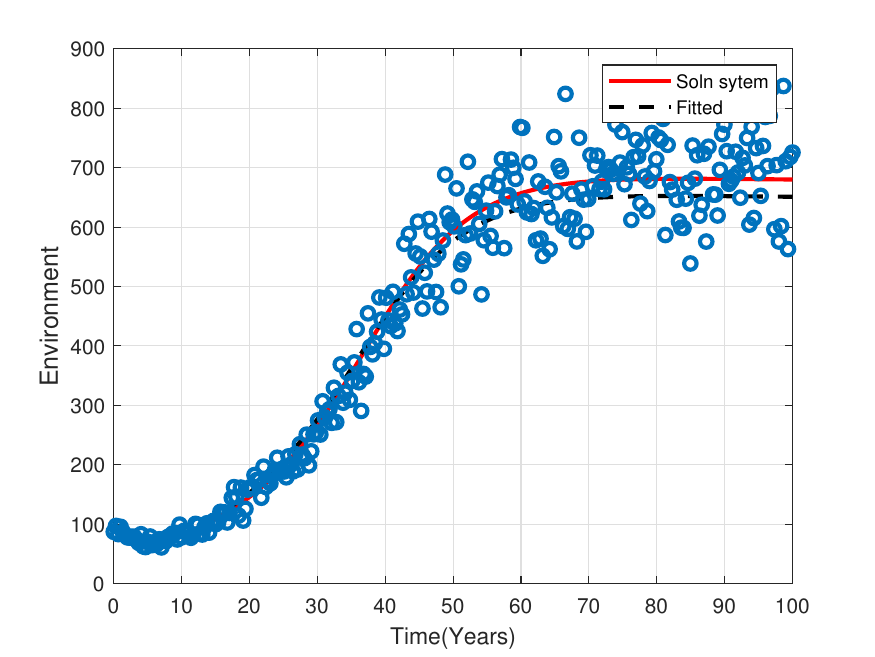}
\centering{(d)}
\end{minipage}
\centering
\caption{Scatter estimated with standard deviation of 0.05 
and numerical simulation (solid) with confidence interval of 95\%.}
\label{Fig9}
\end{figure}

The parameter values $\overline{YY}$ of Table~\ref{table3} 
were determined by minimizing the sum of squared residuals  expressed as 
\begin{equation*}
\begin{array}{ll}
\displaystyle\overline{YY}(\Theta) = \min \sum_{k=1}^{n} (Y_i - Y)^2
\end{array}
\end{equation*}
between the model solutions $(Y)$ obtained through solving the rabies 
$\left(\ref{eqn1}\right)$ model using the real parameters from 
the generated data and the synthetic data $(Y_i)$ generated by 
introducing random Gaussian noise to the model output $RD\left(t_i\;\;\Theta_i\right)$ 
\citep{li2018introduction}. The estimated parameter values were then used to fit 
the data $(Y_{i})$, and the resulting best fits are depicted in Figure~\ref{Fig9}(a)--(d) 
and the resulting  estimated parameters given in Table~\ref{table3}.


\section{Numerical Simulations}
\label{sec:4}

In this section, we employed the \texttt{ode45} method available in \textsf{MATLAB} 
software to numerically solve a model system (\ref{eqn1}) using  parameters  
presented  in Table~\ref{table3} along with the initial conditions 
$S_H\left(0\right) = 142000 $, $E_H\left(0\right) = 40$, $I_H\left(0\right) = 0$, 
$R_H\left(0\right) = 0$, $S_D\left(0\right) = 15000$, $E_D\left(0\right) = 25$,
$I_D\left(0\right) = 0$, $R_D\left(0\right) = 0$, $S_F\left(0\right) = 12500$, 
$E_F\left(0\right) = 20$, $I_F\left(0\right) = 0$, and $M\left(0\right) = 90$. 
The objective is to illustrate the analytical findings discussed in earlier sections.


\subsection{Impact of the periodic infection rate on the occurrence of rabies outbreaks}

\begin{figure}[!ht]
\begin{minipage}[b]{0.45\textwidth}
\includegraphics[height=5.0cm, width=7.5cm]{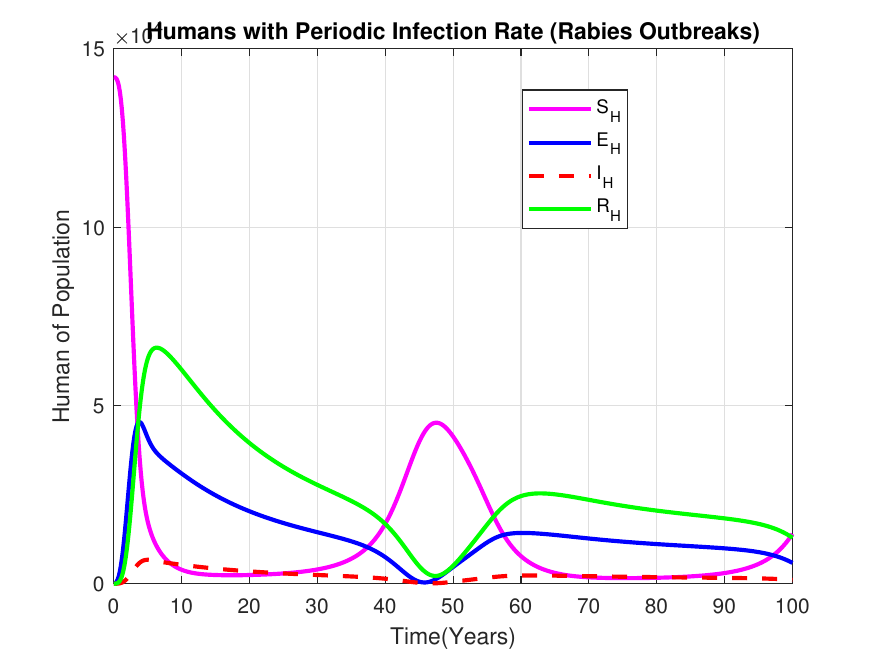}
\centering{(a)}
\end{minipage}
\begin{minipage}[b]{0.45\textwidth}
\includegraphics[height=5.0cm, width=7.5cm]{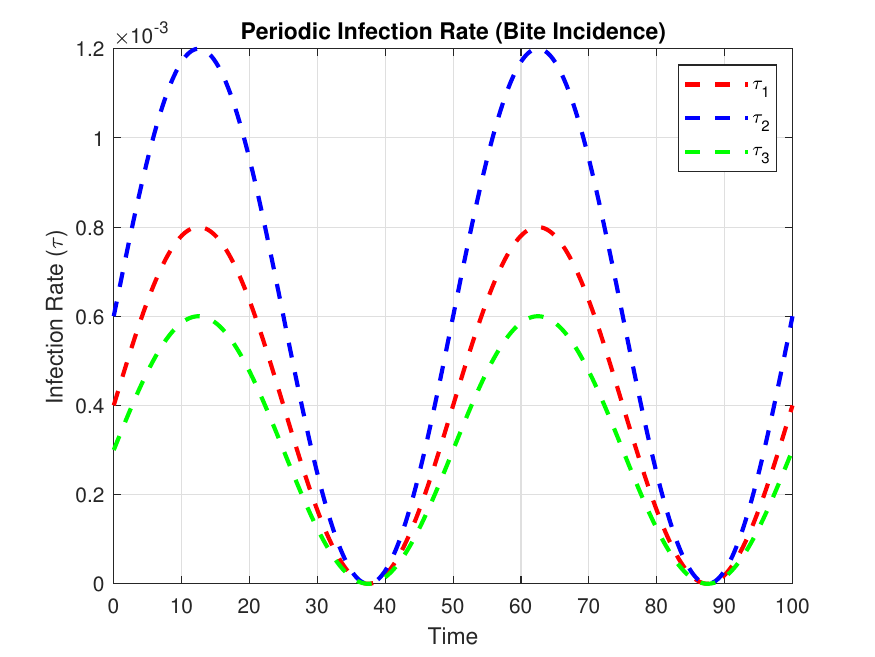}
\centering{(b)}
\end{minipage}
\centering
\caption{The impact of human population bite 
incidence on the occurrence of periodic rabies outbreaks.}
\label{Fig10}
\end{figure}

\begin{figure}[!ht]
\begin{minipage}[b]{0.45\textwidth}
\includegraphics[height=5.0cm, width=7.5cm]{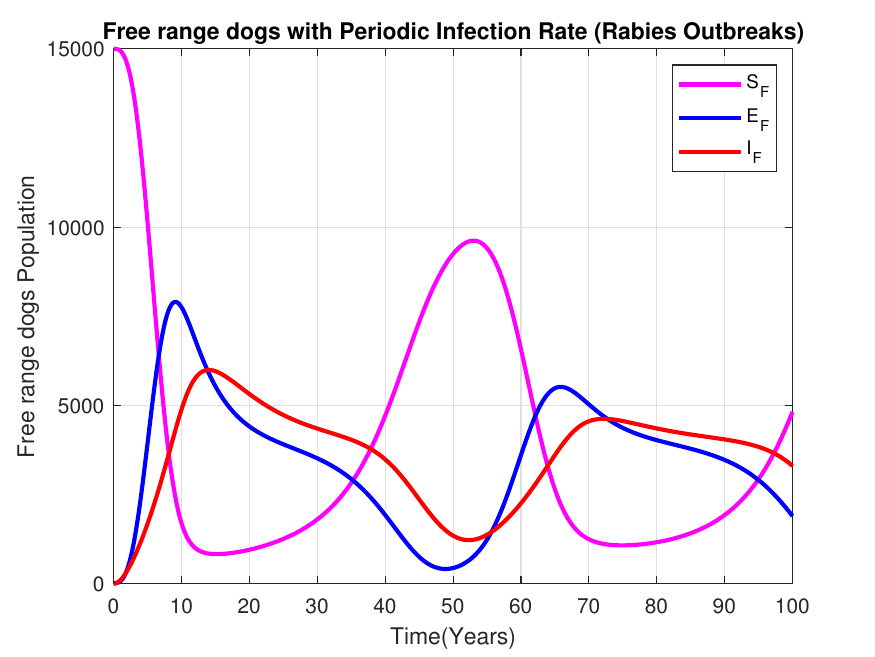}
\centering{(a)}
\end{minipage}
\begin{minipage}[b]{0.45\textwidth}
\includegraphics[height=5.0cm, width=7.5cm]{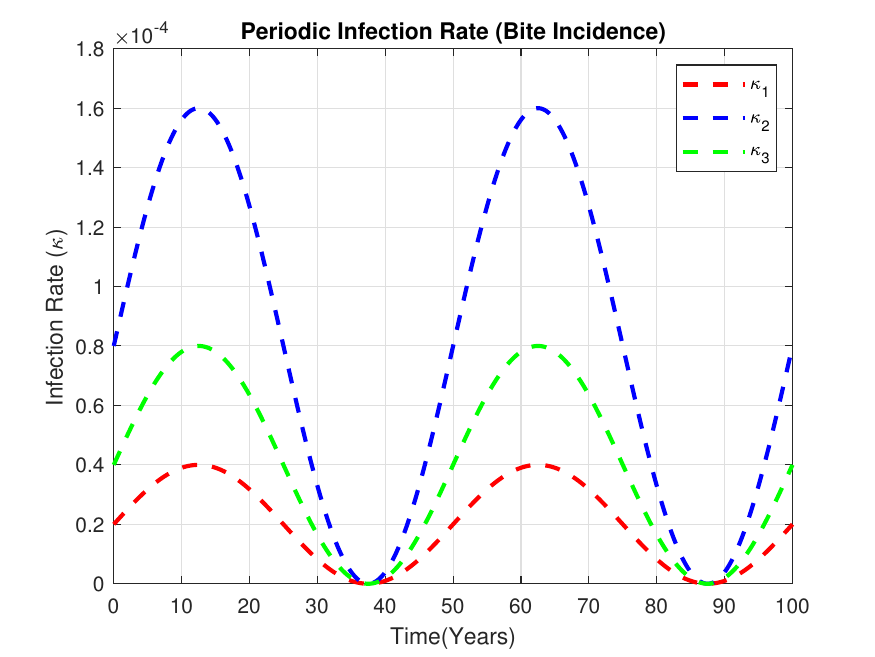}
\centering{(b)}
\end{minipage}
\centering
\caption{The impact of free range dogs bite incidence 
on the occurrence of periodic rabies outbreaks.}
\label{Fig11}
\end{figure}

\begin{figure}[!ht]
\begin{minipage}[b]{0.45\textwidth}
\includegraphics[height=5.0cm, width=7.5cm]{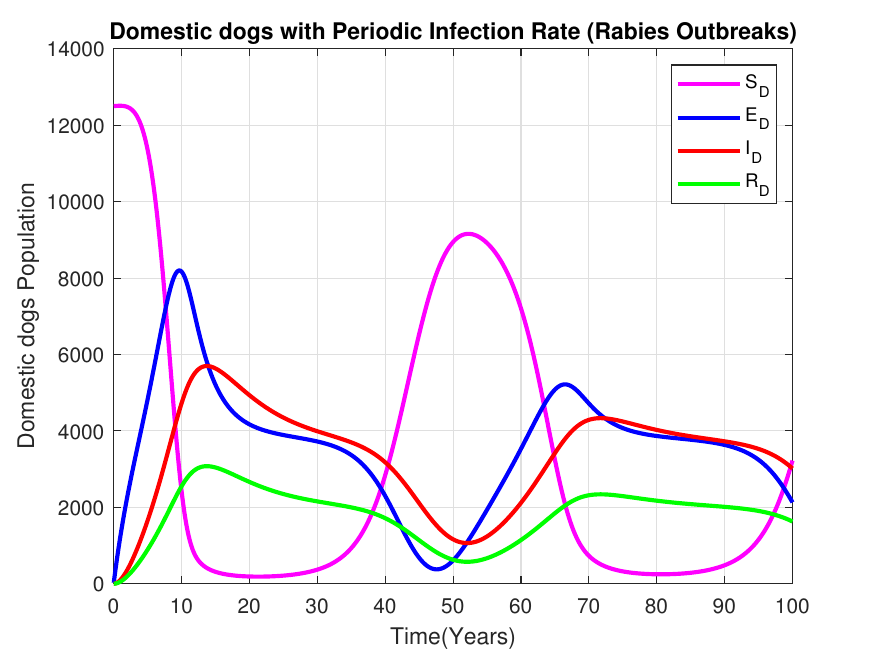}
\centering{(a)}
\end{minipage}
\begin{minipage}[b]{0.45\textwidth}
\includegraphics[height=5.0cm, width=7.5cm]{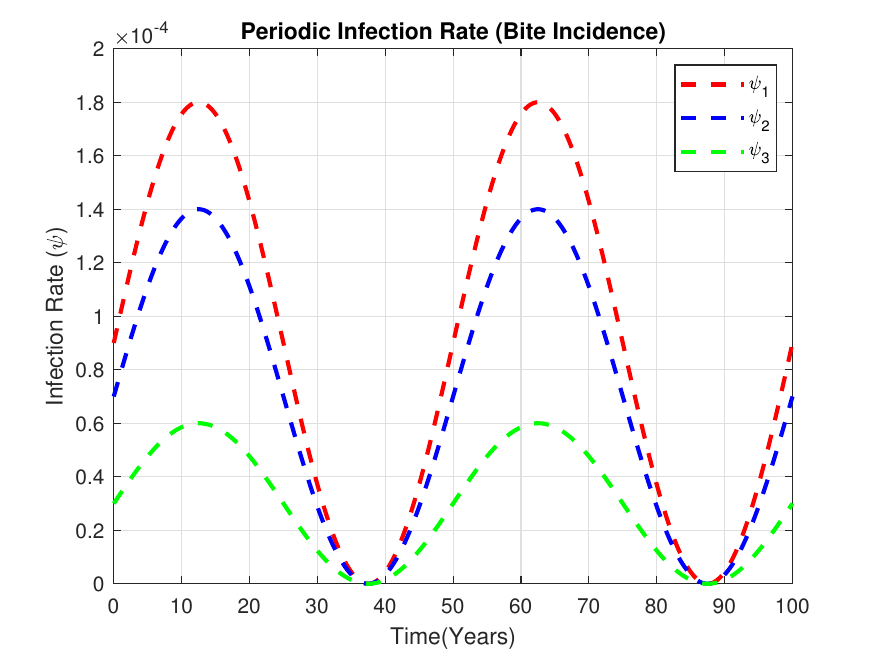}
\centering{(b)}
\end{minipage}
\centering
\caption{The impact of domestic dogs bite incidence 
on the occurrence of periodic rabies outbreaks.}
\label{Fig12}
\end{figure}

\begin{figure}[!ht]
\begin{minipage}[b]{0.45\textwidth}
\includegraphics[height=7.0cm, width=7.5cm]{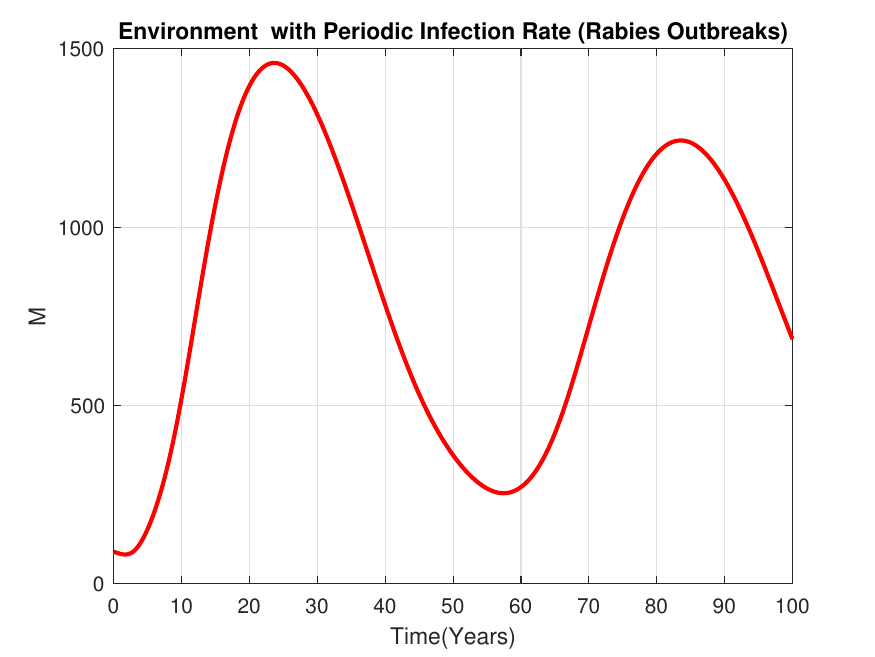}
\centering{(a)}
\end{minipage}
\begin{minipage}[b]{0.45\textwidth}
\includegraphics[height=7.0cm, width=7.5cm]{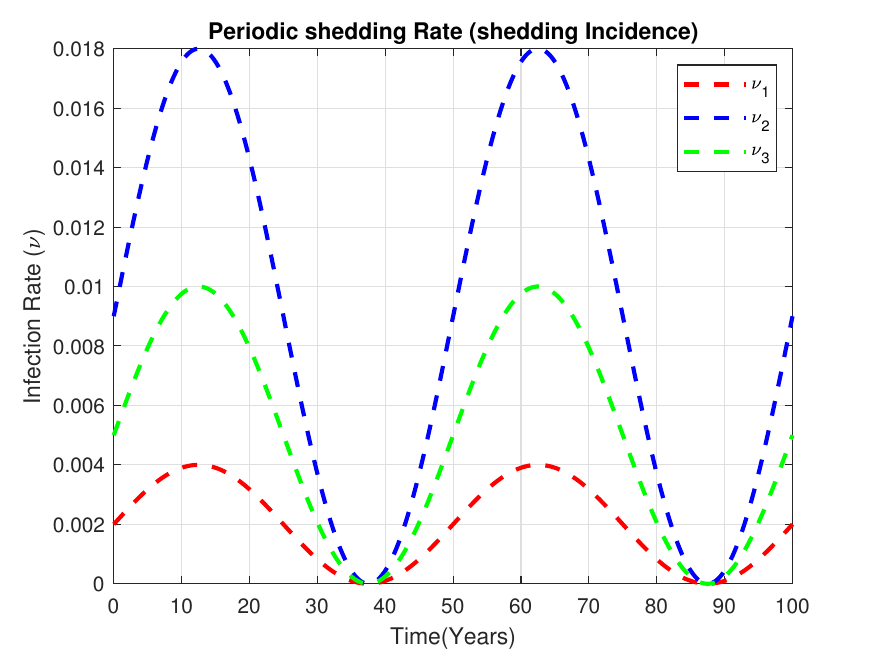}
\centering{(b)}
\end{minipage}
\centering
\caption{The impact of environment shedding incidence 
on the occurrence of periodic rabies outbreaks.}
\label{Fig13}
\end{figure}

The periodic effect of rabies on bite incidence describes the cyclic variation 
in the number of dog bites within a population due to recurrent outbreaks 
of rabies in dogs \cite{abdulmoghni2021incidence}. These outbreaks stem 
from the viral infection's influence on dog behavior, causing increased 
aggression and a propensity to bite. This cyclic pattern arises as rabies 
outbreaks occur intermittently, influenced by factors like seasonal 
fluctuations, vaccination efforts, and animal movement. To investigate 
on the effect of periodic infection for bite incidence,  
we employed the formula
\begin{equation*}
\textup{Bite Incidence}(\textup{t}) 
= \beta_{\textup{mean}}\left(1 + A\sin(2\pi\textup{f} + \phi)\right)\textup{SI},
\end{equation*}
where $\beta_{\textup{mean}}$ is the infection rate, $A$  
and $f$ are the amplitude and the frequency of the sinusoidal 
variation, respectively, $\phi$ is the phase shift, $S$ and $I$ 
are susceptible and infections, respectively, and
\begin{equation*}
\textup{f} = \frac{t}{\text{period of control of outbreak}}=\dfrac{t}{T}.
\end{equation*}
Therefore, to incorporate and unify the changing dynamics of 
rabies, transmission rates $\tau_{1}$, $\tau_{2}$, $\tau_{3}$, 
$\psi_{1}$, $\psi_{2}$, $\psi_{3}$, $\kappa_{1}$, $\kappa_{2}$,
$\kappa_{3}$, $\nu_{1}$, $\nu_{2}$ and  $\nu_{3}$, as applied  
by \citep{dutta2024periodic}, are considered as
\begin{equation*}
\left.
\begin{aligned}
\tau_{i} 
&= \tau_{(\textup{mean})}\left(1+\textup{A}_{i}\sin\left(
\frac{2\pi \textup{t}}{\textup{T}}+\phi\right)\right), 
\quad  \psi_{i} = \psi_{(\textup{mean})}\left(1+\textup{A}_{i}
\sin\left(\frac{2\pi \textup{t}}{\textup{T}}+\phi\right)\right),\\
\kappa_{i}&= \kappa_{(\textup{mean})}\left(1+\textup{A}_{i}
\sin\left(\frac{2\pi \textup{t}}{\textup{T}}+\phi\right)\right), 
\quad \nu_{i}= \nu_{(\textup{mean})}\left(1+\textup{A}_{i}\sin\left(
\frac{2\pi \textup{t}}{\textup{T}}+\phi\right)\right),\\ 
\textup{for} \;\; \textup{i}=\textup{1,  2  and  3}.
\end{aligned}
\right\}
\end{equation*}
The results of bite incidence are presented in Figure~\ref{Fig10}(a)--(b),  
Figure~\ref{Fig11}(a)--(b), Figure~\ref{Fig12}(a)--(b),  
and Figure~\ref{Fig13}(a)--(b).

Figure~\ref{Fig10}(a)--(b), Figure~\ref{Fig11}(a)--(b), 
and Figure~\ref{Fig12}(a)--(b) illustrate that an increase in bite 
incidents leads to a rise in the number of exposed and infected individuals 
while reducing the number of susceptible individuals in both human and dog populations. 
Furthermore, all these figures demonstrate that the rabies outbreak, 
driven by a higher infection rate, remains active within the first 20 years and, 
subsequently, exhibits periodic declines. This decline in the number of infections 
is attributed to the decrease in the infection rate in both populations. 
On the other hand, Figure~\ref{Fig13}(a)--(b) shows that an increase in the rate of 
shedding into the environment results in periodic rises in rabies contamination 
within the environment. These scenarios highlight the significance of effective 
vaccination campaigns, responsible pet ownership, and timely post-exposure 
prophylaxis for individuals who have been bitten. These measures are essential 
for managing the public health impact of this periodic phenomenon, 
underscoring the importance of rabies control strategies.
   	
\begin{figure}[!ht]
\subfigure[Effect of ${\cal R}_0$ with respect to $\kappa_{1}$ and $\kappa_{2}$.]{
\includegraphics[width=0.46\textwidth]{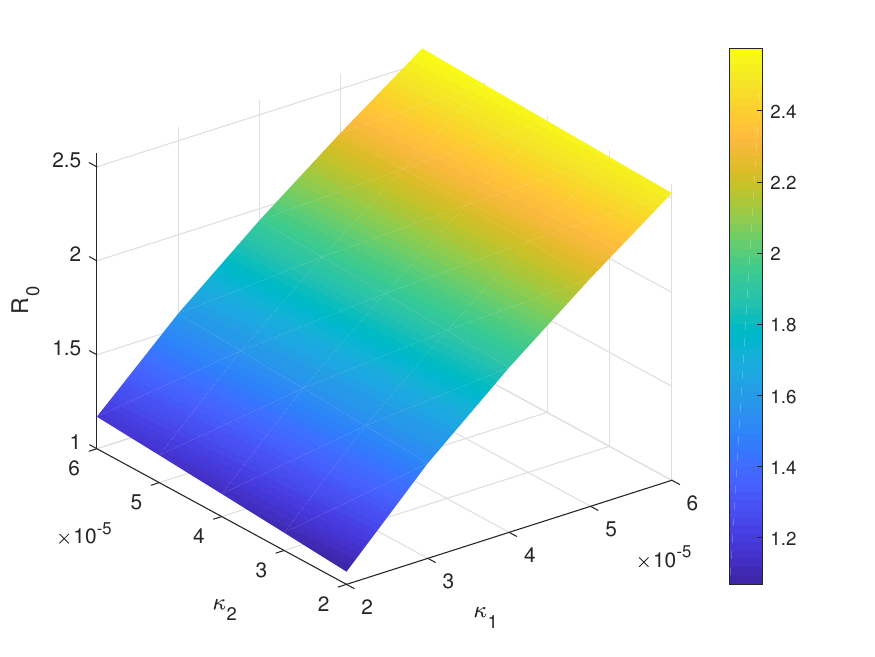}}\hfill
\subfigure[Effect of ${\cal R}_0$ with respect to $\psi_{1}$ and $\psi_{2}$.]{
\includegraphics[width=0.46\textwidth]{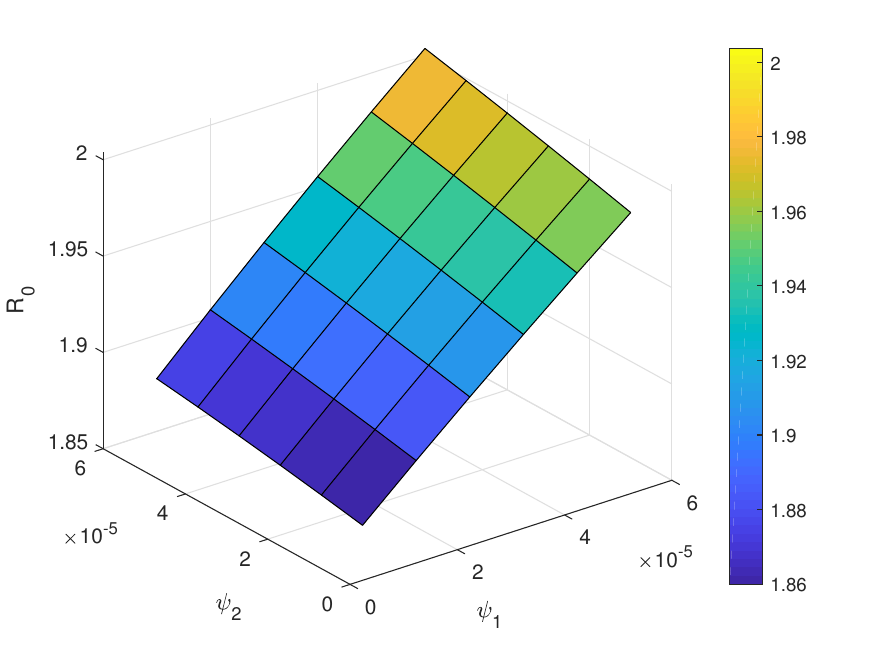}}
\caption{Effect of ${\cal R}_0$ with respect to $\kappa_{1}$, 
$\kappa_{2}$, $\psi_{1}$, and $\psi_{2}$.}
\label{fig:6}
\end{figure}

The results presented in Figure~\ref{fig:6}(a)--(b) demonstrate that the parameters 
$\psi_{1}$, $\psi_{2}$, $\kappa_{1}$, and $\kappa_{2}$ have a positive impact 
on the basic reproductive number, $\mathcal{R}_{0}$. 
The study reveals that changes in these parameters generate varying effects on 
$\mathcal{R}_{0}$, ranging from 1.8 to 2.0 and 1.2 to 2.4, respectively. These 
findings support the estimates provided by \cite{kadowaki2018risk} and suggest 
that intervention strategies can have a significant impact on the incidence of 
rabies in a given population. Furthermore, the parameter values outlined 
in Table~\ref{table3} indicate that increasing $\kappa_{1}$, $\kappa_{2}$, 
$\psi_{1}$, and $\psi_{2}$ corresponds to an increase in $\mathcal{R}_{0}$.


\subsection{Effect of varying the most sensitive parameters}

Now we investigate the impact of the contact rate between 
infectious agent and: (i)~susceptible human; (ii)~susceptible domestic dogs;
and (iii)~free range dogs.

\subsubsection{Impact of contact rate between  infectious agent and  susceptible human}

\begin{figure}[!ht]
\begin{minipage}[b]{0.45\textwidth}
\includegraphics[height=5.0cm, width=7.5cm]{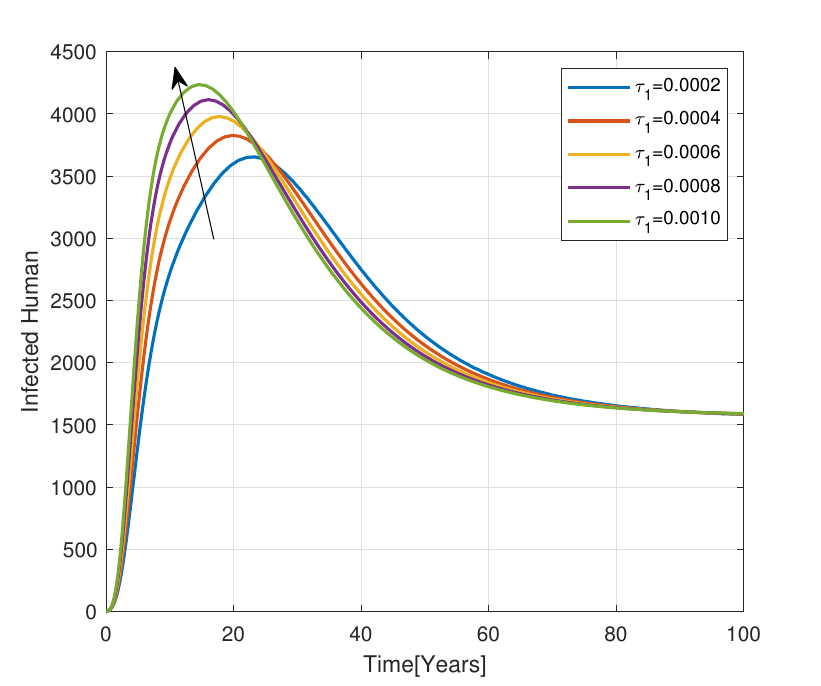}
\centering{(a)}
\end{minipage}
\begin{minipage}[b]{0.45\textwidth}
\includegraphics[height=5.0cm, width=7.5cm]{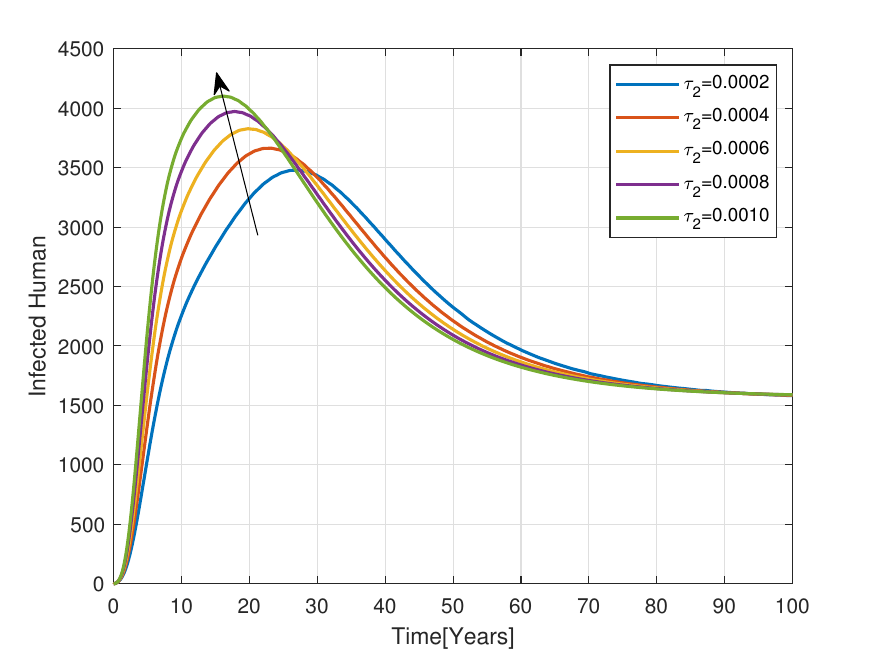}
\centering{(b)}
\end{minipage}
\begin{minipage}[b]{0.45\textwidth}
\includegraphics[height=5.0cm, width=7.5cm]{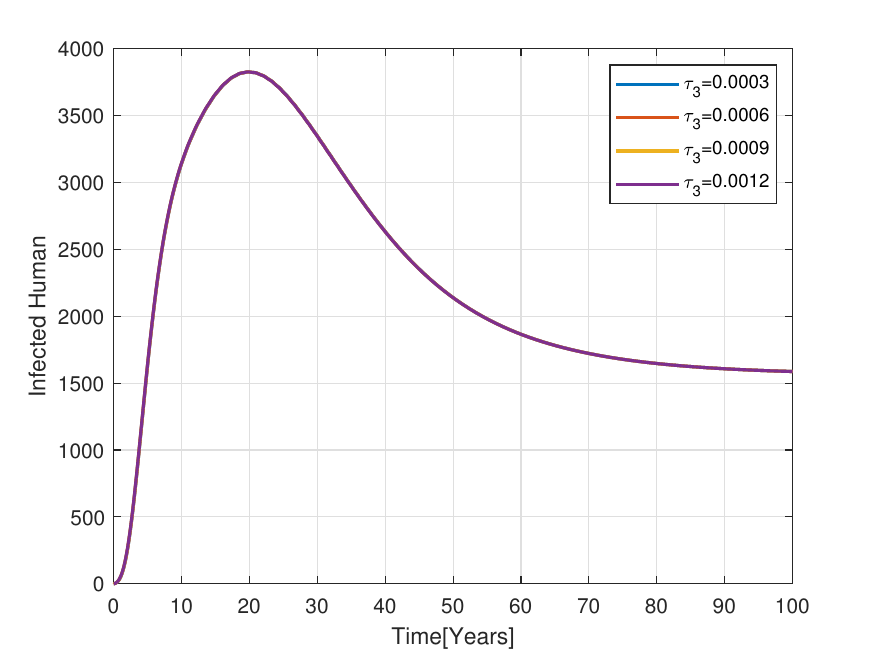}
\centering{(c)}
\end{minipage}
\centering
\caption{Simulation results of model (\ref{eqn1}) 
for $I_H$ with respect to $\tau_{1}$, $\tau_{2}$, and $\tau_{3}$.}
\label{Fig6}
\end{figure}

The findings presented in Figure~\ref{Fig6}(a)--(c) demonstrate 
that the contact rates $\tau_1$, $\tau_2$, and $\tau_3$ exert 
a significant influence on the transmission dynamics among susceptible humans, 
infected free-range and domestic dogs, and the environment. Nevertheless, 
it is noted that these dynamics ultimately reach a stable state after 80 years. 
This underscores the crucial role of education and awareness in mitigating 
the transmission of rabies among the human population by reducing contact 
between susceptible humans and sources carrying the rabies virus.


\subsubsection{Impact of contact rate between  infectious agent and susceptible domestic dogs}

\begin{figure} [!ht]
\begin{minipage}[b]{0.45\textwidth}
\includegraphics[height=4.0cm, width=7.5cm]{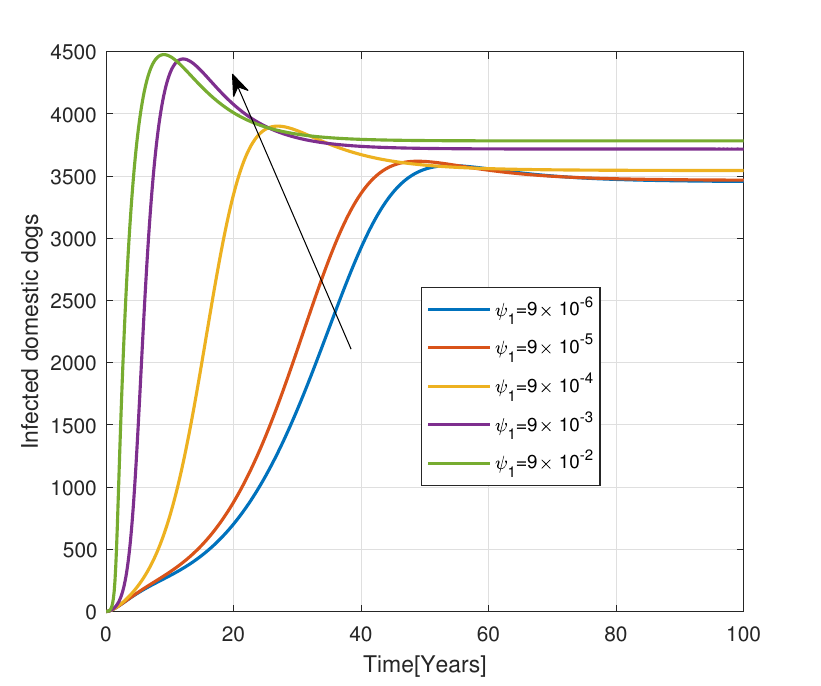}
\centering{(a)}
\end{minipage}
\begin{minipage}[b]{0.45\textwidth}
\includegraphics[height=4.0cm, width=7.5cm]{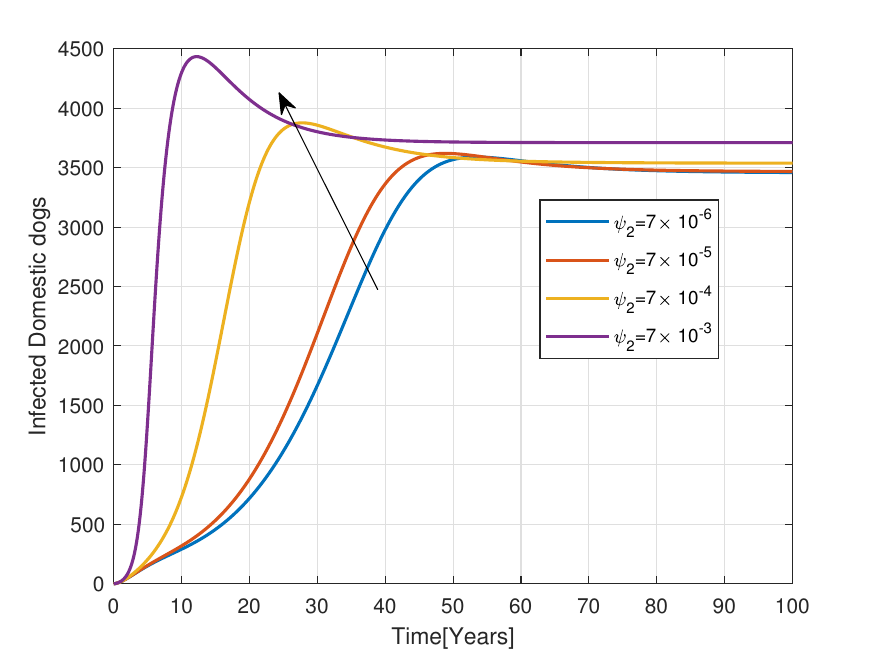}
\centering{(b)}
\end{minipage}
\begin{minipage}[b]{0.45\textwidth}
\includegraphics[height=5.0cm, width=7.5cm]{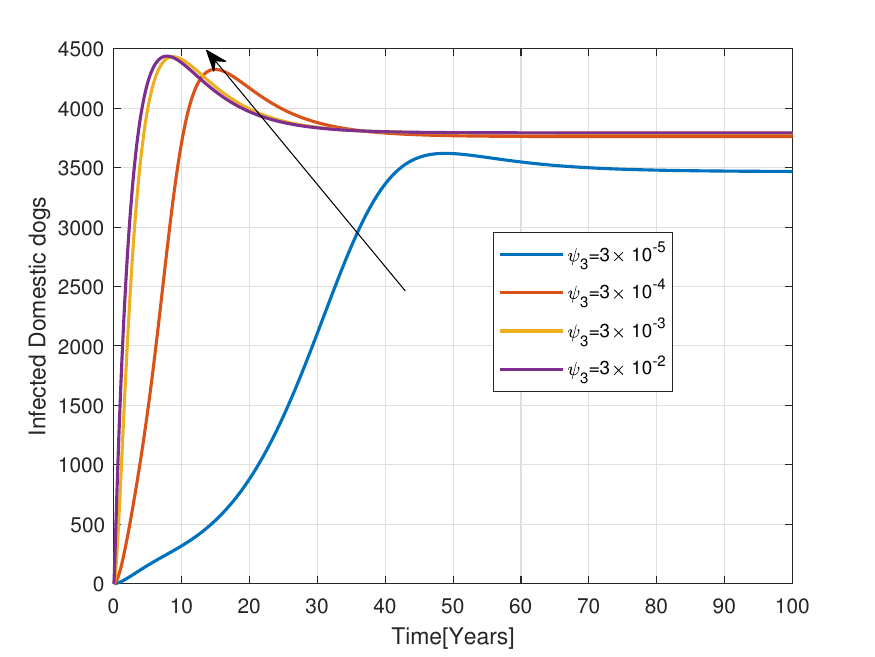}
\centering{(c)}
\end{minipage}
\centering
\caption{Simulation results of model (\ref{eqn1}) for $I_D$ 
with respect to $\psi_{1}$, $\psi_{2}$, and $\psi_{3}$.}
\label{Fig7}
\end{figure}

The findings presented in Figure~\ref{Fig7}(a)--(c) reveal that an increase 
in contact rates, denoted as $\psi_{1}$, $\psi_{2}$, and $\psi_{3}$, results 
in a higher prevalence of rabies in domestic dogs. After approximately 50 years, 
the number of infected dogs reaches a steady state, implying that mitigating 
the contact between susceptible, infected, and free-range dogs and the environment 
carrying the rabies virus is critical to reduce the transmission of the disease.


\subsubsection{Impact of contact rate between infectious agent and free range dogs}

\begin{figure}[!ht]
\begin{minipage}[b]{0.45\textwidth}
\includegraphics[height=4.0cm, width=7.5cm]{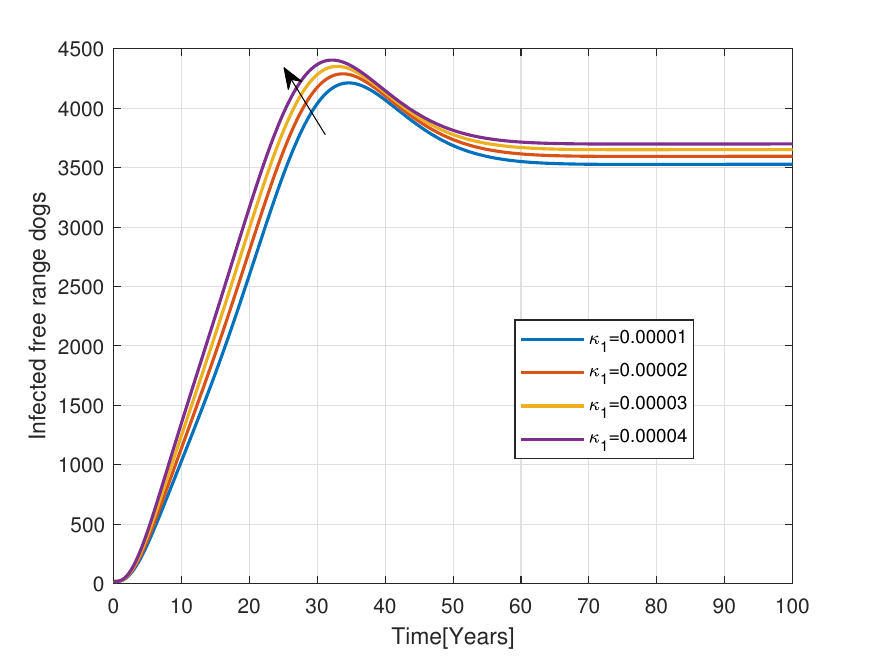}
\centering{(a)}
\end{minipage}
\begin{minipage}[b]{0.45\textwidth}
\includegraphics[height=4.0cm, width=7.5cm]{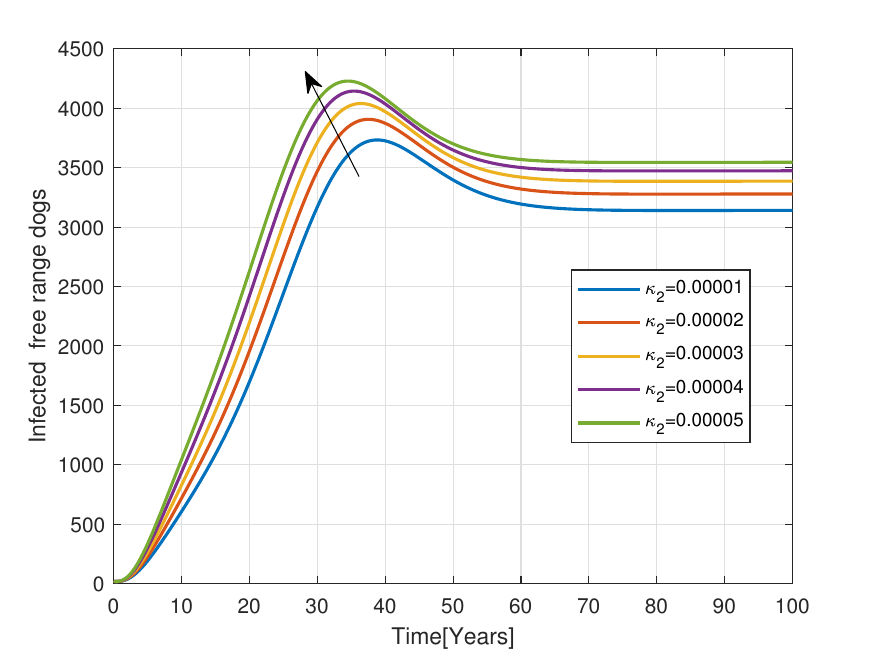}
\centering{(b)}
\end{minipage}
\begin{minipage}[b]{0.45\textwidth}
\includegraphics[height=4.0cm, width=7.5cm]{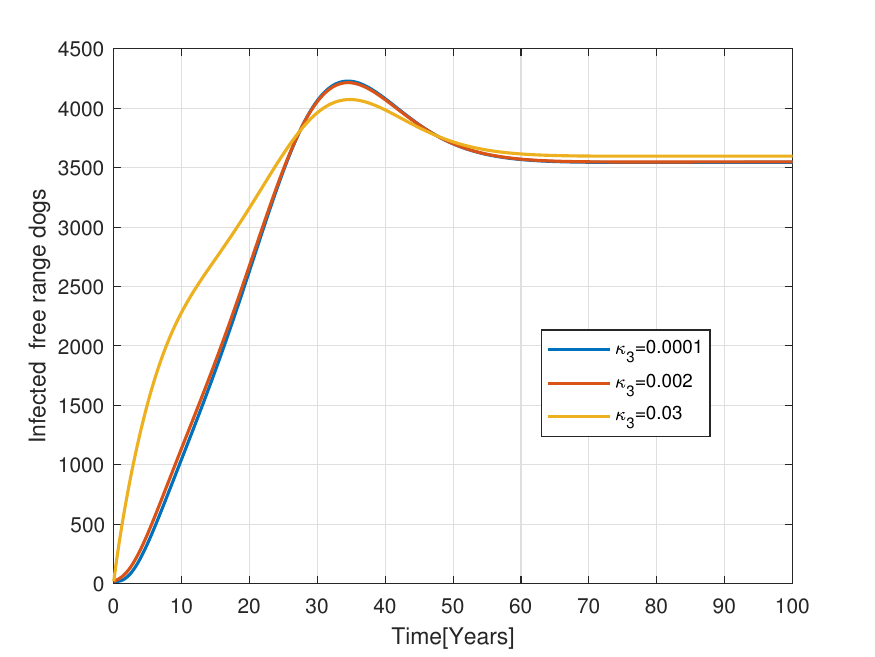}
\centering{(c)}
\end{minipage}
\centering
\caption{Simulation results of model (\ref{eqn1}) for 
$I_F$ with respect to $\tau_{1}$, $\tau_{2}$, and $\tau_{3}$.}
\label{Fig8}
\end{figure}

Figure~\ref{Fig8}(a) presents the finding that an increase in the contact rate $\kappa_{1}$ 
results in a higher number of susceptible free-range dogs becoming infected, which suggests 
an inadequacy of control measures. Conversely, Figure~\ref{Fig8}(b) portrays that a rise 
in the contact rate $\kappa_{2}$ with free-range dogs leads to an increase in carriers 
and symptomatic infections. In addition, Figure~\ref{Fig8}(c) indicates that an increase 
in the contact rate $\kappa_{3}$ between free-range dogs and the environment yields a 
slight upsurge in the number of infectious individuals or carriers.


\section{Discussion}
\label{sec:5}

Rabies is a fatal viral disease that can easily spread from an infected animal to a human, 
making it a significant public health concern worldwide. Dogs are the primary reservoir 
and transmitter of rabies to humans, causing most human cases.
To understand the transmission dynamics of rabies and develop effective prevention 
and control strategies, a study was conducted. The study aimed to create a deterministic 
model to investigate how changes in contact rates and environmental conditions impact 
the spread of rabies. Mathematical tools such as Jacobian and Metzeler matrices were 
used to conduct stability analyses and uncover the underlying dynamics of rabies transmission.
The study also aimed to determine the relationship between contact rates, environmental factors, 
and the basic reproduction number ${\mathcal R}_{0}$, which is a crucial indicator 
of disease spread. By gaining insights into the complex dynamics of rabies transmission, 
the study hopes to contribute to the development of targeted and sustainable 
strategies for its prevention and control.
 

\section{Conclusion} 
\label{sec:6}

The transmission of rabies among humans and dogs is influenced by their contact rate 
and environmental factors. A deterministic model was created to investigate how these 
factors affect the spread of rabies that results from dog bites. Stability analysis was 
conducted using Jacobian and Metzelar matrices. The study's numerical simulations showed 
that the transmission of rabies from dog bites has serious consequences for both human 
and dog populations. Furthermore, the study found that an increase in the contact rate 
(such as $\psi_{{1}}$, $\psi_{{2}}$, $\psi_{3}$, $\kappa_{{1}}$, $\kappa_{{1}}$, 
and $\kappa_{{3}}$) leads to a rise in basic reproduction ${\mathcal R}_{0}$. 
By examining the relationship between contact rate, environmental impact, 
and the transmission dynamics of rabies, this study provides insights into 
the complexity of rabies transmission. Ultimately, it contributes to the 
development of targeted, sustainable strategies for preventing 
and controlling the spread of rabies.


\section*{Data Availability}

The data used in this study is available from 
the corresponding author upon request.


\section*{Funding Statement}

Torres is supported by the Portuguese Foundation for Science and Technology (FCT) 
and CIDMA, project UIDB/04106/2020.


\section*{Conflicts of Interest}

The authors declare that they have no known competing 
financial interests or personal relationships that could 
have appeared to influence the work reported in this paper.


\section*{Acknowledgments}

The authors thank the handling editor and reviewer for helpful suggestions that improved 
this paper's quality. We acknowledge The Nelson Mandela African Institution of Science 
and Technology (NM-AIST) and the College of Business Education (CBE) for providing 
a conducive environment. 


\appendix

\section*{Appendix A. Global Stability of the Endemic Equilibrium ${\mathbb E}^{*}$}

Here we prove the global stability of the endemic equilibrium characterized in Section~\ref{sec:EEp}.

\begin{theorem}
The endemic equilibrium point ${\mathbb E}^{*}$ of the rabies model 
$\left(\ref{eqn1}\right)$ is globally asymptotically stable when 
the disease persists in the short and long term. Any invariant set close to 
${\mathbb E}^{*}$ in $\Omega$ converges to the equilibrium point over an extended period.
\label{Th3}
\end{theorem}

\begin{proof}
To prove Theorem~\ref{Th3}, we adopt the approach of 
\cite{ayoade2023modeling,osman2022analysis,musaili2020mathematical} 
by constructing a Lyapunov function of the form
\begin{equation*}
{\cal L}=\sum^{n}_{i=1}G_{i}\left(y_{i}-y_{i}^{*}+y_{i}^{*}
\ln \left(\frac{y_{i}^{*}}{y_{i}}\right)\right), G_{i}>0 
\;\;\text{for} \;\; i=1,2,3,\ldots,n,
\end{equation*}
where $G_{i}$ represents a positive constant that needs to be determined, 
$y_{i}$ stands for the population variable at compartment $i$, and 
$y^{*}_{i}$ denotes the equilibrium point of the rabies model at compartment $i$ 
for $i \in \{1,2,3,\ldots,12\}$. Therefore, we define the Lyapunov $\mathcal{L}$ 
for model system $\left(\ref{eqn1}\right)$ as follows:
\begin{equation}
\begin{array}{llll}
\mathcal{L}= G_{1}\left(S_{H}-S_{H}^{*}+S_{H}\ln\left(\frac{S_{H}^{*}}{S_{H}}\right)\right)
+G_{2}\left(E_{H}-E_{H}^{*}+E_{H}\ln\left(\frac{E_{H}^{*}}{E_{H}}\right)\right)
+G_{3}\left(I_{H}-I_{H}^{*}+I_{H}\ln\left(\frac{I_{H}^{*}}{I_{H}}\right)\right)\\
+G_{4}\left(R_{H}-R_{H}^{*}+R_{H}\ln\left(\frac{R_{H}^{*}}{R_{H}}\right)\right)
+G_{5}\left(S_{F}-S_{F}^{*}+S_{F}\ln\left(\frac{S_{F}^{*}}{S_{F}}\right)\right)
+G_{6}\left(E_{F}-E_{F}^{*}+E_{F}\ln\left(\frac{E_{F}^{*}}{E_{F}}\right)\right)\\
+G_{7}\left(I_{F}-I_{F}^{*}+I_{F}\ln\left(\frac{I_{F}^{*}}{I_{F}}\right)\right)
+G_{8}\left(S_{D}-S_{D}^{*}+S_{D}\ln\left(\frac{S_{D}^{*}}{S_{D}}\right)\right)
+G_{9}\left(E_{D}-E_{D}^{*}+E_{D}\ln\left(\frac{E_{D}^{*}}{E_{D}}\right)\right)\\
+G_{10}\left(I_{D}-I_{D}^{*}+I_{D}\ln\left(\frac{I_{D}^{*}}{I_{D}}\right)\right)
+G_{11}\left(R_{D}-R_{D}^{*}+R_{D}\ln\left(\frac{R_{D}^{*}}{R_{D}}\right)\right)
+G_{12}\left(M-M^{*}+M\ln\left(\frac{M^{*}}{M}\right)\right).
\end{array}
\label{eqn14}
\end{equation}
Evaluating  equation $\left(\ref{eqn14}\right)$ at the
endemic equilibrium point $\mathbb{E}^{*}$ gives 
\begin{equation*}
\mathcal{L} =\mathbb{E} ^{*}\left(S_{H}^{*}\;, E_{H}^{*}\;, I_{H}^{*}\;, 
R_{H}^{*} \;, S_{F}^{*} \;,E_{F}^{*}\;, I_{F}^{*}\;, S_{D}^{*}\;,
E_{D}^{*} \;,I_{D}^{*} \;, R_{D}^{*},M^{*}\right)=0.
\end{equation*}
Then, using the time derivative of the Lyapunov function $\mathcal{L}$ 
in equation $\left(\ref{eqn14}\right)$ gives
\begin{equation}
\left.
\begin{array}{llll}
\dfrac{d\mathcal{L}}{dt} 
&= G_{1}\left(1-\dfrac{S^{*}_H}{S_{H}} \right)\dfrac{dS_H}{dt}
+G_{2}\left(1-\dfrac{E^{*}_H}{E_H}\right)\dfrac{dE_H}{dt}
+G_{3}\left(1-\dfrac{I^{*}_H}{I_H}\right)\dfrac{dI_H}{dt}
+G_{4}\left(1-\dfrac{R^{*}_H}{R_H}\right)\dfrac{dR_H}{dt}\\ 
&+G_{5}\left(1-\dfrac{S^{*}_F}{S_{F}} \right)\dfrac{dS_F}{dt}
+G_{6}\left(1-\dfrac{E^{*}_F}{E_{F}} \right)\dfrac{dE_F}{dt}
+G_{7}\left(1-\dfrac{I^{*}_F}{I_{F}} \right)\dfrac{dI_F}{dt}
+G_{8}\left(1-\dfrac{S^{*}_D}{S_{D}} \right)\dfrac{dS_D}{dt}\\
&+G_{9}\left(1-\dfrac{E^{*}_D}{E_{D}} \right)\dfrac{dE_D}{dt}
+G_{10}\left(1-\dfrac{I^{*}_D}{I_{D}} \right)\dfrac{dI_D}{dt}
+G_{11}\left(1-\dfrac{R^{*}_D}{R_{D}} \right)\dfrac{dR_D}{dt}
+G_{12}\left(1-\dfrac{M^{*}}{M} \right)\dfrac{dM}{dt}.
\end{array}
\right\}
\label{eqn15}
\end{equation}
Consider the endemic equilibrium point $\left(EEP\right), \mathbb{E}^{*}$ 
of equation $\left(\ref{eqn1}\right)$ such that 
\begin{equation}
\left.
\begin{array}{llll}
\theta_{1}&=\left(\tau_{1}I^{*}_{F}+\tau_{2}I^{*}_{D}
+\tau_{3}\lambda\left(M^{*}\right)\right)S^{*}_{H}
+\mu_{1}S^{*}_{H}-\beta_{3}R^{*}_{H},\;\; \mu_{1}+\beta_{1}+\beta_{2}
=\dfrac{\left(\tau_{1}I^{*}_{F}+\tau_{2}I^{*}_{D}+\tau_{3}
\lambda\left(M^{*}\right)\right)S^{*}_{H}}{E^{*}_{H}},\\ 
& \sigma_{1}+\mu_{1}=\dfrac{\beta_{1}E^{*}_{H}}{I^{*}_{H}},\;\;
\beta_3+\mu_{1} =\dfrac{\beta_{2}E^{*}_{H}}{R^{*}_{H}},\;\;
\theta_{2}=\left(\kappa_{1}I^{*}_{F}+\kappa_{2}I^{*}_{D}
+\kappa_{3}\lambda\left(M^{*}\right)\right)S^{*}_{F}+\mu_{2}S_{F},\\
&\mu_{2}+\gamma=\dfrac{\left(\kappa_{1}I^{*}_{F}+\kappa_{2}I^{*}_{D}
+\tau_{3}\lambda\left(M^{*}\right)\right)S^{*}_{F}}{E^{*}_{F}},\;\;
\sigma_{2}+\mu_{2}=\dfrac{\gamma E^{*}_{F}}{I^{*}_{F}},\\ 
&\theta_{3}=\left(\dfrac{\psi_1I^{*}_{F}}{1+\rho_{1}}
+\dfrac{\psi_2I^{*}_{D}}{1+\rho_{2}}\dfrac{\psi_3\lambda
\left(M^{*}\right)}{1+\rho_{3}}\right)S^{*}_{D}
+\mu_{3}S^{*}_{D}-\gamma_{3}R^{*}_{D},\;\;
\mu_{3}+\gamma_{1}+\gamma_{2}=\dfrac{\left(\dfrac{\psi_1I^{*}_{F}}{1
+\rho_{1}}+\dfrac{\psi_2I^{*}_{D}}{1+\rho_{2}}
\dfrac{\psi_3\lambda\left(M^{*}\right)}{1+\rho_{3}}\right)S^{*}_{D}}{E^{*}_{D}},\\
&\sigma_{3}+\mu_{3}=\dfrac{\gamma_{1} E^{*}_{D}}{I^{*}_{D}},\;\;
\gamma_3+\mu_{3}=\dfrac{\gamma_{2}E^{*}_{D}}{R^{*}_{D}},\;\;
\mu_{4}=\dfrac{\left(\nu_{1}I^{*}_{H}+\nu_{2}I^{*}_{F}
+\nu_{3}I^{*}_{D}\right)}{M^{*}}.
\end{array}
\right\}
\label{eqn16}
\end{equation}
Then, by substituting $\eqref{eqn1}$ into $\eqref{eqn15}$, we have
\begin{equation}
\left.
\begin{array}{llll}
\frac{d\mathcal{L}}{dt} 
&= G_{1}\left(1-\frac{S^{*}_H}{S_{H}} \right)\left(\theta_{1}
+\beta_{3}R_{H}-\mu_{1} S_{H}-\chi_{1}\right)
+G_{2}\left(1-\frac{E^{*}_H}{E_H} \right)\left(\chi_{1}
-\left(\mu_{1}+\beta_{1}+\beta_{2}\right)E_{H}\right)\\
&+G_{3}\left(1-\frac{I^{*}_H}{I_H}\right)\left(\beta_{1}
E_{H}-\left(\sigma_{1}+\mu_{1}\right) I_{H}\right)
+G_{4}\left(1-\frac{R^{*}_H}{R_H}\right)\left(\beta_{2} E_{H}
-\left(\beta_{3}+\mu_{1} \right) R_{H}\right)\\ 
&+G_{5}\left(1-\frac{S^{*}_F}{S_{F}} \right)\left(\theta_{2}
-\chi_{2}-\mu_{2}S_{F}\right)+G_{6}\left(1-\frac{E^{*}_F}{E_{F}} 
\right)\left(\chi_{2}-\left(\mu_{2}+\gamma\right)E_{F}\right)\\
&+G_{7}\left(1-\frac{I^{*}_F}{I_{F}} \right)\left(\gamma E_{F}
-\left(\mu_{2}+\sigma_{2}\right)I_{F}\right)
+G_{8}\left(1-\frac{S^{*}_D}{S_{D}} \right)\left(\theta_{3}
-\mu_{3}S_{D}-\chi_{3}+\gamma_{3}R_{D}\right)\\
&+G_{9}\left(1-\frac{E^{*}_D}{E_{D}} \right)\left(\chi_{3}
-\left(\mu_{3}+\gamma_{1}+\gamma_{2}\right) E_{D}\right)
+G_{10}\left(1-\frac{I^{*}_D}{I_{D}} \right)\left(\gamma_{1}E_{D}
-\left(\mu_{3}+\delta_{3}\right) I_{D}\right)\\
&+G_{11}\left(1-\frac{R^{*}_D}{R_{D}} \right)\left(\gamma_{2}
E_{D}-\left(\mu_{3}+\gamma_{3}\right)R_{D}\right)
+G_{12}\left(1-\frac{M^{*}}{M} \right)\left(\left(\nu_1I_H
+\nu_2I_F+\nu_3I_D\right)-\mu_4M\right).
\end{array}
\right\}
\label{eqn27}
\end{equation}
Using the endemic equilibrium point $\left(EEP\right)$ 
described in equation $\left(\ref{eqn16}\right)$,  
we simplify the equation $\left(\ref{eqn27}\right)$ as 
\begin{equation}
\left.
\begin{array}{llll}
\frac{d\mathcal{L}}{dt} 
&= G_{1}\left(1-\frac{S^{*}_H}{S_{H}}\right)\Bigg(\left(\tau_{1}I^{*}_{F}
+\tau_{2}I^{*}_{D}+\tau_{3}\lambda\left(M^{*}\right)\right)S^{*}_{H}
+\mu_{1}S^{*}_{H}-\beta_{3}R^{*}_{H}-\mu_{1}S_{H}\\
&\quad -\left(\tau_{1}I_{F}+\tau_{2}I_{D}+\tau_{3}
\lambda\left(M\right)\right)S_{H}+\beta_{3}R_{H}\Bigg)\\
&+G_{2}\left(1-\frac{E^{*}_H}{E_{H}}\right)\Bigg(
\left(\tau_{1}I_{F}+\tau_{2}I_{D}+\tau_{3}\lambda\left(M\right)\right)S_{H}
-\dfrac{\left(\tau_{1}I^{*}_{F}+\tau_{2}I^{*}_{D}+\tau_{3}
\lambda\left(M^{*}\right)\right)S^{*}_{H}E_{H}}{E^{*}_{H}}\Bigg)\\
&+G_{3}\left(1-\frac{I^{*}_H}{I_{H}}\right)\left(\beta_{1}E_{H}
-\dfrac{\beta_{1}E^{*}_{H}I_{H}}{I^{*}_{H}}\right)
+G_{4}\left(1-\frac{R^{*}_H}{R_{H}}\right)\left(\beta_{2}E_{H}
-\dfrac{\beta_{2}E^{*}_{H}R_{H}}{R^{*}_{H}}\right)\\
&+G_{5}\left(1-\frac{S^{*}_F}{S_{F}}\right)\Bigg(\left(
\kappa_{1}I^{*}_{F}+\kappa_{2}I^{*}_{D}+\kappa_{3}
\lambda\left(M^{*}\right)\right)S^{*}_{F}+\mu_{2}S^{*}_{F}-\mu_{2}S_{F}\\
&\quad -\left(\kappa_{1}I_{F}+\kappa_{2}I_{D}
+\kappa_{3}\lambda\left(M\right)\right)S_{F}\Bigg)\\
&+G_{6}\left(1-\frac{E^{*}_F}{E_{F}}\right)\Bigg(
\left(\kappa_{1}I_{F}+\kappa_{2}I_{D}+\kappa_{3}
\lambda\left(M\right)\right)S_{F}-\dfrac{\left(\tau_{1}I^{*}_{F}
+\kappa_{2}I^{*}_{D}+\kappa_{3}\lambda\left(M^{*}\right)\right)
S^{*}_{F}E_{F}}{E^{*}_{F}}\Bigg)\\
&+G_{7}\left(1-\frac{I^{*}_F}{I_{F}} \right)\left(\gamma E_{F}
-\dfrac{\gamma E^{*}_{F}I_{F}}{I^{*}_{F}}\right)\\
&+G_{8}\left(1-\frac{S^{*}_D}{S_{D}}\right)\Bigg(
\left(\dfrac{\psi_1I^{*}_{F}}{1+\rho_{1}}+\dfrac{\psi_2I^{*}_{D}}{1
+\rho_{2}}\dfrac{\psi_3\lambda\left(M^{*}\right)}{1+\rho_{3}}\right)
S^{*}_{D}+\mu_{3}S^{*}_{D}-\gamma_{3}R^{*}_{D}\\
&\quad -\left(\dfrac{\psi_1I_{F}}{1+\rho_{1}}+\dfrac{\psi_2I_{D}}{1+\rho_{2}}
\dfrac{\psi_3\lambda\left(M\right)}{1+\rho_{3}}\right)S_{D}
-\mu_{3}S_{D}+\gamma_{3}R_{D}\Bigg)\\
&+G_{9}\left(1-\frac{E^{*}_D}{E_{D}}\right)\Bigg(\left(\dfrac{\psi_1I_{F}}{1
+\rho_{1}}+\dfrac{\psi_2I_{D}}{1+\rho_{2}}\dfrac{\psi_3\lambda\left(M\right)}{1
+\rho_{3}}\right)S_{D}-\dfrac{\left(\dfrac{\psi_1I^{*}_{F}}{1+\rho_{1}}
+\dfrac{\psi_2I^{*}_{D}}{1+\rho_{2}}\dfrac{\psi_3\lambda\left(M^{*}\right)}{1
+\rho_{3}}\right)S^{*}_{D}E_{D}}{E^{*}_{D}}\Bigg)\\
&+G_{10}\left(1-\frac{I^{*}_D}{I_{D}}\right)\left(\gamma_{1}E_{D}
-\dfrac{\gamma_{1}E^{*}_{D}I_{D}}{I^{*}_{D}}\right)+G_{11}\left(1
-\frac{R^{*}_D}{R_{D}} \right)\left(\gamma_{2}E_{D}
-\dfrac{\gamma_{2}E^{*}_{D}R_{D}}{R^{*}_{D}}\right)\\
&+G_{12}\left(1-\frac{M^{*}}{M}\right)\left(\nu_{1}I_{H}+\nu_{2}I_{F}
+\nu_{3}I_{D}-\dfrac{\left(\nu_{1}I^{*}_{H}+\nu_{2}I^{*}_{F}
+\nu_{3}I^{*}_{D}\right)M}{M^{*}}\right).
\end{array}
\right\}
\label{eqn21}
\end{equation}
Then, equation $\eqref{eqn21}$ can be expressed as follows: 
\begin{equation}
\left.
\begin{array}{llll}
\frac{d\mathcal{L}}{dt} = -G_{1}\mu_{1}S_{H}\left(1
-\frac{S^{*}_H}{S_{H}} \right)^{2}+G_{1}\tau_{1}S_{H}
I_{F}\left(1-\frac{S^{*}_H}{S_{H}} \right)\left(
\dfrac{I^{*}_{F}S^{*}_{H}}{I_{F}S_{H}}-1\right)
+G_{1}\tau_{2}S_{H}I_{D}\left(1-\frac{S^{*}_H}{S_{H}}\right)
\left(\dfrac{I^{*}_{D}S^{*}_{H}}{I_{D}S_{H}}-1\right)\\
+G_{1}\tau_{3}S_{H}\lambda\left(M\right)\left(1-\frac{S^{*}_H}{S_{H}}\right)
\left(\dfrac{\lambda\left(M^{*}\right)S^{*}_{H}}{\lambda
\left(M\right)S_{H}}-1\right)+G_{1}\beta_{3}R_{H}\left(1-\frac{S^{*}_H}{S_{H}}\right)
\left(1-\frac{R^{*}_H}{R_{H}} \right)\\
+G_{2}\tau_{1}S_{H}I_{F}\left(1-\frac{E^{*}_H}{E_{H}}\right)
\left(1-\dfrac{I^{*}_{F}S^{*}_{H}E_{H}}{I_{F}S_{H}E^{*}_{H}}\right)
+G_{2}\tau_{2}S_{H}I_{D}\left(1-\frac{E^{*}_H}{E_{H}}\right)
\left(1-\dfrac{I^{*}_{D}S^{*}_{H}E_{H}}{I_{D}S_{H}E^{*}_{H}}\right)\\
+G_{2}\tau_{3}S_{H}\lambda\left(M\right)\left(1-\frac{E^{*}_H}{E_{H}}\right)
\left(1-\dfrac{\lambda\left(M^{*}\right)S^{*}_{H}E_{H}}{
\lambda\left(M\right)S_{H}E^{*}_{H}}\right)\\
+G_{3}\beta_{1}E_{H}\left(1-\frac{I^{*}_H}{I_{H}}\right)
\left(1-\dfrac{E^{*}_{H}I_{H}}{E_{H}I^{*}_{H}}\right)
+G_{4}\beta_{2}E_{H}\left(1-\frac{R^{*}_H}{R_{H}}\right)
\left(-\dfrac{E^{*}_{H}R_{H}}{E_{H}R^{*}_{H}}\right)\\
-G_{5}\mu_{2}S_{F}\left(1-\frac{S^{*}_F}{S_{F}} \right)^{2}
+G_{5}\kappa_{1}S_{F}I_{F}\left(1-\frac{S^{*}_F}{S_{F}}\right)
\left(\dfrac{I^{*}_{F}S^{*}_{F}}{I_{F}S_{F}}-1\right)
+G_{5}\kappa_{2}S_{F}I_{D}\left(1-\frac{S^{*}_F}{S_{F}}\right)
\left(\dfrac{I^{*}_{D}S^{*}_{F}}{I_{D}S_{F}}-1\right)\\
+G_{5}\kappa_{3}S_{F}\lambda\left(M\right)\left(1-\frac{S^{*}_{F}}{S_{F}}\right)
\left(\dfrac{\lambda\left(M^{*}\right)S^{*}_{F}}{\lambda\left(M\right)S_{F}}-1\right)
+ G_{6}\kappa_{1}S_{F}I_{F}\left(1-\frac{E^{*}_F}{E_{F}}\right)
\left(1-\dfrac{I^{*}_{F}S^{*}_{F}E_{F}}{I_{F}S_{F}E^{*}_{F}}\right)\\
+G_{6}\kappa_{2}S_{F}I_{D}\left(1-\frac{E^{*}_F}{E_{F}}\right)
\left(1-\dfrac{I^{*}_{D}S^{*}_{F}E_{F}}{I_{D}S_{F}E^{*}_{F}}\right)\\
+G_{6}\kappa_{3}S_{F}\lambda\left(M\right)\left(1-\frac{E^{*}_F}{E_{F}}\right)
\left(1-\dfrac{\lambda\left(M^{*}\right)S^{*}_{F}E_{F}}{\lambda\left(M\right)
S_{F}E^{*}_{F}}\right)+G_{7}\gamma E_{F}\left(1-\frac{I^{*}_F}{I_{F}}\right)
\left(1-\dfrac{E^{*}_{F}I_{F}}{E_{F}I^{*}_{F}}\right)\\
-G_{8}\mu_{3}S_{D}\left(1-\frac{S^{*}_D}{S_{D}} \right)^{2}
+\dfrac{\psi_{1}S_{D}I_{F}G_{8}}{\left(1+\rho_{1}\right)}
\left(1-\frac{S^{*}_D}{S_{D}}\right)\left(\dfrac{I^{*}_{F}
S^{*}_{D}}{I_{F}S_{D}}-1\right)+\dfrac{\psi_{2}S_{D}I_{F}
G_{8}}{\left(1+\rho_{2}\right)}\left(1-\frac{S^{*}_D}{S_{D}}\right)
\left(\dfrac{I^{*}_{D}S^{*}_{D}}{I_{D}S_{D}}-1\right)\\
+\dfrac{\psi_{3}S_{D}\lambda\left(M\right)G_{8}}{\left(1
+\rho_{3}\right)} \left(1-\frac{S^{*}_D}{S_{D}}\right)
\left(\dfrac{\lambda\left(M^{*}\right)S^{*}_{D}}{\lambda
\left(M\right)S_{D}}-1\right)+G_{8}\gamma_{3}R_{D}\left(
1-\frac{S^{*}_D}{S_{D}} \right)\left(1-\frac{R^{*}_D}{R_{D}} \right)\\
+\dfrac{\psi_{1}S_{D}I_{F}G_{9}}{\left(1+\rho_{1}\right)}
\left(1-\frac{E^{*}_D}{E_{D}} \right)\left(1-\dfrac{I^{*}_{F}S^{*}_{D}
E_{D}}{I_{F}S_{D}E^{*}_{D}}\right)+\dfrac{\psi_{2}S_{D}I_{F}
G_{9}}{\left(1+\rho_{2}\right)}\left(1-\frac{E^{*}_D}{E_{D}}\right)
\left(1-\dfrac{I^{*}_{D}S^{*}_{D}E_{D}}{I_{D}S_{D}E^{*}_{D}}\right)\\
+\dfrac{\psi_{3}S_{D}I_{F}G_{9}}{\left(1+\rho_{3}\right)}
\left(1-\frac{E^{*}_D}{E_{D}} \right)\left(1
-\dfrac{\lambda\left(M^{*}\right)S^{*}_{D}E_{D}}{\lambda
\left(M\right)S_{D}E^{*}_{D}}\right)\\
+G_{10}\gamma_{1}E_{D}\left(1-\frac{I^{*}_D}{I_{D}}\right)
\left(1-\dfrac{E^{*}_{D}I_{D}}{E_{D}I^{*}_{D}}\right)
+G_{11}\gamma_{2}E_{D}\left(1-\frac{R^{*}_D}{R_{D}}\right)
\left(-\dfrac{E^{*}_{D}R_{D}}{E_{D}R^{*}_{D}}\right)\\ \\
+G_{12}\nu_{1}I_{H}\left(1-\frac{M^{*}}{M} \right)
\left(1-\frac{I^{*}_{H}M}{I_{H}M^{*}} \right)
+G_{12}\nu_{2}I_{F}\left(1-\frac{M^{*}}{M} \right)
\left(1-\frac{I^{*}_{F}M}{I_{F}M^{*}} \right)
+G_{12}\nu_{3}I_{D}\left(1-\frac{M^{*}}{M} \right)
\left(1-\frac{I^{*}_{D}M}{I_{D}M^{*}} \right).
\end{array}
\right\}
\label{eqn23}
\end{equation}
Equation $\left(\ref{eqn23}\right)$ can be written as 
\begin{equation*}
\begin{array}{llll}
\frac{d\mathcal{L}}{dt}=\mathcal{Q}+\mathcal{R},
\end{array}
\end{equation*}
where  
\begin{equation*}
\begin{array}{llll}
{\cal R}=-G_{1}\mu_{1} S_{H}\left(1-\frac{S^{*}_H}{S_H}\right)^{2}
-G_{5}\mu_{2} S_{F}\left(1-\frac{S^{*}_F}{S_F}\right)^{2}
-G_{8}\mu_{3} S_{D}\left(1-\frac{S^{*}_D}{S_D}\right)^{2}
\end{array}
\end{equation*} 
and  
\begin{equation}
\left.
\begin{array}{llll}
\mathcal{Q} = G_{1}\tau_{1}S_{H}I_{F}\left(1-\frac{S^{*}_H}{S_{H}} 
\right)\left(\dfrac{I^{*}_{F}S^{*}_{H}}{I_{F}S_{H}}-1\right)
+G_{1}\tau_{2}S_{H}I_{D}\left(1-\frac{S^{*}_H}{S_{H}} \right)
\left(\dfrac{I^{*}_{D}S^{*}_{H}}{I_{D}S_{H}}-1\right)\\
+G_{1}\tau_{3}S_{H}\lambda\left(M\right)\left(1-\frac{S^{*}_H}{S_{H}} \right)
\left(\dfrac{\lambda\left(M^{*}\right)S^{*}_{H}}
{\lambda\left(M\right)S_{H}}-1\right)+G_{1}\beta_{3}R_{H}
\left(1-\frac{S^{*}_H}{S_{H}} \right)\left(1-\frac{R^{*}_H}{R_{H}} \right)\\
+G_{2}\tau_{1}S_{H}I_{F}\left(1-\frac{E^{*}_H}{E_{H}} \right)
\left(1-\dfrac{I^{*}_{F}S^{*}_{H}E_{H}}{I_{F}S_{H}E^{*}_{H}}\right)
+G_{2}\tau_{2}S_{H}I_{D}\left(1-\frac{E^{*}_H}{E_{H}} \right)
\left(1-\dfrac{I^{*}_{D}S^{*}_{H}E_{H}}{I_{D}
S_{H}E^{*}_{H}}\right)\\
+G_{2}\tau_{3}S_{H}
\lambda\left(M\right)\left(1-\frac{E^{*}_H}{E_{H}} \right)
\left(1-\dfrac{\lambda\left(M^{*}\right)S^{*}_{H}
E_{H}}{\lambda\left(M\right)S_{H}E^{*}_{H}}\right)\\
+G_{3}\beta_{1}E_{H}\left(1-\frac{I^{*}_H}{I_{H}}\right)
\left(1-\dfrac{E^{*}_{H}I_{H}}{E_{H}I^{*}_{H}}\right)
+G_{4}\beta_{2}E_{H}\left(1-\frac{R^{*}_H}{R_{H}}\right)
\left(-\dfrac{E^{*}_{H}R_{H}}{E_{H}R^{*}_{H}}\right)\\
+G_{5}\kappa_{1}S_{F}I_{F}\left(1-\frac{S^{*}_F}{S_{F}} \right)
\left(\dfrac{I^{*}_{F}S^{*}_{F}}{I_{F}S_{F}}-1\right)
+G_{5}\kappa_{2}S_{F}I_{D}\left(1-\frac{S^{*}_F}{S_{F}} \right)
\left(\dfrac{I^{*}_{D}S^{*}_{F}}{I_{D}S_{F}}-1\right)\\
+G_{5}\kappa_{3}S_{F}\lambda\left(M\right)\left(1-\frac{S^{*}_{F}}{S_{F}}\right)
\left(\dfrac{\lambda\left(M^{*}\right)S^{*}_{F}}{\lambda\left(M\right)S_{F}}
-1\right)+ G_{6}\kappa_{1}S_{F}I_{F}\left(1-\frac{E^{*}_F}{E_{F}} \right)
\left(1-\dfrac{I^{*}_{F}S^{*}_{F}E_{F}}{I_{F}S_{F}E^{*}_{F}}\right)\\
+G_{6}\kappa_{2}S_{F}I_{D}\left(1-\frac{E^{*}_F}{E_{F}} \right)
\left(1-\dfrac{I^{*}_{D}S^{*}_{F}E_{F}}{I_{D}S_{F}E^{*}_{F}}\right)\\
+G_{6}\kappa_{3}S_{F}\lambda\left(M\right)\left(1-\frac{E^{*}_F}{E_{F}}\right)
\left(1-\dfrac{\lambda\left(M^{*}\right)S^{*}_{F}
E_{F}}{\lambda\left(M\right)S_{F}E^{*}_{F}}\right)+G_{7}\gamma
E_{F}\left(1-\frac{I^{*}_F}{I_{F}}\right)
\left(1-\dfrac{E^{*}_{F}I_{F}}{E_{F}I^{*}_{F}}\right)\\
+\dfrac{\psi_{1}S_{D}I_{F}G_{8}}{\left(1+\rho_{1}\right)}
\left(1-\frac{S^{*}_D}{S_{D}}\right)\left(\dfrac{I^{*}_{F}S^{*}_{D}}{I_{F}S_{D}}
-1\right)+\dfrac{\psi_{2}S_{D}I_{F}G_{8}}{\left(1+\rho_{2}\right)}
\left(1-\frac{S^{*}_D}{S_{D}}\right)\left(\dfrac{I^{*}_{D}
S^{*}_{D}}{I_{D}S_{D}}-1\right)\\
+\dfrac{\psi_{3}S_{D}\lambda\left(M\right)G_{8}}{\left(1
+\rho_{3}\right)} \left(1-\frac{S^{*}_D}{S_{D}}\right)
\left(\dfrac{\lambda\left(M^{*}\right)S^{*}_{D}}{\lambda
\left(M\right)S_{D}}-1\right)+G_{8}\gamma_{3}
R_{D}\left(1-\frac{S^{*}_D}{S_{D}} \right)\left(1-\frac{R^{*}_D}{R_{D}} \right)\\
+\dfrac{\psi_{1}S_{D}I_{F}G_{9}}{\left(1+\rho_{1}\right)}\left(1-\frac{E^{*}_D}{E_{D}}\right)
\left(1-\dfrac{I^{*}_{F}S^{*}_{D}E_{D}}{I_{F}S_{D}E^{*}_{D}}\right)
+\dfrac{\psi_{2}S_{D}I_{F}G_{9}}{\left(1+\rho_{2}\right)}\left(1
-\frac{E^{*}_D}{E_{D}} \right)
\left(1-\dfrac{I^{*}_{D}S^{*}_{D}E_{D}}{I_{D}S_{D}E^{*}_{D}}\right)\\
+\dfrac{\psi_{3}S_{D}I_{F}G_{9}}{\left(1+\rho_{3}\right)}
\left(1-\frac{E^{*}_D}{E_{D}} \right)\left(1-\dfrac{\lambda\left(M^{*}\right)S^{*}_{D}
E_{D}}{\lambda\left(M\right)S_{D}E^{*}_{D}}\right)\\
+G_{10}\gamma_{1}E_{D}\left(1-\frac{I^{*}_D}{I_{D}}\right)
\left(1-\dfrac{E^{*}_{D}I_{D}}{E_{D}I^{*}_{D}}\right)
+G_{11}\gamma_{2}E_{D}\left(1-\frac{R^{*}_D}{R_{D}}\right)
\left(-\dfrac{E^{*}_{D}R_{D}}{E_{D}R^{*}_{D}}\right)\\ \\
+G_{12}\nu_{1}I_{H}\left(1-\frac{M^{*}}{M} \right)
\left(1-\frac{I^{*}_{H}M}{I_{H}M^{*}} \right)+G_{12}\nu_{2}I_{F}
\left(1-\frac{M^{*}}{M} \right)\left(1-\frac{I^{*}_{F}M}{I_{F}
M^{*}} \right)+G_{12}\nu_{3}I_{D}\left(1-\frac{M^{*}}{M} \right)
\left(1-\frac{I^{*}_{D}M}{I_{D}M^{*}} \right).
\end{array}
\right\}
\label{eqn33}
\end{equation}
To simplify $\left(\ref{eqn33}\right)$, let
\begin{equation*}
\begin{aligned}
a = \dfrac{S_{H}}{S^{*}_{H}}, \;
b = \dfrac{E_{H}}{E^{*}_{H}}, \;
c = \dfrac{I_{H}}{I^{*}_{H}}, \;
d = \dfrac{R_{H}}{R^{*}_{H}}, \;
e = \dfrac{S_{F}}{S^{*}_{F}}, \;
f = \dfrac{E_{F}}{E^{*}_{F}}, \;
g = \dfrac{I_{F}}{I^{*}_{F}}, \\
h = \dfrac{S_{D}}{S^{*}_{D}}, \;
r = \dfrac{E_{D}}{E^{*}_{D}}, \;
n = \dfrac{I_{D}}{I^{*}_{D}}, \;
m=\dfrac{\lambda\left(M\right)}{\lambda\left(M^{*}\right)},\;
l = \dfrac{R_{D}}{R^{*}_{D}}, \;\text{and}\;
k = \dfrac{M}{M^{*}}.
\end{aligned}
\end{equation*}
One gets from $\left(\ref{eqn33}\right)$ that
\begin{equation}
\left.
\begin{array}{llll}
\mathcal{Q} = \tau_{1}S_{H}I_{F}\left(1-\dfrac{1}{a} \right)
\left(\dfrac{1}{ab}-1\right)+\tau_{2}S_{H}I_{D}
\left(1-\dfrac{1}{a} \right)\left(\dfrac{1}{an}-1\right)
+\tau_{3}S_{H}\lambda\left(M\right)\left(1-\dfrac{1}{a} \right)
\left(\dfrac{1}{am}-1\right)\\+ \beta_{3}R_{H}
\left(1-\dfrac{1}{a} \right)\left(1-\dfrac{1}{d} \right)
+ \tau_{1}S_{H}I_{F}\left(1-\dfrac{1}{b} \right)
\left(1-\dfrac{b}{af}\right)+\tau_{2}S_{H}I_{D}
\left(1-\dfrac{1}{b} \right)\left(1-\dfrac{b}{an}\right)\\
+\tau_{3}S_{H}\lambda\left(M\right)\left(1-\frac{1}{b} \right)
\left(1-\dfrac{b}{am}\right)+\beta_{1}E_{H}\left(1-\dfrac{1}{c}\right)
\left(1-\dfrac{b}{c}\right)+\beta_{2}E_{H}\left(1-\dfrac{1}{d}\right)
\left(1-\dfrac{b}{d}\right)\\
+\kappa_{1}S_{F}I_{F}\left(1-\dfrac{1}{e} \right)
\left(\dfrac{1}{ef}-1\right)+\kappa_{2}S_{F}I_{D}
\left(1-\dfrac{1}{e} \right)\left(\dfrac{1}{ne}-1\right)
+\kappa_{3}S_{F}\lambda\left(M\right)\left(1-\dfrac{1}{e}\right)
\left(\dfrac{1}{me}-1\right)\\+ \kappa_{1}S_{F}I_{F}
\left(1-\dfrac{1}{f} \right)\left(1-\dfrac{f}{en}\right)
+\kappa_{2}S_{F}I_{D}\left(1-\dfrac{1}{f} \right)
\left(1-\dfrac{f}{en}\right)+\kappa_{3}S_{F}\lambda
\left(M\right)\left(1-\dfrac{1}{f} \right)
\left(1-\dfrac{f}{me}\right)\\
+\gamma E_{F}\left(1-\dfrac{1}{f}\right)\left(1-\dfrac{g}{f}\right)
+\dfrac{\psi_{1}S_{D}I_{F}}{\left(1+\rho_{1}\right)}\left(1-\dfrac{1}{h}\right)
\left(\dfrac{1}{hg}-1\right)+\dfrac{\psi_{2}S_{D}I_{F}}{\left(1+\rho_{2}\right)}
\left(1-\dfrac{1}{h} \right)\left(\dfrac{1}{hn}-1\right)\\
+\dfrac{\psi_{3}S_{D}\lambda\left(M\right)}{\left(1+\rho_{3}\right)} 
\left(1-\dfrac{1}{h} \right)\left(\dfrac{1}{mh}-1\right)
+\gamma_{3}R_{D}\left(1-\frac{1}{h} \right)\left(1-\dfrac{1}{l} \right)
\dfrac{\psi_{1}S_{D}I_{F}}{\left(1+\rho_{1}\right)}
\left(1-\dfrac{1}{r} \right)\left(1-\dfrac{r}{hg}\right)\\
+\dfrac{\psi_{2}S_{D}I_{F}}{\left(1+\rho_{2}\right)}
\left(1-\dfrac{1}{r} \right)\left(1-\dfrac{r}{hn}\right)
+\dfrac{\psi_{3}S_{D}I_{F}}{\left(1+\rho_{3}\right)}
\left(1-\dfrac{1}{r} \right)\left(1-\dfrac{r}{mh}\right)
+\gamma_{1}E_{D}\left(1-\dfrac{1}{n}\right)\left(1-\dfrac{n}{r}\right)\\
+\gamma_{2}E_{D}\left(1-\dfrac{1}{l}\right)\left(1-\dfrac{l}{r}\right)
+\nu_{1}I_{H}\left(1-\frac{1}{k} \right)\left(1-\frac{k}{c} \right)
+\nu_{2}I_{F}\left(1-\frac{1}{k} \right)\left(1-\frac{k}{g} \right)
+\nu_{3}I_{D}\left(1-\frac{1}{k} \right)\left(1-\frac{k}{n} \right).
\end{array}
\right\}
\label{eqn35}
\end{equation}
We express the  equation $\left(\ref{eqn35}\right)$ as
\begin{equation}
\left.
\begin{array}{llll}
\mathcal{Q} = \tau_{1}S_{H}I_{F}\left(1-\dfrac{b}{ac}+\dfrac{1}{ac}
-\dfrac{1}{b} \right)+\tau_{2}S_{H}I_{D}\left(\dfrac{1}{an}-1+\dfrac{1}{a^{2}n}
+\dfrac{1}{a}\right)+\tau_{3}S_{H}\lambda\left(M\right)
\left(\dfrac{1}{am}-1+\dfrac{1}{a^{2}m}+\dfrac{1}{a}\right)\\
+ \beta_{3}R_{H}\left(1-\dfrac{1}{d}-\dfrac{1}{a} +\dfrac{1}{ad}\right)
+ \tau_{1}S_{H}I_{F}\left(1-\dfrac{b}{af}-\dfrac{1}{b} 
+\dfrac{1}{af}\right)+\tau_{2}S_{H}I_{D}\left(1-\dfrac{b}{an}-\dfrac{1}{b}
+\dfrac{1}{an}\right)\\+\tau_{3}S_{H}\lambda\left(M\right)
\left(1-\dfrac{b}{af}-\dfrac{1}{b} +\dfrac{1}{af}\right)
+\tau_{2}S_{H}I_{D}\left(1-\dfrac{b}{am}-\dfrac{1}{b} 
+\dfrac{1}{am}\right)+\beta_{1}E_{H}\left(1-\dfrac{b}{c}
-\dfrac{1}{c}+\dfrac{b}{c^{2}}\right)\\+\beta_{2}E_{H}
\left(1-\dfrac{b}{d}-\dfrac{1}{d}+\dfrac{b}{d^{2}}\right)
+\kappa_{1}S_{F}I_{F}\left(\dfrac{1}{ef}-1-\dfrac{1}{e^{2}f}
+\dfrac{1}{e}\right)+\kappa_{2}S_{F}I_{D}\left(\dfrac{1}{en}
-1-\dfrac{1}{e^{2}n}+\dfrac{1}{e}\right)\\
+\kappa_{3}S_{F}\lambda\left(M\right)\left(\dfrac{1}{em}-1
-\dfrac{1}{e^{2}m}+\dfrac{1}{e}\right)
+\kappa_{1}S_{F}I_{F}\left(1-\dfrac{f}{en}-\dfrac{1}{f}
+\dfrac{1}{en}\right)+\kappa_{2}S_{F}I_{D}\left(1-\dfrac{f}{en}
-\dfrac{1}{f}+\dfrac{1}{en}\right)\\+\kappa_{3}S_{F}\lambda
\left(M\right)\left(1-\dfrac{f}{me}-\dfrac{1}{f}+\dfrac{1}{me}\right)
+\gamma E_{F}\left(1-\dfrac{g}{f}-\dfrac{1}{f}+\dfrac{g}{f^{2}}\right)
+\dfrac{\psi_{1}S_{D}I_{F}}{\left(1+\rho_{1}\right)}\left(1-\dfrac{1}{h}
-\dfrac{1}{h^{2}g}+\dfrac{1}{h}\right)\\+\dfrac{\psi_{2}S_{D}I_{F}}{\left(1
+\rho_{2}\right)}\left(1-\dfrac{1}{h}-\dfrac{1}{h^{2}g}+\dfrac{1}{h}\right)
+\dfrac{\psi_{3}S_{D}\lambda\left(M\right)}{\left(1+\rho_{3}\right)}
\left(\dfrac{1}{mh}-1-\dfrac{1}{h^{2}m}+\dfrac{1}{h}\right)
+\gamma_{3}R_{D}\left(1-\dfrac{1}{l}-\dfrac{1}{h}+\dfrac{1}{hl}\right)\\
+\dfrac{\psi_{1}S_{D}I_{F}}{\left(1+\rho_{1}\right)}
\left(1-\dfrac{r}{hg}-\dfrac{1}{r}+\dfrac{1}{hg}\right)
+\dfrac{\psi_{2}S_{D}I_{F}}{\left(1+\rho_{2}\right)}
\left(1-\dfrac{r}{hn}-\dfrac{1}{r}+\dfrac{1}{hn}\right)
+\dfrac{\psi_{3}S_{D}I_{F}}{\left(1+\rho_{3}\right)}
\left(1-\dfrac{r}{hm}-\dfrac{1}{r}+\dfrac{1}{hm}\right)\\
+\gamma_{1}E_{D}\left(1-\dfrac{l}{r}-\dfrac{1}{l}+\dfrac{1}{r}\right)
+\gamma_{2}E_{D}\left(1-\dfrac{r}{hg}-\dfrac{1}{r}+\dfrac{1}{hg}\right)
+\nu_{1}I_{H}\left(1-\dfrac{k}{c}-\dfrac{1}{k}+\dfrac{1}{c}\right)
+\nu_{2}I_{F}\left(1-\dfrac{k}{g}-\dfrac{1}{k}+\dfrac{1}{g}\right)\\
+\nu_{3}I_{D}\left(1-\dfrac{k}{n}-\dfrac{1}{k}+\dfrac{1}{n}\right). 
\end{array}
\right\}
\label{eqn36}
\end{equation}

Now we make use of the following basic inequality.

\begin{proposition}
If $\epsilon(y) = 1 - y + \ln y $, then $\epsilon(y) \leq 0$  
such  that  $1 - y \leq -\ln y$ if and only if $y > 0$.
\label{def1}
\end{proposition}

From equation $\left(\ref{eqn36}\right)$, we have
\begin{equation}
\begin{aligned}
1- \dfrac{1}{d}-\dfrac{1}{a}+\dfrac{1}{ad}
=\left(1-\dfrac{1}{d}\right)+\left(1-\dfrac{1}{a}\right)
-\left(1-\dfrac{1}{ad}\right).
\end{aligned}
\label{def2}
\end{equation}
Using Proposition~\ref{def1} and the concept of geometric mean, 
equation $\left(\ref{def2}\right)$ can be written as
\begin{equation}
\begin{aligned}
\left(1-\dfrac{1}{d}\right)+\left(1-\dfrac{1}{a}\right)
-\left(1-\dfrac{1}{ad}\right)
\le -\ln\left(\dfrac{1}{d}\right)
-\ln\left(\dfrac{1}{a}\right)+\ln\left(\dfrac{1}{ad}\right)\\
\le \ln\left(a\times d\times\dfrac{1}{ad} \right)=\ln\left(1\right)=0.
\end{aligned}
\label{def22}
\end{equation}
Following similar procedures in $\left(\ref{def22}\right)$, we get
\begin{equation*}
\begin{aligned}
1-\dfrac{c}{b}-\dfrac{1}{c}+\dfrac{1}{b}\le 0,
\quad 1-\dfrac{p}{m} -\dfrac{1}{p}+\dfrac{1}{m}\le 0,
\quad 1-\dfrac{d}{b}-\dfrac{1}{d}+\dfrac{1}{b} \le 0.
\end{aligned}
\end{equation*}
From  equation $\left(\ref{eqn23}\right)$, the global stability holds only if 
$\dfrac{d\mathcal{L}}{dt}\leq0$. Now, if $\mathcal{R}<\mathcal{Q}$, then 
$\dfrac{d\mathcal{L}}{dt}$ will be negative definite, which implies that 
$\dfrac{d\mathcal{L}}{dt}<0$ and $\dfrac{d\mathcal{L}}{dt}=0$ only at the 
endemic equilibrium point $\mathbb{E}^{*}$. Hence, by LaSalle's 
invariance principle \cite{lasalle1976stability}, the only invariant set in 
$\left\{\left(S_{H}(t),E_{H}(t),I_{H}(t),R_{H}(t),S_{F}(t),E_{F}(t),I_{F}(t),
S_{D}(t),E_{D}(t),I_{D}(t),R_{H}(t)\right)\in\mathbb{R}_{+}^{12}\right\}:$
$\left\{\left(S_{H}(t),E_{H}(t),I_{H}(t),R_{H}(t),S_{F}(t),E_{F}(t),
I_{F}(t),S_{D}(t),E_{D}(t),I_{D}(t),R_{H}(t)\right)
\rightarrow {\mathbb E}^{*}\right\}$ 
is the singleton endemic point ${\mathbb E}^{*}$. Thus, any solution 
to the rabies model $\left(\ref{eqn1}\right)$ which intersect the interior 
$\mathbb{R}_{+}^{12}$ limits to ${\mathbb E}^{*}$ is globally asymptotically 
stable whatever ${\cal R}_0>1$.
\end{proof}
 



\begin{thebibliography}{10}
	
\bibitem{de2022importance}
Paola De~Benedictis, Stefania Leopardi, Wanda Markotter, and Andres
Velasco-Villa.
\newblock The importance of accurate host species identification in the
framework of rabies surveillance, control and elimination.
\newblock {\em Viruses}, 14(3):492, 2022.
	
\bibitem{kumar2023canine}
Anil Kumar, Sonam Bhatt, Ankesh Kumar, and Tanmoy Rana.
\newblock Canine rabies: An epidemiological significance, pathogenesis,
diagnosis, prevention and public health issues.
\newblock {\em Comparative Immunology, Microbiology and Infectious Diseases},
page 101992, 2023.
	
\bibitem{slathia2023rabies}
Pallvi Slathia, Riya Abrol, Satuti Sharma, and Sakshi Sharma.
\newblock Rabies: A review on clinical signs, prevention and control.
\newblock {\em The Pharma Innovation Journal}, 12(5):1675--1680, 2023.
	
\bibitem{mcmahon2018ecosystem}
Barry~J McMahon, Serge Morand, and Jeremy~S Gray.
\newblock Ecosystem change and zoonoses in the anthropocene.
\newblock {\em Zoonoses and Public Health}, 65(7):755--765, 2018.
	
\bibitem{nigg2009overview}
Andrea~Julia Nigg and Pamela~L Walker.
\newblock Overview, prevention, and treatment of rabies.
\newblock {\em Pharmacotherapy: The Journal of Human Pharmacology and Drug Therapy}, 
29(10):1182--1195, 2009.
	
\bibitem{johnson2010immune}
Nicholas Johnson, Adam~F Cunningham, and Anthony~R Fooks.
\newblock The immune response to rabies virus infection and vaccination.
\newblock {\em Vaccine}, 28(23):3896--3901, 2010.
	
\bibitem{hailemichael2022effect}
Demsis~Dejene Hailemichael, Geremew~Kenassa Edessa, and Purnachandra~Rao Koya.
\newblock Effect of vaccination and culling on the dynamics of rabies
transmission from stray dogs to domestic dogs.
\newblock {\em Journal of Applied Mathematics}, 2022, 2022.
	
\bibitem{bilal2021rabies}
A~Bilal.
\newblock Rabies is a zoonotic disease: a literature review.
\newblock {\em Occup. Med. Health Aff}, 9(2), 2021.
	
\bibitem{mbilo2021dog}
C{\'e}line Mbilo, Andre Coetzer, Bassirou Bonfoh, Ang{\'e}lique Angot, Charles
Bebay, Bernardo Cassam{\'a}, Paola De~Benedictis, Moina~Hasni Ebou, Corneille
Gnanvi, Vessaly Kallo, et~al.
\newblock Dog rabies control in west and central africa: A review.
\newblock {\em Acta tropica}, 224:105459, 2021.
	
\bibitem{tian2018transmission}
Huaiyu Tian, Yun Feng, Bram Vrancken, Bernard Cazelles, Hua Tan, Mandev~S Gill,
Qiqi Yang, Yidan Li, Weihong Yang, Yuzhen Zhang, et~al.
\newblock Transmission dynamics of re-emerging rabies in domestic dogs of rural China.
\newblock {\em PLoS Pathogens}, 14(12):e1007392, 2018.
	
\bibitem{sambo2013burden}
Maganga Sambo, Sarah Cleaveland, Heather Ferguson, Tiziana Lembo, Cleophas
Simon, Honorati Urassa, and Katie Hampson.
\newblock The burden of rabies in tanzania and its impact on local communities.
\newblock {\em PLoS neglected tropical diseases}, 7(11):e2510, 2013.
	
\bibitem{abdulmajid2021analysis}
Shafiu Abdulmajid and Adamu~Shitu Hassan.
\newblock Analysis of time delayed rabies model in human and dog populations with controls.
\newblock {\em Afrika Matematika}, 32(5-6):1067--1085, 2021.
	
\bibitem{amoako2021rabies}
YA~Amoako, P~El-Duah, AA~Sylverken, M~Owusu, R~Yeboah, R~Gorman, T~Adade,
J~Bonney, W~Tasiame, K~Nyarko-Jectey, et~al.
\newblock Rabies is still a fatal but neglected disease: a case report.
\newblock {\em Journal of Medical Case Reports}, 15(1):1--6, 2021.
	
\bibitem{abrahamian2022rhabdovirus}
Fredrick~M Abrahamian and Charles~E Rupprecht.
\newblock Rhabdovirus: Rabies.
\newblock In {\em Viral Infections of Humans: 
Epidemiology and Control}, 1--49, Springer, 2022.
	
\bibitem{ruan2017modeling}
Shigui Ruan.
\newblock Modeling the transmission dynamics and control of rabies in china.
\newblock {\em Mathematical biosciences}, 286:65--93, 2017.
	
\bibitem{ayoade2019saturated}
Abayomi~Ayotunde Ayoade, Olumuyiwa~James Peter, Tokunbo~Aderemi Ayoola, and
AA~Victor.
\newblock A saturated treatment model for the transmission dynamics of
rabies/ayoade abayomi ayotunde…[et al.].
\newblock {\em Malaysian Journal of Computing (MJoC)}, 4(1):201--213, 2019.
	
\bibitem{chapwanya2022environment}
Michael Chapwanya and Phindile Dumani.
\newblock Environment considerations on the spread of rabies among 
African wild dogs (lycaon pictus) with control measures.
\newblock {\em Mathematical Methods in the Applied Sciences}, 45(8):4124--4139, 2022.
	
\bibitem{kadowaki2018risk}
H~Kadowaki, K~Hampson, K~Tojinbara, A~Yamada, and K~Makita.
\newblock The risk of rabies spread in Japan: a mathematical modelling assessment.
\newblock {\em Epidemiology \& Infection}, 146(10):1245--1252, 2018.
	
\bibitem{ayoade2023modeling}
Abayomi~Ayotunde Ayoade and Mohammed~Olanrewaju Ibrahim.
\newblock Modeling the dynamics and control of rabies in dog population within
and around {L}agos, {N}igeria.
\newblock {\em The European Physical Journal Plus}, 138(5):397, 2023.
	
\bibitem{tulu2017impact}
Aberu~Mengistu Tulu and Purnachandra~Rao Koya.
\newblock The impact of infective immigrants on the spread of dog rabies.
\newblock {\em American Journal of Applied Mathematics}, 5(3):68, 2017.
	
\bibitem{pantha2021modeling}
Buddhi Pantha, Sunil Giri, Hem~Raj Joshi, and Naveen~K Vaidya.
\newblock Modeling transmission dynamics of rabies in Nepal.
\newblock {\em Infectious Disease Modelling}, 6:284--301, 2021.
	
\bibitem{ega2015sensitivity}
Tesfaye~Tadesse Ega, Livingstone Luboobi, Dmitry Kuznetsov, and Abraham~Haile Kidane.
\newblock Sensitivity analysis and numerical simulations for the mathematical
model of rabies in human and animal within and around {A}ddis {A}baba.
\newblock {\em Asian Journal of Mathematics and Applications}, 2015:Art. ID ama0271, 23~pp, 2015.
	
\bibitem{esposito2023impact}
Michelle~Marie Esposito, Sara Turku, Leora Lehrfield, and Ayat Shoman.
\newblock The impact of human activities on zoonotic infection transmissions.
\newblock {\em Animals}, 13(10):1646, 2023.
	
\bibitem{hampson2019potential}
Katie Hampson, Francesco Ventura, Rachel Steenson, Rebecca Mancy, Caroline
Trotter, Laura Cooper, Bernadette Abela-Ridder, Lea Knopf, Moniek Ringenier,
Tenzin Tenzin, et~al.
\newblock The potential effect of improved provision of rabies post-exposure
prophylaxis in gavi-eligible countries: a modelling study.
\newblock {\em The Lancet Infectious Diseases}, 19(1):102--111, 2019.
	
\bibitem{rulli2021land}
Maria~Cristina Rulli, Paolo D’Odorico, Nikolas Galli, and David~TS Hayman.
\newblock Land-use change and the livestock revolution increase the risk of
zoonotic coronavirus transmission from rhinolophid bats.
\newblock {\em Nature Food}, 2(6):409--416, 2021.
	
\bibitem{peng2022transmission}
Hao Peng, Qing Yang, Xinhong Zhang, and Daqing Jiang.
\newblock Transmission dynamics of a high dimensional rabies epidemic model in
a Markovian random environment.
\newblock {\em Qualitative Theory of Dynamical Systems}, 21(2):46, 2022.
	
\bibitem{allen2008basic}
Linda~JS Allen and Pauline van~den Driessche.
\newblock The basic reproduction number in some discrete-time epidemic models.
\newblock {\em Journal of difference equations and applications},
14(10-11):1127--1147, 2008.
	
\bibitem{lasalle1976stability}
Joseph~P LaSalle.
\newblock Stability theory and invariance principles.
\newblock In {\em Dynamical systems}, pages 211--222. Elsevier, 1976.
	
\bibitem{yang2014basic}
Hyun~Mo Yang.
\newblock The basic reproduction number obtained from Jacobian and next
generation matrices--a case study of dengue transmission modelling.
\newblock {\em Biosystems}, 126:52--75, 2014.
	
\bibitem{saha2021dynamics}
Sangeeta Saha and Guruprasad Samanta.
\newblock Dynamics of an epidemic model under the influence of environmental
stress.
\newblock {\em Mathematical Biology and Bioinformatics}, 16(2):201--243, 2021.
	
\bibitem{dharmaratne2020estimation}
Samath Dharmaratne, Supun Sudaraka, Ishanya Abeyagunawardena, Kasun
Manchanayake, Mahen Kothalawala, and Wasantha Gunathunga.
\newblock Estimation of the basic reproduction number ($R_0$) for the novel
coronavirus disease in Sri Lanka.
\newblock {\em Virology Journal}, 17:1--7, 2020.
	
\bibitem{zhang2011analysis}
Juan Zhang, Zhen Jin, Gui-Quan Sun, Tao Zhou, and Shigui Ruan.
\newblock Analysis of rabies in China: transmission dynamics and control.
\newblock {\em PLoS one}, 6(7):e20891, 2011.
	
\bibitem{asamoah2017modelling}
Joshua Kiddy~K Asamoah, Francis~T Oduro, Ebenezer Bonyah, Baba Seidu, et~al.
\newblock Modelling of rabies transmission dynamics using optimal control analysis.
\newblock {\em Journal of Applied Mathematics}, 2017, 2017.
	
\bibitem{world2010working}
World~Health Organization et~al.
\newblock {\em Working to overcome the global impact of neglected tropical
diseases: first {WHO} report on neglected tropical diseases}.
\newblock World Health Organization, 2010.
	
\bibitem{world2013expert}
World~Health Organization.
\newblock {\em WHO expert consultation on rabies: second report}, volume 982.
\newblock World Health Organization, 2013.
	
\bibitem{addo2012seir}
Kwaku~Mari Addo.
\newblock {\em An {SEIR} Mathematical Model for Dog Rabies; Case Study: {B}ongo
District, {G}hana}.
\newblock PhD thesis, Kwame Nkrumah University of Science and Technology, 2012.
	
\bibitem{ruan2017spatiotemporal}
Shigui Ruan.
\newblock Spatiotemporal epidemic models for rabies among animals.
\newblock {\em Infectious disease modelling}, 2(3):277--287, 2017.
	
\bibitem{kamgang2008computation}
Jean~Claude Kamgang and Gauthier Sallet.
\newblock Computation of threshold conditions for epidemiological models and
global stability of the disease-free equilibrium (DFE).
\newblock {\em Mathematical Biosciences}, 213(1):1--12, 2008.
	
\bibitem{li2018introduction}
Michael~Y Li.
\newblock {\em An introduction to mathematical modeling of infectious diseases}, 
\newblock Springer, 2018.
	
\bibitem{myung2003tutorial}
In~Jae Myung.
\newblock Tutorial on maximum likelihood estimation.
\newblock {\em Journal of Mathematical Psychology}, 47(1):90--100, 2003.
	
\bibitem{abdulmoghni2021incidence}
Rihana~Taher Abdulmoghni, Ahmed~Hasan Al-Ward, Khaled~Abdullah Al-Moayed,
Mohammed~Abdullah Al-Amad, and Yousef~S Khader.
\newblock Incidence, trend, and mortality of human exposure to rabies in yemen,
2011-2017: observational study.
\newblock {\em JMIR Public Health and Surveillance}, 7(6):e27623, 2021.
	
\bibitem{dutta2024periodic}
Protyusha Dutta, Guruprasad Samanta, and Juan~J Nieto.
\newblock Periodic transmission and vaccination effects in epidemic dynamics: a
study using the sivis model.
\newblock {\em Nonlinear Dynamics}, pages 1--29, 2024.
	
\bibitem{osman2022analysis}
Shaibu Osman, Getachew~Teshome Tilahun, Seleshi~Demie Alemu, and Winnie~Mokeira Onsongo.
\newblock Analysis of the dynamics of rabies in north shewa, ethiopia.
\newblock {\em Italian J. Pure Appl. Math}, 48:877--902, 2022.
	
\bibitem{musaili2020mathematical}
Jane~S Musaili and Isaac Chepkwony.
\newblock A mathematical model of rabies transmission dynamics in dogs
incorporating public health education as a control strategy-a case study of Makueni county.
\newblock {\em Journal of Advances in Mathematics and Computer Science},
35(1):1--11, 2020.
	
\end{thebibliography}
\end{document}